\magnification1050

  
  %
  %

  %
  %
  \let\footnotea=\footnote
  \def\anote#1#2{\footnotea{\hbox{$^{#1}$}}{\eightpoint#2}}  
  \catcode`@=12 

 \def\defrefnote#1{\definexref{#1}{{\the\footnotenumber}}{refnotes}}

  %
  %


 \input eplain.tex
\makeatletter
\def\numberedfootnote{%
ÊÊ\global\advance\footnotenumber by 1
ÊÊ\@eplainfootnote{{\number\footnotenumber}}%
}%
\def\makecolumns#1/#2 {\par \begingroup
ÊÊ \@columndepth = #1
ÊÊ \advance\@columndepth by -1
ÊÊ \divide \@columndepth by #2
ÊÊ \advance\@columndepth by 1
ÊÊ \@linestogoincolumn = \@columndepth
ÊÊ \@linestogo = #1
ÊÊ \currentcolumn = 1
ÊÊ \def\@endcolumnactions{%
ÊÊÊÊÊÊ\ifnum \@linestogo<2
ÊÊÊÊÊÊÊÊ \the\crtok \egroup \endgroup \par 
ÊÊÊÊÊÊ\else
ÊÊÊÊÊÊÊÊ \global\advance\@linestogo by -1
ÊÊÊÊÊÊÊÊ \ifnum\@linestogoincolumn<2
ÊÊÊÊÊÊÊÊÊÊÊÊ\global\advance\currentcolumn by 1
ÊÊÊÊÊÊÊÊÊÊÊÊ\global\@linestogoincolumn = \@columndepth
ÊÊÊÊÊÊÊÊÊÊÊÊ\the\crtok
ÊÊÊÊÊÊÊÊ \else
ÊÊÊÊÊÊÊÊÊÊÊÊ&\global\advance\@linestogoincolumn by -1
ÊÊÊÊÊÊÊÊ \fi
ÊÊÊÊÊÊ\fi
ÊÊ }%
ÊÊ \makeactive\^^M
ÊÊ \letreturn \@endcolumnactions
ÊÊ \@columnwidth = \hsize
ÊÊÊÊ \advance\@columnwidth by -\parindent
ÊÊÊÊ \divide\@columnwidth by #2
ÊÊ \penalty\abovecolumnspenalty
ÊÊ \noindent 
ÊÊ \valign\bgroup
ÊÊÊÊ &\hbox to \@columnwidth{\strut \hsize = \@columnwidth ##\hfil}\cr
}%
\makeatother

\lefteqnumbers
   \def\testd{oui}
   \def\choixlat{\ifx\numadroite\testd\righteqnumbers
            \else  \lefteqnumbers\fi}
    \choixlat

\catcode`@=\letter
\def\@eplainfootnote#1{\let\@sf\empty 
  \ifhmode\edef\@sf{\spacefactor\the\spacefactor}\/\fi
  \global\advance\hlfootlabelnumber by 1
  \hlstart@impl{foot}{\hlfootlabel}%
  \hldest@impl{footback}{\hlfootbacklabel}%
  \hbox{$^{(#1)}$}%
  \hlend@impl{foot}%
  \@sf\vfootnote{#1.}%
}%
\catcode`@=\other

  \interfootnoteskip=0pt
  \let\note=\numberedfootnote
  \everyfootnote={\eightpoint\leftskip=5truemm\rightskip5truemm}
  
  \hsize150truemm\vsize 240truemm\hoffset=5truemm
  \def\dimstand{\hsize 150truemm\vsize 240truemm\hoffset=5truemm\voffset=0truemm}
  
  \pretolerance=500\tolerance=1000\brokenpenalty=5000
  \parindent3mm
  
  \countdef\temps=170
  \temps=\time
  \countdef\nminutes=171{\nminutes = \time}
  \countdef\nheure=172
  \def\heure{\begingroup                     
     \temps = \time \divide\temps by 60
     \nheure = \temps                        
     \nminutes = \time
     \multiply\temps by 60
     \advance\nminutes by -\temps            
     \ifnum\nminutes<10 \toks1 = {0}%
     \else\toks1 = {}%
     \fi
     \number\nheure h\the\toks1 \number\nminutes  
  \endgroup}%

  \newcount\chstart
  \chstart=\pageno
 \headline={\ifnum\pageno=\chstart {\hfill} \else {\hss \tenrm --\ \folio\ --\hss}\fi}
  \footline={\hfill}
  \normalbaselines
  \frenchspacing
    \def\dater{\vglue-10mm\rightline{(\the\day/\the\month/\the\year)}}
  \def\dateheure{\vglue-10mm\rightline{(\the\day/\the\month/\the\year,\ \heure)}}

  \newif\ifpagetitre \pagetitretrue
\newtoks\hautpagetitre \hautpagetitre={\hfill}
\newtoks\baspagetitre \baspagetitre={\hfill}
\newtoks\auteurcourant \auteurcourant={\hfill}
\newtoks\titrecourant \titrecourant={\hfill}
\newtoks\hautpagegauche
\newtoks\hautpagedroite
\newtoks\hautpagemilieu
\hautpagemilieu={\tenrm\hfil -- \folio\ -- \hfil}
\hautpagegauche={\ifx\midfolio\oui\the\hautpagemilieu\else\tenrm\folio\hfill\the\auteurcourant\hfill\fi}
\hautpagedroite={\ifx\midfolio\oui\the\hautpagemilieu\else\hfill\the\titrecourant\hfill\tenrm\folio\fi}
\newtoks\baspagegauche \baspagegauche={\hfil}
\newtoks\baspagedroite \baspagedroite={\hfil}
\headline={\ifpagetitre\the\hautpagetitre
\else\ifodd\pageno\the\hautpagedroite\else\the\hautpagegauche\fi\fi }
\footline={\ifpagetitre\the\baspagetitre
\else\ifodd\pageno\the\baspagedroite
\else\the\baspagegauche\fi\fi \global\pagetitrefalse}

\def\pageblanche{\vfill\eject\pagetitretrue
\null\vfill\eject
\pagetitretrue
}
\def\chgtpage{\ifodd\pageno \else
\pageblanche \fi \pagetitretrue\titreun=0\footnotenumber=0}

\def\chgtpageincrtitreun{\ifodd\pageno \else
\pageblanche \fi \pagetitretrue\footnotenumber=0}

\def\majnombres{\ifodd\pageno \else
\pageblanche \fi \pagetitretrue\hautpoly\titreun=0\footnotenumber=0}

\def\hautspages#1#2{\auteurcourant={\ninepcap#1}\titrecourant={\nineit#2}}

\ifnum\chstart=\pageno \pagetitretrue\fi
  


  \def\leftnote#1{\vadjust{\setbox1=\vtop{\hsize 20mm\parindent=0pt\eightpoint
  \baselineskip=9pt\rightskip=4mm plus 4mm\vskip-4.7mm#1}\hbox{\kern-2cm\smash{\box1}}}}

  
  \def\raggedcenter{\leftskip=20pt plus 10em  
       \rightskip=\leftskip 
        \parfillskip=0pt 
         \spaceskip=.3333em \xspaceskip=.5em 
          \pretolerance=9999 \tolerance=9999
           \hyphenpenalty=9999 \exhyphenpenalty=9999 }
           
  \def\titrecentre#1{{\parindent0mm\raggedcenter
       \spaceskip=.6em plus .2em minus .2em\xspaceskip=.6em plus .2em minus .2em
        \tit#1\par}}
        


  \def\oui{oui}
  
   \def\fontetitreun{\ifx\paradouze\oui\douzepts\gpdouze\twelvebf\textfont1=\twelveib\else
\quatorzepts\gpquatorze\fourteenbf\fi}

\def\fontetitreunl{\douzepts\textfont1=\twelveib\scriptfont1=\tenib\fourteenti}
 
 \def\fontetitredeux{\textfont1=\eleveni\ifx\paradouze\oui\onzepts\scriptfont1=\ninei\elevenit\else
                        \douzepts\twelveit\fi}
 
   \def\fontetitredeuxb{\ifx\paradouze\oui\onzepts\eleventi\gponze\textfont1=\elevenib\scriptfont1=\nineib
                         \else\douzepts\twelveti\scriptfont1=\twelveib\scriptfont1=\tenib\gpdouze\fi}
                         
\def\fontetitredeuxl{\onzepts\textfont1=\elevenbf\scriptfont1=\ninebf\twelvebf}
  
\def\fontetitretrois{\textfont0=\elevenrm\scriptfont0=\eightrm\textfont1=\eleveni
                      \scriptfont1=\eighti\scriptscriptfont1=\sixi\elevenit}
                      
\def\fontetitrequatre{\textfont0=\elevenrm\scriptfont0=\eightrm\textfont1=\eleveni
                      \scriptfont1=\eighti\scriptscriptfont1=\sixi\elevenrm}
  
  \newcount\titreun\titreun=0
  \newcount\titredeux\titredeux=0
  \newcount\titretrois\titretrois=0
  \newcount\titrequatre\titrequatre=0
  \newcount\enonce\enonce=0
  
  \def\incr#1{\global\advance#1 by 1 {\the #1}}
  \def\avance#1{\global\advance#1 by 1}
  \def\init#1{\global#1=0}
  
  \long\def\Indentation#1#2{\setbox10=\hbox{\fontetitreun#1}
  	                    \ifdim\wd10 < 4mm
                         \setbox10=\hbox to 4mm{\box10\hfill}
                       \else\ifdim\wd10 < 6mm
                         \setbox10=\hbox to 6mm{\box10\hfill}
  	                    \else\ifdim\wd10 < 8mm
                         \setbox10=\hbox to 8mm{\box10\hfill}
                       \else\ifdim\wd10 < 12mm
                         \setbox10=\hbox to 12mm{\box10\hfill}
                       \fi\fi\fi\fi
                       \dimen10=\hsize
                       \advance \dimen10 by -\wd10
                       \noindent \box10 %
                       \ignorespaces
                       \hbox{\vtop{\hsize=\dimen10\raggedright\noindent\fontetitreun#2}}}

  \long\def\paraun#1{\removelastskip\par\medskip\goodbreak\vskip0pt plus.01\vsize\penalty-100
                \vskip0pt plus-.01\vsize
  	              \init{\titredeux}\ifnum\optionparag=1{\init\eqnumber\init\enonce}\else{}\fi
                  \goodbreak{\fontetitreun
  	                \Indentation{\incr{\titreun}.\ }{\fontetitreun #1\par}}\nobreak\medskip}

 %
 %
 \long\def\paraunc#1{\removelastskip\par\bigskip\goodbreak\vskip0pt plus.01\vsize\penalty-100
                \vskip0pt plus-.01\vsize
  	              \init{\titredeux}
                 \ifnum\optionparag=1{\init{\eqnumber}\init\enonce}\else{}\fi
                  \goodbreak
  	                {\parindent0mm\raggedcenter\fontetitreun\incr{\titreun}.\ 
                     \fontetitreun #1\par}\nobreak\medskip}
                     
\newtoks\titreunl
\titreunl={\ifnum\titreun=1{I}\fi%
\ifnum\titreun=2{II}\fi%
\ifnum\titreun=3{III}\fi%
\ifnum\titreun=4{IV}\fi%
\ifnum\titreun=5{V}\fi%
\ifnum\titreun=6{VI}\fi%
\ifnum\titreun=7{VII}\fi%
\ifnum\titreun=8{VIII}\fi%
\ifnum\titreun=9{IX}\fi%
\ifnum\titreun=10{X}\fi%
\ifnum\titreun=11{XI}\fi%
\ifnum\titreun=12{XII}\fi%
\ifnum\titreun=13{XIII}\fi%
}
\long\def\paraunl#1{\removelastskip\par\bigskip\bigskip\goodbreak\vskip0pt plus.01\vsize\penalty-100
                \vskip0pt plus-.01\vsize
  	              \init{\titredeux}\ifnum\optionparag=1{\init\eqnumber\init\enonce}\else{}\fi
                  \goodbreak{\fontetitreunl
  	                \Indentation{\global\advance\titreun by 1{\the\titreunl}.\ }{\fontetitreunl #1\par}}\nobreak\smallskip}

  
  \long\def\paradeux#1{\init{\titretrois}\vskip0pt plus.01\vsize\penalty-10
                \vskip0pt plus-.01\vsize\ifx \elie\oui\medskip\ifnum\titredeux>0\medskip\fi\fi
                 \Indentation{\fontetitredeux\the\titreun${\cdot}$\incr{\titredeux}.}
                              {\fontetitredeux\textfont1=\eleveni#1}\nobreak\par }
  
  \long\def\paradeuxb#1{\init{\titretrois}\vskip0pt plus.001\vsize\penalty-10
                \vskip0pt plus-.01\vsize{\ifx \elie\oui\medskip\ifnum\titredeux>0\medskip\fi\fi
                  \Indentation
  {\fontetitredeuxb\the\titreun${\cdot}$\incr{\titredeux}.}{ \fontetitredeuxb#1}}\nobreak
\smallskip}

\newtoks\titredeuxl
\titredeuxl={\ifnum\titredeux=1{A}\fi%
\ifnum\titredeux=2{B}\fi%
\ifnum\titredeux=3{C}\fi%
\ifnum\titredeux=4{D}\fi%
\ifnum\titredeux=5{E}\fi%
\ifnum\titredeux=6{F}\fi%
\ifnum\titredeux=7{G}\fi%
\ifnum\titredeux=8{H}\fi%
\ifnum\titredeux=9{I}\fi%
\ifnum\titredeux=10{J}\fi%
\ifnum\titredeux=11{K}\fi%
\ifnum\titredeux=12{L}\fi%
\ifnum\titredeux=13{M}\fi%
}
 \long\def\paradeuxl#1{\init{\titretrois}\vskip0pt plus.001\vsize\penalty-10
                \vskip0pt plus-.01
                \vsize \bigskip%
                  \Indentation
     {\fontetitredeuxl\global\advance\titredeux by 1
  \quad \the\titreunl${\cdot}$\the\titredeuxl.}{ \fontetitredeuxl#1}
  \removelastskip\nobreak\smallskip}
  

  \long\def\paratrois#1{\init{\titrequatre}\ifdim\lastskip<\smallskipamount
                \removelastskip\smallskip\fi
                 \vskip0pt plus.01\vsize\penalty-10
                  \vskip0pt
plus-.01\vsize{\ifx \elie\oui\ifnum\titretrois>0\medskip\fi\fi
\Indentation{\fontetitretrois\the\titreun${\cdot}$\the\titredeux${\cdot}$\incr{\titretrois}.\ }
  {\hskip0mm\baselineskip=14pt\fontetitretrois#1}\nobreak\smallskip}}
  
  
  \long\def\paratroisl#1{\init{\titrequatre}\ifdim\lastskip<\smallskipamount
                \removelastskip\fi
                 \vskip0pt plus.01\vsize\penalty-10
                  \vskip0pt
plus-.01\vsize\ifx \elie\oui\bigskip
\fi
\Indentation{\fontetitretrois\quad \quad \the\titreunl{${\cdot}$}\the\titredeuxl${\cdot}$\incr{\titretrois}.\ }
  {\hskip0mm\fontetitretrois#1}\nobreak\smallskip}


  \long\def\paraquatre#1{\ifdim\lastskip<\smallskipamount
                \removelastskip\smallskip\fi
                 \vskip0pt plus.01\vsize\penalty-10
                  \vskip0pt
                  plus-.01\vsize\par
 
\Indentation{\fontetitrequatre \the\titreun{${\cdot}$}\the\titredeux${\cdot}$\the\titretrois${\cdot}$\incr{\titrequatre}.\ }
{\hskip0mm\fontetitrequatre#1}\nobreak\smallskip}


\newtoks\titrequatrel
\titrequatrel={\ifnum\titrequatre=1{a}\fi%
\ifnum\titrequatre=2{b}\fi%
\ifnum\titrequatre=3{c}\fi%
\ifnum\titrequatre=4{d}\fi%
\ifnum\titrequatre=5{e}\fi%
\ifnum\titrequatre=6{f}\fi%
\ifnum\titrequatre=7{g}\fi%
\ifnum\titrequatre=8{h}\fi%
\ifnum\titrequatre=9{i}\fi%
\ifnum\titrequatre=10{j}\fi%
\ifnum\titrequatre=11{k}\fi%
\ifnum\titrequatre=12{l}\fi%
\ifnum\titrequatre=13{m}\fi%
}
\long\def\paraquatrel#1{\ifdim\lastskip<\smallskipamount
                \removelastskip\smallskip\fi
                 \vskip0pt plus.01\vsize\penalty-10
                  \vskip0pt
                  plus-.01\vsize{\bigskip
\Indentation{\global\advance\titrequatre by 1
\fontetitrequatre\quad \quad \quad \the\titreunl${\cdot}$\the\titredeuxl${\cdot}$\the\titretrois${\cdot}$\the\titrequatrel.\ }
{\hskip0mm\fontetitrequatre#1}\nobreak\smallskip}}

\ifx\optionkeys\oui
\def\drefun#1{\definexref{¤#1}{{\the\titreun}}{}} 
\def\drefdeux#1{\definexref{¤#1}{{\the\titreun}.{\the\titredeux}}{}}
\def\dreftrois#1{\definexref{¤#1}{{\the\titreun}.{\the\titredeux}.{\the\titretrois}}{}}
\else
\def\drefun#1{\definexref{prg#1}{{\the\titreun}}{}} 
\def\drefdeux#1{\definexref{prg#1}{{\the\titreun}.{\the\titredeux}}{}}
\def\dreftrois#1{\definexref{prg#1}{{\the\titreun}.{\the\titredeux}.{\the\titretrois}}{}}
\fi

%


  \long\def\propdeux#1#2#3#4{%
       \avance{\enonce}
       \leavevmode\edef\temp{#2}%
         \ifx\temp\empty 
          \else
           \definexref{#2}{#1~{\the\titreun.\the\enonce}}{enonces}
            \definexref{s#2}{{\the\titreun.\the\enonce}}{enonces}
             \fi
\smallskip
      \noindent{\bf#1\ {\bf\the\titreun.\the\enonce{#3}.}\enspace}{\sl#4\par}%
      \ifdim\lastskip<\medskipamount \removelastskip\penalty55\par \fi
   }

  \long\def\propun#1#2#3#4{%
      \avance{\enonce}
       \leavevmode\edef\temp{#2}%
        \ifx\temp\empty 
          \else
           \definexref{#2}{#1~{\the\enonce}}{enonces}
            \definexref{{s#2}}{{\the\enonce}}{enonces}
             \fi
   \par 
     \noindent{\bf#1\ {\bf\the\enonce{#3}.}\enspace}{\sl#4\par}%
     \ifdim\lastskip<\medskipamount \removelastskip\penalty55\medskip\fi
  }
  
  \long\def\prop#1#2#3#4{\ifnum\optionparag=1
                          \propdeux{#1}{#2}{\textfont1=\elevenib#3}{#4} \else\propun{#1}{#2}{\textfont1=\elevenib#3}{#4}\fi}

  \long\def\propt#1#2#3{\ifx\tpf\oui \prop{Th\'eo\-r\`eme}{#1}{#2}{#3}\par
                       \else\prop{Theorem}{#1}{#2}{#3}\par\fi}
  \long\def\Propt#1#2{\propt{#1}{}{#2}}
  \long\def\propl#1#2#3{\ifx\tpf\oui\prop{Lem\-me}{#1}{#2}{#3}\par
                         \else\prop{Lemma}{#1}{#2}{#3}\par\fi}
  \long\def\Propl#1#2{\propl{#1}{}{#2}}
  \long\def\propc#1#2#3{\ifx\tpf\oui\prop{Corol\-laire}{#1}{#2}{#3}\par
                         \else\prop{Corollary}{#1}{#2}{#3}\par\fi}
  \long\def\Propc#1#2{\propc{#1}{}{#2}}
  \long\def\propp#1#2#3{\prop{Pro\-po\-si\-tion}{#1}{#2}{#3}\par}
  \long\def\Propp#1#2{\propp{#1}{}{#2}} 
  \long\def\propd#1#2#3{\ifx\tpf\oui\prop{D\'efi\-nition}{#1}{#2}{#3}\par
                       \else\prop{Definition}{#1}{#2}{#3}\par\fi} 
  
  \long\def\proptd#1#2#3{\ifx\tpf\oui\prop{Th\'eor\`eme et d\'efi\-nition}{#1}{#2}{#3}\par
                       \else\prop{Theorem and definition}{#1}{#2}{#3}\par\fi}


  
  \newcount\optionparag\optionparag=1
  
  \long\def\section#1#2{\ifnum\optionparag=1 \paraun{#2} 
                        \else\goodbreak{\fontetitreun
  	                \Indentation{#1.\ }{#2}}\nobreak\smallskip\fi}

  \def\eqconstruct#1{\ifnum\optionparag=1{\the\titreun\hbox{$\cdot$}#1}\else{#1}\fi}

  
  
  \def\numref{oui}  
  
  \newcount\mesref\mesref=0 
  \def\defbib#1{\ifx\numref\oui\global\advance\mesref by 1 \definexref{#1}{{\the
                 \mesref}}{}\else\definexref{#1}{#1}{}\fi}
  \def\bibtem#1{\defbib{#1}\item{\citer{#1}}}
  \def\citer#1{[\ref{#1}]}

  
  \font\seventeenmsa=msam10 at 17pt    
  \font\fourteenmsa=msam10 at 14pt
  \font\twelvemsa=msam10 at 12pt
  \font\tenmsa=msam10                 
  \font\ninemsa=msam10 at 9pt 
  \font\eightmsa=msam10 at 8pt 
  \font\sevenmsa=msam7 
  \font\sixmsa=msam10 at 6pt
  \font\fivemsa=msam5
  \newfam\msafam\textfont\msafam=\tenmsa\scriptfont\msafam=\sevenmsa\scriptscriptfont\msafam=\fivemsa
  
  \font\seventeenbb=msbm10 at 17pt     
  \font\fourteenbb=msbm10 at 14pt
  \font\twelvebb=msbm10 at 12pt
  \font\tenbb=msbm10                   
  \font\ninebb=msbm10 at 9pt 
  \font\eightbb=msbm10 at 8pt 
  \font\sevenbb=msbm7 
  \font\sixbb=msbm10 at 6pt
  \font\fivebb=msbm5 
  \newfam\bbfam\textfont\bbfam=\tenbb\scriptfont\bbfam=\sevenbb\scriptscriptfont\bbfam=\fivebb
  \def\bb{\fam\bbfam\tenbb}%

  \font\seventeenscaln=eusm10 at 17pt   
  \font\twelvescaln=eusm10 at 12pt
  \font\tenscaln=eusm10                
  \font\ninescaln=eusm10 scaled 900
  \font\eightscaln=eusm10 scaled 800
  \font\sevenscaln=eusm10 scaled 700
  \font\sixscaln=eusm10 scaled 600
   
  \newfam\scalnfam\textfont\scalnfam=\tenscaln\scriptfont\scalnfam=\sevenscaln\scriptscriptfont\scalnfam=\sixscaln
  \def\scaln{\fam\scalnfam\tenscaln}%
  \def\scal{\scaln}
  
  \font\tenscalb=eusb10                

  \font\sevenscalb=eusb10 scaled 700

  \newfam\scalbfam\textfont\scalbfam=\tenscalb\scriptfont\scalbfam=\sevenscalb
  %
  
  %
  %
  \font\fourteenrm=cmr12 scaled 1200
  \font\elevenrm=cmr10 at 11pt
  \font\twelverm=cmr12
  \font\ninerm=cmr9
  \font\eightrm=cmr8      
  \font\sevenrm=cmr7
  \font\sixrm=cmr6

  \font\seventeenpcap=cmcsc10 at 17pt
  \font\tenpcap=cmcsc10                        
  \font\ninepcap=cmcsc9
  \font\eightpcap=cmcsc8
  \font\sevenpcap=cmcsc10 scaled 700
  
  \newfam\pcapfam\textfont\pcapfam=\tenpcap\scriptfont\pcapfam=\sevenpcap
  \def\pcap{\fam\pcapfam\tenpcap}
  
  \font\seventeenrm=cmbx12 scaled 1400

  \font\fourteenbf=cmbx10 scaled 1400
  
  \font\twelvebf=cmbx12
  \font\elevenbf=cmbx10 at 11pt
  \font\ninebf=cmbx9  
  \font\eightbf=cmbx8
  \font\sixbf=cmbx6
  
  \font\tengot=eufm10                           
   
  \font\eightgot=eufm10 at 8truept 
  \font\sevengot=eufm7 
  \font\sixgot=eufm10 at 6 truept 
   
  \newfam\gotfam
  \textfont\gotfam=\tengot\scriptfont\gotfam=\sevengot\scriptscriptfont\gotfam=\sixgot
  \def\got{\fam\gotfam\tengot}%

  
  \def\tit{%
  \textfont0=\seventeenrm\scriptfont0=\tenrm\def\rm{\fam0\seventeenrm}%
  \textfont1=\seventeenib\scriptfont1=\twelveib%
  \textfont2=\seventeensy\scriptfont2=\twelvesy\scriptscriptfont2=\ninesy
  \textfont3=\seventeenex
  \textfont\itfam=\seventeenti
  \def\it{\fam\itfam\seventeenti}%
  \textfont\bbfam=\seventeenbb \scriptfont\bbfam=\twelvebb
  \def\bb{\fam\bbfam\seventeenbb}%
  \textfont\msafam=\seventeenmsa\scriptfont\msafam=\twelvemsa
  \textfont\scalnfam=\seventeenscaln
  \def\pcap{\fam\pcapfam\seventeenpcap}
  \normalbaselineskip=25pt\normalbaselines\rm}

  \font\seventeenti=cmbxti10 scaled 1680
  
  \font\fourteenti=cmbxti10 at 14pt
  
  \font\twelveti=cmbxti10 scaled 1200
  \font\eleventi=cmbxti10 at 11pt

  %
  %
  \font\twelveit=cmti12	
  \font\elevenit=cmti10 scaled 1100
  \font\nineit=cmti9
  \font\eightit=cmti8
  \font\sevenit=cmti7

  %
  %
 
 \font\seventeenib=cmmib10 scaled 1680
  \font\fourteenib=cmmib10 scaled 1400
  \font\twelveib=cmmib10 scaled 1200
  \font\elevenib=cmmib10 scaled 1100
  \font\tenib=cmmib10
\font\eightib=cmmib10 scaled 800
  \font\nineib=cmmib10 scaled 900
\font\sevenib=cmmib10 scaled 700
\font\sixib=cmmib10 scaled 600
\font\fiveib=cmmib10 scaled 500

\ifx\ITAN\oui
\else
\innernewfam\cmmibfam
\textfont\cmmibfam=\tenib
\scriptfont\cmmibfam=\sevenib
\scriptscriptfont\cmmibfam=\fiveib
\def\ib{\fam\cmmibfam\tenib}
\fi

  %
  %
  
  \font\eleveni=cmmi10 scaled 1100
  \font\ninei=cmmi9
  \font\eighti=cmmi8 
  \font\seveni=cmmi7 			                
  \font\sixi=cmmi6
  
  \font\ninesl=cmsl9                    
  \font\eightsl=cmsl8 
  \font\sevensl=cmsl10 at 7pt

  \font\ninett=cmtt9                    
  \font\eighttt=cmtt8
  \font\seventt=cmtt10 scaled 700

  \font\seventeensy=cmsy10 scaled 1680    
  \font\fourteensy=cmsy10 scaled 1400
  \font\twelvesy=cmsy10 scaled 1176
  
  \font\ninesy=cmsy9                      
  \font\eightsy=cmsy8
  \font\sixsy=cmsy6
  \font\seventeenex=cmex10 at 17pt
  \font\fourteenex=cmex10 at 14pt
  \font\twelveex=cmex10 at 12pt
  \font\nineex=cmex10 at 9pt
  \font\eightex=cmex10 at 8pt
  \font\sevenex=cmex10 at 7pt
  \font\sixex=cmex10 at 6pt
  \font\fiveex=cmex10 at 5pt
  
   
  \font\fourteengp=cmmi10 at 14pt
  
  \font\twelvegp=cmmib10 at 12pt
  \font\elevengp=cmmib10 at 11pt
  \font\tengp=cmmib10                          
  \font\ninegp=cmmib10 at 9pt 
  \font\eightgp=cmmib8 
   
  \font\sixgp=cmmib6


  \def\gponze{\textfont0=\elevenbf\scriptfont0=\eightbf\scriptscriptfont0=\sixbf
           \textfont1=\elevengp\scriptfont1=\eightgp\scriptscriptfont1=\sixgp}
  \def\gpdouze{\textfont0=\twelvebf\scriptfont0=\tenbf\scriptscriptfont0=\ninebf
           \textfont1=\twelvegp\scriptfont1=\tengp\scriptscriptfont1=\ninegp}        
  
 \def\gpquatorze{\textfont0=\fourteenbf\scriptfont0=\twelvebf\scriptscriptfont0=\elevenbf
           \textfont1=\fourteengp\scriptfont1=\twelvegp\scriptscriptfont1=\elevengp}

  
  \expandafter\chardef\csname pre amssym.def at\endcsname=\the\catcode`\@
  \catcode`\@=11
  \def\undefine#1{\let#1\undefined}
  \def\newsymbol#1#2#3#4#5{\let\next@\relax
   \ifnum#2=\@ne\let\next@\msafam@\else
   \ifnum#2=\tw@\let\next@\bbfam@\fi\fi
   \mathchardef#1="#3\next@#4#5}
  \def\mathhexbox@#1#2#3{\relax
   \ifmmode\mathpalette{}{\m@th\mathchar"#1#2#3}%
   \else\leavevmode\hbox{$\m@th\mathchar"#1#2#3$}\fi}
  \def\hexnumber@#1{\ifcase#1 0\or 1\or 2\or 3\or 4\or 5\or 6\or 7\or 8\or
   9\or A\or B\or C\or D\or E\or F\fi}
  
  \def\setboxz@h{\setbox\z@\hbox}
  \def\wdz@{\wd\z@}
  \def\boxz@{\box\z@}
  
  \edef\msafam@{\hexnumber@\msafam}
  \mathchardef\dabar@"0\msafam@39
  
  \edef\bbfam@{\hexnumber@\bbfam}
  \def\widehat#1{\setboxz@h{$\m@th#1$}%
   \ifdim\wdz@>\tw@ em\mathaccent"0\bbfam@5B{#1}%
   \else\mathaccent"0362{#1}\fi}
  \def\widetilde#1{\setboxz@h{$\m@th#1$}%
   \ifdim\wdz@>\tw@ em\mathaccent"0\bbfam@5D{#1}%
   \else\mathaccent"0365{#1}\fi}
  \newsymbol\leqq 1335          
  \newsymbol\leqslant 1336
  \newsymbol\lessgtr 1337       
  \newsymbol\backprime 1038     
  \newsymbol\risingdotseq 133A  
  \newsymbol\fallingdotseq 133B 
  \newsymbol\succcurlyeq 133C   
  \newsymbol\geqq 133D          
  \newsymbol\geqslant 133E
  \newsymbol\nmid 232D
  \newsymbol\nexists 2040
  \newsymbol\smallsetminus 2272
  \newsymbol\varnothing 203F
  
  \catcode`\@=\active

  \catcode`\@=11
  \newcount\typofr\typofr=1
  
  \catcode`\;=\active
  \def;{\ifnum\typofr=1\relax\ifhmode\ifdim\lastskip>\z@\unskip\fi
     \kern.2em\fi\string;\else\string;\fi}
  
  \catcode`\:=\active
  \def:{\ifnum\typofr=1\relax\ifhmode\ifdim\lastskip>\z@\unskip\fi
  \penalty\@M\ \fi\string:\else\string:\fi}
  
  \catcode`\!=\active
  \def!{\ifnum\typofr=1\relax\ifhmode\ifdim\lastskip>\z@\unskip\fi
     \kern.2em\fi\string!\else\string!\fi}
  
  \catcode`\?=\active
  \def?{\ifnum\typofr=1\relax\ifhmode\ifdim\lastskip>\z@\unskip\fi
     \kern.2em\fi\string?\else\string?\fi}

  \def\francais{\typofr=1\def\tpf{oui}}
  \def\anglais{\typofr=2\def\tpf{non}\def\english{oui}}
  \def\oui{oui}
  \francais
  
  \catcode`\@=12
  

\ifx\textures\oui
\def\raye #1|{\leavevmode\setbox1=\hbox{#1}%
\raise .5pt\hbox to \wd1{\xleaders\hbox{\rge{$ \sct / $}%
\kern 1pt}\hfill\kern -1pt }\kern -\wd1 \unhbox1\relax }

\def\barre #1|{\leavevmode\setbox1=\hbox{#1}%
\rlap{\Red\vrule height 2.4pt depth -1.2pt width \wd1}\Black \unhbox1\relax}
\else
\def\raye #1|{\leavevmode\setbox1=\hbox{#1}%
\raise .5pt\hbox to \wd1{\xleaders\hbox{\rge{$ \sct / $}%
\kern 1pt}\hfill\kern -1pt }\kern -\wd1 \unhbox1\relax }

\def\barre #1|{\leavevmode\setbox1=\hbox{#1}%
\rlap{\color{red}\vrule height 2.4pt depth -1.2pt width \wd1}\color{black} \unhbox1\relax}

\fi
  

  
  \def\og{\leavevmode\raise.24ex\hbox{$\scriptscriptstyle\langle\!\langle\>$}}    
  \def\fg{\leavevmode\raise.24ex\hbox{$\scriptscriptstyle\>\rangle\!\rangle$}}    

  \def\d{\,{\rm d}}
  \def\dd{{\rm d}}
  \def\dt{\d t}

  \def\dy{\d y}
  \def\z{{\bb Z}}
  \def\r{{\bb R}}
  \def\CC{{\bb C}}
  \def\N{{\bb N}}

  \def\D{{\scal D}}

  \def\HH{{\scal H}}
  \def\I{{\scal I}}
  \def\J{{\scal J}}
  \def\K{{\scal K}}
  
  \def\L{{\scal L}}
  \def\M{{\scal M}}
  
  \def\O{{\scal O}}
  \def\P{{\scaln P}}
  
  \def\R{{\scal R}}

  \def\frac#1#2{{#1\over #2}}
  \def\di#1#2{\sct#1\atop{\sct#2}}

  \def\wt{\widetilde}
  \def\numero{n$^{\rm o}\thinspace$}
\def\numeros{n$^{\rm os}\thinspace$}
  \def\page#1{\rm p.\thinspace#1}
  \def\theoreme#1{\rm th.\thinspace#1}
  
  \def\qedbox{$\rlap{$\sqcap$}\sqcup$}           
  \def\qed{\nobreak\hfill\penalty250 \hbox{}\nobreak\hfill\qedbox\par }

  \def\fas{fonc\-tions arith\-m\'etiques}

  \def\sumast{\mathop{{\sum}^*}}
  
  \def\numero{n$^{\rm o}\thinspace$}
  \def\np{nombre premier}
  \def\nps{nombres premiers}
  \def\tnp{th\'eo\-r\`eme des nombres premiers}

  \def\¤{\S\thinspace}

  \def\¥{$\bullet$ }
  
  
  \def\e{{\rm e}}
  \def\mod{\mathop{\rm mod}\nolimits}

  \def\epsilon{\varepsilon}

  \def\phi{\varphi}
  \def\theta{\vartheta}
  \def\rho{\varrho}
  \def\dm{{\textstyle{1\over 2}}}
  \def\txt{\textstyle}
  \def\dsp{\displaystyle}
  \def\sct{\scriptstyle}
  \def\pf{\noi{\it Proof. }}
  \def\nid{\ifnum\typofr=1\par\noindent{\it D\'emonstration. }\else\pf\fi}
  \def\noi{\noindent}
  \def\rem{\ifnum\typofr=1\noi{\it Remarque.}\ \else\noi{\it Remark.}\ \fi}
  \def\rems{\ifnum\typofr=1\noi{\it Remarques.}\ \else\noi{\it Remarks.}\ \fi}
  \def\re{{\Re e\,}}
  \def\im{{\Im m\,}}
  \def\ov{\overline}
  
  \def\sset{\smallsetminus}

  \def\1{{\bf 1}}
  \def\|{\Vert}

  \def\le{\leqslant}\def\leq{\leqslant}
  \def\ge{\geqslant}\def\geq{\geqslant}
  \def\wh{\widehat}
  \def\cf{{cf.}}
  
  \def\eg{{e.g.}}
  

  \def\fl#1{\left\lfloor #1 \right\rfloor}

  \def\log{\mathop{\rm log}\nolimits}
  \def\ft#1#2{{\txt{#1\over #2}}}
  \def\fs#1#2{{\scriptstyle{#1\over #2}}}
  


\def\Vbs#1{\bigg|#1\bigg|}
\def\abs#1{\left|#1\right|}


  \def\pmb#1{\setbox0=\hbox{#1}%
  \kern-.025em\copy0\kern-\wd0\kern.05em\copy0\kern-\wd0\kern-.025em\raise .0433em\box0 }

  
  \skewchar\eighti='177 \skewchar\sixi='177
  \skewchar\eightsy='60 \skewchar\sixsy='60
  
  \def\eightpoint{%
  \textfont0=\eightrm\scriptfont0=\sixrm\scriptscriptfont0=\fiverm
  \def\rm{\fam0\eightrm}%
  \textfont1=\eighti\scriptfont1=\sixi
  \scriptscriptfont1=\fivei\def\oldstyle{\fam1\seveni}%
  \textfont2=\eightsy\scriptfont2=\sixsy\scriptscriptfont2=\fivesy
  \textfont3=\eightex\scriptfont3=\sixex
  \textfont\itfam=\eightit
  \def\it{\fam\itfam\eightit}%
  \textfont\slfam=\eightsl
  \def\sl{\fam\slfam\eightsl}%
  \textfont\bbfam=\eightbb \scriptfont\bbfam=\sixbb\scriptscriptfont\bbfam=\fivebb
  \def\bb{\fam\bbfam\eightbb}%
  \textfont\msafam=\eightmsa\scriptfont\msafam=\sixmsa
  \textfont\scalnfam=\eightscaln
  \def\scaln{\fam\scalnfam\eightscaln}
  \textfont\ttfam=\eighttt
  \def\tt{\fam\ttfam\eighttt}%
\textfont\gotfam=\eightgot
  \textfont\bffam=\eightbf\scriptfont\bffam=\sixbf\scriptscriptfont\bffam=\fivebf
  \def\bf{\fam\bffam\eightbf}%
  \ifx\ITAN\oui\else\textfont\cmmibfam=\eightib
       \scriptfont\cmmibfam=\sixib
        \scriptscriptfont\cmmibfam=\fiveib
         \def\ib{\fam\cmmibfam\eightib}
   \fi
  \textfont\pcapfam=\eightpcap
  \def\pcap{\fam\pcapfam\eightpcap}
  \abovedisplayskip=2pt plus2pt minus 2pt
  \belowdisplayskip=2pt plus1pt minus 2pt
  \abovedisplayshortskip= 1pt plus 2pt minus 1pt
  \belowdisplayshortskip= 1pt plus 2pt minus 1pt
  \smallskipamount=2pt plus 1pt minus 2pt
  \medskipamount=3pt plus 2pt minus 2pt
  \bigskipamount=7pt plus 3pt minus 3pt
  \setbox\strutbox=\hbox{\vrule height 5pt depth 2pt width 0pt}%
  \normalbaselineskip=9pt\normalbaselines\rm}

  \def\({\left(}
  \def\){\right)}
  
  \def\footnoterule{\kern -2pt\hrule width 7truecm\kern 2.4pt}
  
  \def\xnotedef#1{\definexref{#1}{\noexpand\number\footnotenumber}{Note}}%

  
  
  \def\ninepoint{%
  \textfont0=\ninerm\scriptfont0=\sixrm\scriptscriptfont0=\fiverm
  \def\rm{\fam0\ninerm}%
  \textfont1=\ninei\scriptfont1=\sixi
  \scriptscriptfont1=\fivei\def\oldstyle{\fam1\ninei}%
  \textfont2=\ninesy\scriptfont2=\sixsy\scriptscriptfont2=\fivesy
  \textfont3=\nineex\scriptfont3=\sixex
  \textfont\itfam=\nineit
  \def\it{\fam\itfam\nineit}%
  \textfont\slfam=\ninesl
  \def\sl{\fam\slfam\ninesl}%
  \textfont\bbfam=\ninebb\scriptfont\bbfam=\sixbb\scriptscriptfont\bbfam=\fivebb
  \def\bb{\fam\bbfam\ninebb}%
  \textfont\msafam=\ninemsa\scriptfont\msafam=\sixmsa\scriptscriptfont\msafam=\fivemsa
  \textfont\scalnfam=\ninescaln
  \def\scaln{\fam\scalnfam\ninescaln}
  \textfont\ttfam=\ninett
  \def\tt{\fam\ttfam\ninett}%
  \textfont\bffam=\ninebf\scriptfont\bffam=\sixbf\scriptscriptfont\bffam=\fivebf
  \def\bf{\fam\bffam\ninebf}%
  \abovedisplayskip=3pt plus2pt minus 2pt
  \belowdisplayskip=3pt plus1pt minus 2pt
  \abovedisplayshortskip= 2pt plus 2pt minus 1pt
  \belowdisplayshortskip= 2pt plus 2pt minus 1pt
  \smallskipamount=2pt plus 1pt minus 2pt
  \medskipamount=3pt plus 2pt minus 2pt
  \bigskipamount=7pt plus 3pt minus 3pt
  \setbox\strutbox=\hbox{\vrule height 5pt depth 2pt width 0pt}%
  \normalbaselineskip=10.5pt plus.3pt minus.3pt\normalbaselines\rm}

  \def\sevenpoint{%
  \textfont0=\sevenrm\scriptfont0=\sixrm\scriptscriptfont0=\fiverm
  \def\rm{\fam0\sevenrm}%
  \textfont1=\seveni\scriptfont1=\sixi
  \scriptscriptfont1=\fivei\def\oldstyle{\fam1\seveni}%
  \textfont2=\sevensy\scriptfont2=\sixsy\scriptscriptfont2=\fivesy
  \textfont3=\sevenex\scriptfont3=\fiveex
  \textfont\itfam=\sevenit
  \def\it{\fam\itfam\sevenit}%
  \textfont\slfam=\sevensl
  \def\sl{\fam\slfam\sevensl}%
  \textfont\bbfam=\sevenbb \scriptfont\bbfam=\sixbb\scriptscriptfont\bbfam=\fivebb
  \def\bb{\fam\bbfam\sevenbb}%
  \textfont\msafam=\sevenmsa\scriptfont\msafam=\sixmsa
  \textfont\scalnfam=\sevenscaln
  \def\scaln{\fam\scalnfam\sevenscaln}
  \textfont\bffam=\sevenbf\scriptfont\bffam=\sixbf\scriptscriptfont\bffam=\fivebf
  \def\bf{\fam\bffam\sevenbf}%
  \textfont\ttfam=\seventt
  \abovedisplayskip=2pt plus2pt minus 2pt
  \belowdisplayskip=2pt plus1pt minus 2pt
  \abovedisplayshortskip= 1pt plus 2pt minus 1pt
  \belowdisplayshortskip= 1pt plus 2pt minus 1pt
  \smallskipamount=2pt plus 1pt minus 2pt
  \medskipamount=3pt plus 2pt minus 2pt
  \bigskipamount=7pt plus 3pt minus 3pt
  \setbox\strutbox=\hbox{\vrule height 5pt depth 2pt width 0pt}%
  \normalbaselineskip=9pt\normalbaselines\rm}

 \def\onzepts{%
 \textfont0=\elevenrm\scriptfont0=\ninerm
 \textfont1=\elevenib\scriptfont1=\ninei}

\def\douzepts{%
  \textfont0=\twelverm\scriptfont0=\tenrm\def\rm{\fam0\twelverm}%
  \textfont1=\twelveib\scriptfont1=\teni%
  \textfont2=\twelvesy\scriptfont2=\tensy\scriptscriptfont2=\eightsy
  \textfont3=\twelveex
  \textfont\itfam=\twelveti
  \def\it{\fam\itfam\twelveti}%
  \textfont\bffam=\twelvebf\scriptfont\bffam=\tenbf\scriptscriptfont\bffam=\eightbf
  \def\bf{\fam\bffam\twelvebf}%
  \textfont\bbfam=\twelvebb \scriptfont\bbfam=\tenbb
  \def\bb{\fam\bbfam\twelvebb}%
  \textfont\msafam=\twelvemsa\scriptfont\msafam=\tenmsa
  \textfont\scalnfam=\twelvescaln
  \normalbaselineskip=15pt\normalbaselines\rm}

\def\quatorzepts{%
  \textfont0=\fourteenrm\scriptfont0=\twelverm\def\rm{\fam0\fourteenrm}%
  \textfont1=\fourteenib\scriptfont1=\twelveib%
  \textfont2=\fourteensy\scriptfont2=\twelvesy\scriptscriptfont2=\tensy
  \textfont3=\fourteenex
  \textfont\itfam=\fourteenti
  \def\it{\fam\itfam\fourteenti}%
  \textfont\bffam=\fourteenbf\scriptfont\bffam=\twelvebf\scriptscriptfont\bffam=\tenbf
  \def\bf{\fam\bffam\fourteenbf}%
  \textfont\bbfam=\fourteenbb \scriptfont\bbfam=\twelvebb
  \def\bb{\fam\bbfam\fourteenbb}%
  \textfont\msafam=\fourteenmsa\scriptfont\msafam=\twelvemsa
  \textfont\scalnfam=\twelvescaln
  \normalbaselineskip=18pt\normalbaselines\rm}


\def\AA{{\it Acta Arith.}}

\def\picture #1 by #2 (#3){\leavevmode\vbox to #2{
     \hrule width #1 height 0pt depth 0pt
      \vfill
       \special{picture #3}}}

\def\illustration #1 by #2 (#3) scaled #4{\dimen1=#2
  \divide\dimen1 by 1000
  \multiply\dimen1 by #4
  \vtop to \dimen1{\dimen1=#1
  \divide\dimen1 by 1000
  \multiply\dimen1 by #4
  \hsize=\dimen1\vss
  \noindent\special{illustration #3 scaled #4}}}

\ifx\optionkeymacros\undefined\else \fi

\catcode`\Œ=\active\defŒ{{\aa}}       
\catcode`\º=\active\defº{\int}        
\catcode`\=\active\def{\c c}        
\catcode`\¶=\active\def¶{\partial}    
\catcode`\Ä=\active\defÄ{\oint}       
\catcode`\Æ=\active\defÆ{\triangle}   
\catcode`\Â=\active\defÂ{\neg}        
\catcode`\µ=\active\defµ{\mu}         
\catcode`\¿=\active\def¿{{\o}}        
\catcode`\¹=\active\def¹{\pi}         
\catcode`\Ï=\active\defÏ{{\oe}}       
\catcode`\§=\active\def§{{\ss}}       
\catcode`\ =\active\def {\dagger}     
\catcode`\Ã=\active\defÃ{\sqrt}       
\catcode`\·=\active\def·{\Sigma}      
\catcode`\Å=\active\defÅ{\approx}     
\catcode`\½=\active\def½{\Omega}      
\catcode`\£=\active\def£{{\it\$}}     
\catcode`\°=\active\def°{\infty}      
\catcode`\¤=\active\def¤{{\S}}        
\catcode`\¦=\active\def¦{{\P}}        
\catcode`\¥=\active\def¥{\bullet}     
\catcode`\»=\active\def»{\leavevmode\raise.585ex\hbox{\b a}}      
\catcode`\¼=\active\def¼{\leavevmode\raise.6ex\hbox{\b o}}        
\catcode`\­=\active\def­{\not=}       
\catcode`\²=\active\def²{\leq}        
\catcode`\³=\active\def³{\geq}        
\catcode`\Ö=\active\defÖ{\div}        
\catcode`\É=\active\defÉ{{\dots}}     
\catcode`\¾=\active\def¾{{\ae}}       
\catcode`\Ç=\active\defÇ{\og}         
\catcode`\Ò=\active\defÒ{``}          
\catcode`\Á=\active\defÁ{!`}          
\catcode`\¢=\active\def¢{\rlap/c}     
\catcode`\Ô=\active\defÔ{`}           
\catcode`\Õ=\active\defÕ{'}           


\catcode`\=\active\def{{\AA}}       
\catcode`\'=\active\def'{\c C}        
\catcode`\¯=\active\def¯{{\O}}        
\catcode`\¸=\active\def¸{\Pi}         
\catcode`\Î=\active\defÎ{{\OE}}       
\catcode`\®=\active\def®{{\AE}}       
\catcode`\×=\active\def×{\diamond}    
\catcode`\¡=\active\def¡{\accent'27}  
\catcode`\Ó=\active\defÓ{''}          
\catcode`\±=\active\def±{\pm}         
\catcode`\È=\active\defÈ{\fg}         
\catcode`\À=\active\defÀ{?`}          
\catcode`\Ð=\active\defÐ{--}          
\catcode`\Ñ=\active\defÑ{---}         


\catcode`\Š=\active\defŠ{\"a}        
\catcode`\'=\active\def'{\"e}        
\catcode`\•=\active\def•{\"{\i}}     
\catcode`\š=\active\defš{\"o}        
\catcode`\Ÿ=\active\defŸ{\"u}        
\catcode`\Ø=\active\defØ{\"y}        
\catcode`\å=\active\defå{\^A}        
\catcode`\€=\active\def€{\"A}        
\catcode`\…=\active\def…{\"O}        
\catcode`\†=\active\def†{\"U}        
\catcode`\‡=\active\def‡{\'a}        
\catcode`\Ž=\active\defŽ{\'e}        
\catcode`\'=\active\def'{\'{\i}}     
\catcode`\—=\active\def—{\'o}        
\catcode`\œ=\active\defœ{\'u}        
\catcode`\ƒ=\active\defƒ{\'E}        
\catcode`\æ=\active\defæ{\^E}        
\catcode`\é=\active\defé{\`E}        %
\catcode`\ˆ=\active\defˆ{\`a}        
\catcode`\=\active\def{\`e}        
\catcode`\"=\active\def"{\`{\i}}     
\catcode`\˜=\active\def˜{\`o}        
\catcode`\=\active\def{\`u}        
\catcode`\Ë=\active\defË{\`A}        
\catcode`\‹=\active\def‹{\~a}        
\catcode`\–=\active\def–{\~n}        
\catcode`\›=\active\def›{\~o}        
\catcode`\Ì=\active\defÌ{\~A}        
\catcode`\"=\active\def"{\~N}        
\catcode`\Í=\active\defÍ{\~O}        
\catcode`\‰=\active\def‰{\^a}        
\catcode`\=\active\def{\^e}        
\catcode`\"=\active\def"{\^{\i}}     
\catcode`\™=\active\def™{\^o}        
\catcode`\ž=\active\defž{\^u}        

\let\optionkeymacros\null

\optionparag=1
\def\paradouze{oui}

\ifx\montrerlabels\oui
\input montrerlabels.tex
\fi

\dimstand
\voffset=-4mm
\abovedisplayskip=7pt plus4pt minus 4pt
 \belowdisplayskip=7pt plus3pt minus 3pt
  \abovedisplayshortskip= 3pt plus 2pt minus 1.5pt
  \belowdisplayshortskip= 3pt plus 2pt minus 1.5pt
  \smallskipamount=3pt plus 1pt minus 2pt
  \medskipamount=4pt plus 2pt minus 2.5pt
  \bigskipamount=8pt plus 2pt minus 4pt

\newcount\numchapi\numchapi=0
\def\paraunn#1{\paraun{#1}\writetocentry{section}{#1}}
\def\paradeuxn#1{\paradeuxb{#1}\writetocentry{subsection}{#1}}

\scriptscriptfont0=\fiverm

    \font\tenrsfs=rsfs7 at 10pt

    \font\sevenrsfs=rsfs7
    \font\sixrsfs=rsfs7 at 6pt

\newfam\rsfsfam\textfont\rsfsfam=\tenrsfs\scriptfont\rsfsfam=\sevenrsfs\scriptscriptfont\rsfsfam=\sixrsfs
    \def\rsfs{\fam\rsfsfam\tenrsfs}%

\def\L{{\rsfs L}}

\def\wt{\widetilde}
\def\rd{{\rm d}}
\def\rst{\Big\{\epsilon^{\varrho_0}+\Big({\log y\over\log x}\Big)^{\varrho_0}\Big\}}
\def\rstt#1{\Big\{\epsilon^{#1}+\Big({\log y\over\log x}\Big)^{#1}\Big\}}

\def\pnu{p^\nu}
\def\pmu{p^\mu}
\def\qmu{q^\mu}
\def\pk{p^k}
\def\fmu{fonc\-tion mul\-ti\-pli\-cative}
\def\fmus{fonc\-tions mul\-ti\-pli\-catives}
\def\fad{fonc\-tion ad\-di\-tive}

\def\I{{\scal I}}
\def\J{{\scal J}}
\def\ga{{\got a}}
\def\gb{{\got b}}
\def\gc{{\got c}}

\def\gh{{\got h}}

\def\gt{{\got t}}
\def\gE{{\got E}}
\def\gG{{\got G}}
\def\gR{{\got R}}
\def\gS{{\got S}}
\def\gT{{\got T}}
\def\gV{{\got V}}
\def\gW{{\got W}}
\def\gpp{{\got p}}
\def\Halasz{Hal‡sz}
\def\Sarkozy{S‡rkšzy}
{\leftskip-5mm\obeylines
{\it Ramanujan J. 
\bf44}, \numero3 (2017), 641-701;
Corrig. {\bf 51}, \numero1 (2020), 243-244.\par }
\hautspages{GŽrald Tenenbaum}{Moyennes effectives de fonctions multiplicatives complexes}
\vskip0mm
\titrecentre{Moyennes effectives de fonctions multiplicatives complexes\anote{*}{\hskip-4mm Nous incluons ici quelques corrections mineures relativement ˆ la version publiŽe.\hfill}}
\bigskip\smallskip
\centerline{GŽrald Tenenbaum}
\bigskip\bigskip\medskip
%
%
\newcount\paras\paras=0
\newcount\sparas\sparas=0
\def\tocsectionentry#1#2{\init{\sparas}\avance\paras
       {\quad\bf\the\paras\quad }{\hskip-2mm #1\dotfill\hskip3mm\rm#2}
\par}%
\def\tocsubsectionentry#1#2{\avance\sparas
       {\qquad\eightpoint\the\paras.\the\sparas}
{\hskip-2mm\eightpoint#1\dotfill\hskip3mm\rm#2}\par}%

{\eightpoint\leftskip1cm\rightskip1cm
\anglais
\noi{\bf Abstract.} We establish effective mean-value estimates for a wide class of multiplicative arithmetic functions, thereby providing (essentially optimal) quantitative versions of Wirsing's classical estimates and extending those of \Halasz. Several applications are derived, including: estimates for the difference of mean-values of so-called pretentious functions, local laws for the distribution of prime factors in an arbitrary set, and weighted distribution of additive functions.
\medskip
\noi \bf Keywords. \rm Quantitative estimates, multiplicative functions, effective mean-value theorems, weighted distribution of additive functions.
\par \smallskip
\bf \noi 2010 AMS Classification. \rm Primary 11N56, Secondary 11K65, 11N37, 11N60, 11N64.
\par }
\medskip

\noi{\qquad\quad \bf Sommaire}\par \smallskip
{\eightpoint\leftskip.6cm\rightskip.8cm
\baselineskip=8.5pt\readtocfile}

\bigskip\bigskip

\paraunn{Introduction et ŽnoncŽ des rŽsultats}
Les estimations de valeurs moyennes de fonctions multiplicatives constituent un outil pri\-vilŽgiŽ de la thŽorie probabiliste des  nombres. Elles permettent notamment d'apprŽhender la rŽpartition des fonctions additives sur les $N$ premiers entiers via leurs fonctions caractŽristiques et, partant, d'obtenir des thŽormes de convergence avec contr™le de l'approximation.
Les deux succs historiques de la thŽorie sont respectivement dus ˆ Wirsing \citer{Wi67} et \Halasz~\citer{Ha68}.\par  DŽsignons par $\M(A,B)$ la classe des \fmus\ vŽrifiant 
$$\max_{p}|f(p)|\leqslant A,\quad  \sum_{p,\,\nu\geqslant 2}{|f(p^\nu)|\log p^\nu\over p^\nu}\leqslant B.\eqdef{C0}
$$
Ici et dans la suite, nous rŽservons la lettre~$p$ pour dŽsigner un \np.\par 
Dans son remarquable article \citer{Wi67}, Wirsing Žtablit notamment que, si $r\in\M(A,B)$, $r\geqslant 0$, et s'il existe $\varrho>0$ tel que
$$\sum_{p\leqslant x}{r(p)\log p\over p}\sim \varrho\log x\qquad (x\to\infty),\eqdef{Wi1}$$
alors toute \fmu\ rŽelle $f$ telle que $|f|\leqslant r$ vŽrifie, lorsque $x\to\infty$,
$$M(x;f):=\sum_{n\leqslant x}f(n)=\bigg\{{\e^{-\gamma\varrho}\over \Gamma(\varrho)}\prod_{p}{\sum_{\nu\geqslant 0}{f(\pnu)/\pnu}\over \sum_{\nu\geqslant 0}r(\pnu)/\pnu}+o(1)\bigg\}{x\over\log x}\prod_{p} \sum_{\pnu\leqslant x}{r(\pnu)\over \pnu},\eqdef{Wi2}$$
o le produit infini est considŽrŽ comme nul lorsqu'il diverge.\note{Le second produit est en fait fini.} Ici et dans la suite, nous notons $\gamma$ la constante d'Euler. \par \goodbreak
Le cas $r=\1$ (la fonction constante Žgale ˆ 1), $\varrho=1$, confirme une cŽlbre conjecture d'Erd\H os selon laquelle une \fmu\ rŽelle ˆ valeurs dans $[-1,1]$ possde nŽcessairement une valeur moyenne. 

\par 
Dans \citer{Ha68}, \Halasz\ a ŽlucidŽ le comportement asymptotique des \fmus\ complexes $f$ ˆ valeurs dans le disque unitŽ. Son rŽsultat principal Žtablit une dichotomie: soit il existe $\tau\in\r$ tel que
$$\sum_{p}{1-\re\{f(p)/p^{i\tau}\}\over p}<\infty,\eqdef{Ha0}$$
et nous avons alors
$$M(x;f)\sim {x^{1+i\tau}\over 1+i\tau}\prod_{p\leqslant x}\Big(1-{1\over p}\Big)\sum_{\nu\geqslant 0}{f(\pnu)\over p^{\nu(1+i\tau)}}+o(x)\qquad (x\to\infty),\eqdef{Ha1}$$
soit la sŽrie diverge pour tout $\tau\in\r$, et l'on a $M(x;f)=o(x)$ lorsque $x\to\infty$. Une prŽcision supplŽmentaire est que, sous l'hypothse \eqref{Ha0}, le terme principal de \eqref{Ha1} est de la forme $K_\tau x^{1+i\tau}L_\tau(\log x)$ o $L_\tau$ est une fonction de module unitŽ ˆ croissance lente au sens de Karamata \citer{Ka30}, c'est-ˆ-dire telle que $L_\tau(u)/L_\tau(v)\to1$ lorsque $u$ et $v$ tendent vers l'infini sous la condition~$u\asymp v$.\par \goodbreak
Une preuve alternative du cas de convergence a ŽtŽ obtenue par Delange via une mŽthode reposant, ˆ partir d'une idŽe de RŽnyi, sur l'inŽgalitŽ de Tur‡n--Kubilius: voir \citer{Te15}, th.\thinspace III.4.4, pour une dŽmonstration de ce rŽsultat non publiŽe par Delange. 
\par \goodbreak
Indlekofer, K‡tai \& Wagner \citer{IKW01} ont gŽnŽralisŽ ces rŽsultats\note{Un rŽsultat qualitatif antŽrieur,  de mme nature mais valide sous des hypothses plus fortes, est dž ˆ Levin \& Timofeev \citer{LT82}.} en Žtablissant que, pour toutes \fmus\ $r,\,f$, telles que $r\in\M(A,B)$, $f$ ˆ valeurs complexe, $|f|\leqslant r$, et sous l'hypothse \eqref{Wi1}, nous avons $M(x;f)=o\big(M(x;r)\big)$ lorsque $x\to\infty$ si la condition
$$\sum_{p}{r(p)-\re\{f(p)/p^{i\tau}\}\over p}<\infty\eqdef{Ha2}$$ n'est rŽalisŽe pour aucune valeur de $\tau\in\r$.
Ils Žnoncent Žgalement que, sous la condition \eqref{Ha2} nous avons, lorsque $x\to\infty$,
$$M(x;f)=\bigg\{\prod_{p\leqslant x}{\sum_{\nu\geqslant 0}{f(\pnu)/p^{\nu(1+i\tau)}}\over \sum_{\nu\geqslant 0}r(\pnu)/\pnu}+o(1)\bigg\}{x^{i\tau}M(x;r)\over 1+i\tau}\cdot\eqdef{IKW}$$
Cette formule asymptotique  rŽsulte en fait implicitement, sous la mme hypothse \eqref{Ha2}, du travail de Wirsing \citer{Wi67}. Notons Žgalement que les estimations de Wirsing impliquent
  l'Žqui\-valence de \eqref{IKW} et
$$M(x;f)={\e^{-\gamma\varrho}x^{1+i\tau}\over (1+i\tau)\Gamma(\varrho)\log x}\bigg\{\prod_{p\leqslant x}\sum_{\nu\geqslant 0}{f(\pnu)\over p^{\nu(1+i\tau)}}+o\big((\log x)^\varrho\big)\bigg\}.\eqdef{IKW2}$$
De plus, le cas gŽnŽral se rŽduit aisŽment, par sommation d'Abel, au cas $\tau=0$.
\par 

Les dŽveloppements ultŽrieurs de la thŽorie ont principalement visŽ ˆ rendre les rŽsultats prŽcŽdents effectifs, c'est-ˆ-dire ˆ expliciter les majorations dans le cas de divergence et les termes d'erreur dans le cas de convergence. La nŽcessitŽ de telles estimations Žtant naturellement plus impŽrieuse dans le cas o le terme principal est nul, les recherches se sont d'abord orientŽes dans cette direction. Les travaux de Hal‡sz~\citer{Ha71}, prŽcisŽs par Montgomery \citer{Mo78}, Elliott \citer{El80}, puis par l'auteur \citer{Te15} (ch. III.4), fournissent ainsi des majorations explicites dans le cas de fonctions ˆ valeurs dans le disque unitŽ ou dont les  valeurs sur les \nps\ Žvitent un secteur fixe.\note{Voir Žgalement \citer{MV01} pour des variantes relatives ˆ des sommes pondŽrŽes.}
 \par\goodbreak
  Lorsque les $f(p)$ sont confinŽs ˆ une ellipse de diamtre 2 strictement incluse dans le disque unitŽ, Hall \& Tenenbaum \citer{HT91} Žtablissent la majoration effective
$$M(x;f)\ll x\exp\Big\{-K\sum_{p\leqslant x}{1-\re f(p)\over p}\Big\},$$ 
o la constante $K$ est optimale. Ce rŽsultat a ŽtŽ gŽnŽralisŽ par Hall \citer{Ha95} au cas o les $f(p)$ sont confinŽs ˆ sous-ensemble strict et fermŽ du disque unitŽ Žvitant au moins un point de module 1. Des complŽments et raffinements sont proposŽs dans l'article exhaustif de Granville et Soundararajan \citer{GS03}, qui contient Žgalement des dŽveloppements relatifs au comportement local les moyennes $M(x;f)$.
\par 
Dans la voie d'une version quantitative des estimations de Wirsing, un rŽsultat de \Halasz\ (\citer{Ha71}, \theoreme 3) fournit  une formule asymptotique avec terme principal non nul lorsque les valeurs aux \nps\ sont proches de $1$. Une version relative aux fonctions ˆ valeurs dans $[-1,1]$ a Žgalement ŽtŽ donnŽe par Indlekofer \citer{In09} --- \cf\ le \ref{Wieff} {\it infra}. 
\par \medskip
Nous nous proposons ici de prŽciser les rŽsultats antŽrieurs dans deux directions : d'une part en Žtendant les majorations effectives aux fonctions des classes $\M(A,B)$, d'autre part en fournissant des versions quantitatives de l'ensemble des estimations de type Wirsing telles qu'Žtablies sous forme qualitative par Indlekofer, K‡tai \& Wagner dans \citer{IKW01}.\par 
Un tel programme suppose que les \fmus\ considŽrŽes soient autorisŽes ˆ dŽpendre du paramtre $x$ gouvernant la taille de la moyenne prise en considŽration. Dans la suite, nous considŽrons donc, pour tous paramtres $A>0$, $B>0$, $j=0,1$, $x\geqslant 2$, la classe $\M_j(x;A,B)$ des fonctions multiplicatives complexes $f$ vŽrifiant
$$\max_{p\leqslant x}|f(p)|\leqslant A,\quad  \sum_{\di{p^\nu\le
x}{\nu\geqslant 2}}{|f(p^\nu)|(\log p^\nu)^{1+j}\over p^\nu}\leqslant B.\eqdef{C1}
$$\par\goodbreak 
Notre premier rŽsultat Žtend aux fonctions de $\M_1(x;A,B)$ le thŽorme III.4.7 de \citer{Te15}, restreint aux fonctions  ˆ valeurs dans le disque unitŽ. Pour toute fonction multiplicative $f$ dont la sŽrie de Dirichlet
$\sum_{n\geqslant 1}{f(n)/ n^s}$
converge dans le demi-plan $\re s>1$, nous posons
$$\eqalign{  &v_f(s)=v_f(s;x):=\sum_{p\leqslant x}{f(p)\over p^s}\qquad (s\in\CC),\cr& H_T(\alpha)^2:=\sum_{\di{k\in\z}{|k|\leqslant T}}{1\over k^2+1}\max_{\di{\sigma=1+\alpha}{|\tau-k|\leqslant 1/2}}|\e^{v_f(s;x)}|^2\quad(\alpha>0,\,T\geqslant 1).\cr}\eqdef{MphiH}$$
Dans tout ce travail, nous dŽfinissons implicitement les parties rŽelle et imaginaire d'un nombre complexe $s$ par la formule $s=\sigma+i\tau$.
\par Nous posons encore
$$Z(y;f):=\sum_{p\leqslant y}{f(p)\over p}\qquad (2\leqslant y\leqslant x).\eqdef{defZZ}$$

\Propt{casdiv}{Soient $A>0,\, B>0$. Sous les hypothses $x\geqslant 3$, $f,r\in\M_1(x;A,B)$, $|f|\leqslant r$, et $T\geqslant 1$, nous avons uniformŽment
$$M(x;f)\ll {x\over \log x}\bigg\{\int_{1/\log x}^1{H_T(\alpha)\over \alpha}\d\alpha+{\e^{Z(x;r)}\over \sqrt{T}}+{\e^{Z(x;r)}\over \log x}+{\e^{Z(x;r)}\log_2x\over T}\bigg\}.\eqdef{majeff}$$
De plus, pour tout $c>0$ fixŽ, et sous l'hypothse supplŽmentaire$$Z(x;r)-Z(y;r)\geqslant c\log \Big({\log x\over \log y}\Big)+O(1)\qquad (2\leqslant y\leqslant x),\eqdef{hypsupdiv}$$ le dernier terme dans l'accolade de \eqref{majeff} peut tre omis.}\goodbreak
Le rŽsultat suivant fournit une version effective des formules asymptotiques \eqref{IKW} et \eqref{IKW2}. Confor\-mŽment ˆ une remarque effectuŽe plus haut, nous nous restreignons, sans perte de gŽnŽralitŽ, au cas $\tau=0$.\par 
ƒtant donnŽe une fonction arithmŽtique multiplicative complexe~$f$,
nous posons $w_f:=1$ si $f$ est rŽelle, et $w_f:=\dm$ dans le cas gŽnŽral. 
\par 
\goodbreak
\Propt{cascv}{Soient $$\eqalign{&\ga\in]0,\dm],\ \gb\in\,[\ga,1[,\  A\geqslant2\gb,\  B>0,\cr
&x\geqslant 1,\ \varrho=\varrho_x\in\,[2\gb,A],\
\gpp:={\pi\varrho\over A},\ \beta:=1-{\sin\gpp\over \gpp},\quad\gh:={1-\gb\over \min(1,\varrho)-\gb}\cdot\cr}\eqdef{defpar}$$  Pour tout $\varepsilon=\varepsilon_x\in]1/\sqrt{\log x}, \dm]$, les assertions suivantes relatives aux \fmus\ $f$, $r$ telles que $|f|\leqslant r$ sont vŽrifiŽes. \par 
Sous les hypothses 
$$\leqalignno{\sum_{p\leqslant x}{r(p)-\re f(p)\over p}&\leqslant \ft12\beta\gb\log (1/\epsilon), &\eqdef{pqr}\cr
\sum_{x^\varepsilon<p\leqslant y}{\{r(p)-\re f(p)\}^\gh\log p\over p}&\ll \varepsilon^{\delta_1\gh}\log y\qquad (x^\varepsilon<y\leqslant x), &\eqdef{cdpr1}\cr
\sum_{p\leqslant y}{\{r(p)-\varrho\}\log p\over p}&\ll\varepsilon\log y\qquad 
(x^{\epsilon}< y\leqslant x),&\eqdef{C2}\cr
}$$ 
o
$\delta_1\in]0,\ft23\beta\gb]$, 
 nous avons, uniformŽment
pour $x\geqslant 2$, $r\in\M_0(x;A,B)$,
$$M(x;f)={\e^{-\gamma\varrho}x\over \Gamma(\varrho)\log
x}\bigg\{\prod_{p}\sum_{\pnu\leqslant 
x} {f(p^\nu)\over p^\nu}+O\Big(\varepsilon^\delta\,\e^{Z(x;f)}\Big)\bigg\},\eqdef{moy}$$
o l'on a posŽ $
\delta:=w_f\delta_1.$\par 
 La constante implicite dans \eqref{moy} dŽpend au plus  de $A$, $B$, $\ga$, 
 $\gb$, et des constantes implicites de \eqref{cdpr1} et \eqref{C2}. }
\par\medskip
\rems(i) Le \ref{cascv} gŽnŽralise bien le thŽorme de Wirsing.

\par 
(ii) On a $\gh=1$, ds que $\varrho\geqslant 1$.\par 
(iii)
L'hypothse \eqref{cdpr1} est trivialement impliquŽe par la condition
$$\sum_{x^\varepsilon<p\leqslant x}{\{r(p)-\re f(p)\}^\gh\over p}\ll  \epsilon^{\delta_1\gh},\eqdef{cdpr}$$
et, bien entendu, Žgalement par la majoration uniforme 
$$\max_{x^\varepsilon<p\leqslant x}\{r(p)-\re f(p)\}\ll \varepsilon^{\delta_1}.\eqdef{cdpr2}$$

\par

(iv) Les hypothses du \ref{cascv} impliquent
$$\prod_{p}\sum_{\pnu\leqslant 
x} {f(p^\nu)\over p^\nu}\ll\e^{Z(x;f)}$$ 
alors que les deux membres sont du mme ordre de grandeur ds que
$$\min_{p,x}\Vbs{\sum_{0\leqslant \nu\leqslant (\log x)/\log p} {f(p^\nu)\over p^\nu}}\gg1.\eqdef{hypgen}$$\goodbreak
 \noi Sous cette condition gŽnŽriquement vŽrifiŽe, la formule \eqref{moy} devient
$$M(x;f)=\big\{1+O(\varepsilon^\delta)\big\}{\e^{-\gamma\varrho}x\over \Gamma(\varrho)\log
x}\prod_{p}\sum_{\pnu\leqslant 
x} {f(p^\nu)\over p^\nu}\cdot\eqdef{moy2}$$
\par 
 \medskip\goodbreak
L'hypothse \eqref{C2} reprŽsente une contrainte significative pour la rŽpartition des valeurs~$r(p)$.  Nous pouvons la remplacer par une minoration en moyenne sur de petits intervalles.
\par 
Nous posons
\vskip-3mm
$$\beta_0=\beta_0(\gb,A):=1-{\sin(2\pi\gb/A)\over 2\pi\gb/A},\qquad  \delta_0(\gb)=\delta_0(\gb,A):=\ft13\gb\beta_0.\eqdef{delta0}$$\par 
\vskip-3mm
\goodbreak
\Propt{fr}{Soient $$\ga\in]0,\ft14], \quad\gb\in[\ga,\ft12[,  \quad A\geqslant 2\gb,\quad B>0, \quad \beta:=\beta_0(\gb,A),\quad  x\geqslant 2, \quad 1/\sqrt{\log x}<\varepsilon\leqslant \dm.$$ Supposons que les \fmus\ $f$, $r$, telles que \hbox{$r\in\M_0(x;2A,B)$}, $|f|\leqslant r$,  vŽrifient les conditions \eqref{pqr}, \eqref{cdpr1} avec $\gh:=(1-\gb)/\gb$, \eqref{cdpr} avec $\gh=1$, et 
$$\sum_{y<p\leqslant y^{1+\varepsilon_1}}{r(p)\log p\over p}\geqslant 4\gb\varepsilon_1\log y\qquad \big(\e^{1/\varepsilon_1}\leqslant y\leqslant x^{1/(1+\varepsilon_1)}\big)\eqdef{moyrp}$$
o l'on a posŽ $\varepsilon_1:=\sqrt{\varepsilon}$. Supposons de plus que $\delta_1\in]0,\delta_0(\gb)] $. Nous avons alors
$$M(x;f)=M(x;r)\prod_{p}{\sum_{\pnu\leqslant x}f(\pnu)/\pnu\over \sum_{\pnu\leqslant x}r(\pnu)/\pnu}+O\bigg({x\,\varepsilon^{\delta}\e^{Z(x;r)-\gc Z(x;|f|-f)}\over \log x}\bigg)\eqdef{Mf/Mr}$$
o $\delta:=w_f\delta_1$, et $\gc:=\gb/A$.
La constante implicite de \eqref{Mf/Mr} dŽpend au plus de $A$, $B$, $\ga$ et $\gb$.}
\par \smallskip
\rems (i) Le rŽcent et trs ŽlŽgant article d'Elliott \citer{El16}, finalisŽ simultanŽment au prŽsent travail, fournit une condition suffisante pour la validitŽ de \eqref{IKW} dans laquelle l'hypothse \eqref{Ha2} est remplacŽe par une minoration en moyenne de mme nature que \eqref{moyrp}.
\par 
(ii) Les ThŽormes \ref{scascv} et \ref{sfr}  prŽcisent le thŽorme 2.6 de l'article \citer{Ma17a}, apparu en ligne postŽrieurement ˆ la diffusion du prŽsent travail sur le rŽseau.
\par \medskip
Il est possible de supprimer la condition \eqref{cdpr} au prix d'un renforcement des hypothses sur les nombres~$r(p)$.
\vskip-3mm
\Propt{comprf}{Dans les hypothses du \ref{fr}, la formule asymptotique \eqref{Mf/Mr} persis\-te sans la condition \eqref{cdpr} avec $\gh=1$  si l'hypothse \eqref{moyrp} est remplacŽe par $\dsp\min_{x^\varepsilon<p\leqslant x}r(p)\geqslant 4\gb$.} 
Nous fournissons les dŽtails au \S\thinspace\ref{prgcomp}.
\medskip
Les mŽthodes dŽveloppŽes dans le prŽsent travail reposent principalement sur l'approche de \Halasz\ \citer{Ha68}, assortie de raffinements introduits dans \citer{Mo78} et \citer{Te15}. La preuve du \ref{fr} fait usage d'une version pondŽrŽe de l'inŽgalitŽ de Tur‡n--Kubilius (\ref{TKpond}) et d'arguments de convolution des \fas.  Ainsi que l'attestent les formules asymptotiques obtenues aux  Corollaires \ref{sOmE} et \ref{srepad}, les termes d'erreur effectifs des ThŽormes \ref{scascv}, \ref{sfr} et \ref{scomprf} sont essentiellement optimaux sous les hypothses effectuŽes. 
\medskip
\noi{\bf Notation.} Dans tout ce travail, nous employons la notation de Vinogradov $f\ll g$ pour signifier qu'il existe une constante $C$ telle que $|f|\leqslant C|g|$
 dans le domaine indiquŽ. Nous attirons l'attention du lecteur sur le fait que, selon une pratique largement rŽpandue, nous Žtendons l'usage de cette notation au cas de quantitŽs complexes.
\bigskip\medskip
\paraunn{Applications}
Nous Žnonons ici, de faon non exhaustive, quelques applications des rŽsultats prŽsentŽs plus haut.\par  La premire est une consŽquence simplifiŽe du \ref{casdiv} analogue ˆ une majoration de \Halasz\ \citer{Ha75} valable pour les fonctions de module au plus 1, et dont une version optimale est Žtablie au cor. III.4.12 de \citer{Te15}. Lorsque $r\in\M_0(x;A,B)$ et $f$ est une \fmu\ telle que $|f|\leqslant r$, nous posons
$$m_f(y;T):=\min_{|\tau|\leqslant T}\sum_{p\leqslant y}{r(p)-\re\big( f(p)/p^{i\tau}\big)\over p}\qquad (y\geqslant 2).\eqdef{defmfyT}$$
\Propc{majm}{Soient $A>0,\, B>0$, $\gb>0$. Sous les conditions $x\geqslant 3$, $f,r\in\M_1(x;A,B)$, $|f|\leqslant r$, $T\geqslant 1$,  et 
$$\sum_{y<p\leqslant x}{r(p)\over p}\geqslant \gb\log\Big({\log x\over \log y}\Big)+O(1)\qquad (2\leqslant y\leqslant x),\eqdef{rp/p}$$
nous avons uniformŽment 
$$M(x;f)\ll M(x;r)\bigg\{{1+m_f(x;T)\over \e^{m_f(x;T)}}+{1\over \sqrt{T}}+{1\over \log x}\bigg\}.\eqdef{majfm}$$}
\par \goodbreak
Une illustration simple du \ref{cascv} peut tre obtenue de la faon suivante. Soient $A$ et $B$ des constantes positives et $f$ une \fmu\ complexe de $\M_0(A,B)$ telle que $\max_p|f(p)|\leqslant \varrho$. Supposons encore que l'analogue de \eqref{Ha0} est satisfait avec $\tau=0$, autrement~dit
$$\sum_{p}{\varrho-\re f(p)\over p}<\infty.\eqdef{Ha00}$$
Il rŽsulte alors de \eqref{IKW2} et de la formule de Mertens que
$$M(x;f)=\bigg\{\prod_{p\leqslant x}\Big(1-{1\over p}\Big)^\varrho\sum_{\nu\geqslant 0}{f(\pnu)\over \pnu}+o(1)\bigg\}{x(\log x)^{\varrho-1}\over\Gamma(\varrho)}\cdot\eqdef{faW}$$ 
\par 
Nous pouvons ˆ prŽsent prŽciser le terme d'erreur en fonction de la vitesse de convergence de la sŽrie \eqref{Ha00}. En effet, l'hypothse \eqref{Ha00} implique immŽdiatement, par sommation d'Abel, que $$\sum_{p\leqslant x}{\{\varrho-\re f(p)\}\log p\over p}\leqslant \eta_x\log x\qquad (x\to\infty)\eqdef{etax}$$
pour une fonction convenable $\eta_x$ tendant vers 0 ˆ l'infini. Choisissons $\gb:=\dm\min(1,\varrho)$, $A:=\varrho$, de sorte que $\gpp=\pi$, $\beta=1$, $\gh=2/\min(1,\varrho)-1$.
Posant $$\delta_1:=\ft13\min(1,\varrho)\leqslant \ft23\gb,\eqdef{pb-appli}$$ 
 nous avons alors $\delta_1\gh=\ft23-\ft13\min(1,\varrho)$ et
$$\sum_{x^\varepsilon<p\leqslant y}{\{\varrho-\re f(p)\}^\gh\log p\over p}\ll \varepsilon^{\delta_1\gh}\log y\qquad (x^\varepsilon<y\leqslant x)$$
pour le choix $\varepsilon:=\eta_x^{1/(1+\delta_1\gh)}+1/\sqrt{\log x}$. Comme la condition \eqref{pqr} dŽcoule immŽdiatement de \eqref{Ha00} pour $x$ assez grand, nous pouvons Žnoncer le rŽsultat suivant. 
\Propc{Wieff}{Soient $A>0$, $B>0$, $\varrho>0$, et $f\in\M_0(A,B)$ une \fmu\ complexe telle que $|f(p)|\leqslant \varrho$ pour tout \np\ $p$. Sous l'hypothse \eqref{Ha00} et avec les notations \eqref{etax}, \eqref{pb-appli}, nous avons
$$M(x;f)={\e^{-\gamma\varrho}x\over\Gamma(\varrho)\log x}\bigg\{\prod_{p}\sum_{\pnu\leqslant x}{f(\pnu)\over \pnu}+O\bigg(\eta_x^a\e^{Z(x;f)}+{\e^{Z(x;f)}\over (\log x)^{b}}\bigg)\bigg\}\eqdef{faW+}$$
o l'on a posŽ $a:=w_f\min(1,\varrho)/\{5-\min(1,\varrho)\}$, $b:=w_f\min(1,\varrho)/6$.}
\par \goodbreak
\rem Sous la condition \eqref{hypgen}, nous dŽduisons de \eqref{faW+} que l'on peut remplacer le terme d'erreur de \eqref{faW} par $O\big(\eta_x^a+1/(\log x)^{b}\big)$. 
\par
\medskip Lorsque $f$ est rŽelle et $\varrho=1$, nous avons donc $a=\ft14$, $b=\ft16$. Cela prŽcise un rŽsultat, mentionnŽ plus haut, de Indlekofer~\citer{In09},  qui obtient dans ce cas  $a=\ft1{36}$ par une mŽthode reposant sur une technique de convolution des \fas.
\medskip
Dans le mme esprit, nous pouvons Žnoncer le rŽsultat suivant. Nous notons $P^+(n)$ le plus grand facteur premier d'un entier $n\geqslant 1$ avec la convention $P^+(1)=1$ et rappelons la notation $\delta_0(\gb)$ dŽfinie en \eqref{delta0}.
\par \vskip-3mm\goodbreak
\Propc{fgbar}{Soient $x\geqslant 2$, $r$ une fonction multiplicative positive ou nulle satisfaisant aux hypothses du \ref{fr} et $f$, $g$ deux \fmus\ telles que \hbox{$|f|^2\leqslant r,\,|g|^2\leqslant r$}. Supposons que, pour une fonction $\eta_x$ tendant vers 0 lorsque $x\to\infty$ et telle que $\eta_x\sqrt{\log x}\gg1$, les majorations $$\sum_{p\leqslant x}{r(p)-h(p)\over p}\leqslant\ft38\beta_0\gb\log \Big({1\over \eta_x}\Big),\qquad  \sum_{p\leqslant x}{r(p)-h(p)\over p}\log p\leqslant \eta_x\log x,\eqdef{rh}$$
o $\beta_0$ est dŽfini en \eqref{delta0},
aient lieu  pour $h=|f|^2$, $h=|g|^2$ et $h=\re f\ov g$. Nous avons alors
$$M(x;|f-g|^2)=M(x;r)\Bigg\{{\sum_{P^+(n)\leqslant x}{|f(n)-g(n)|^2/n}\over \sum_{P^+(n)\leqslant x}r(n)/n}+O\big(\eta_x^a\big)\Bigg\}\eqdef{f-g2}$$
avec $a:=\ft12\beta_0/\{3+(1-\gb)\beta_0\}$.}
\par \goodbreak
Pour Žtablir cette assertion, il suffit d'observer que les fonctions $r$ et $f\ov g$ (respec\-tive\-ment~$|f|^2$,~$|g|^2$) satisfont aux hypothses du \ref{fr} pour le choix $\varepsilon:=\eta_x^{1/(1+\delta_1\gh)}$, \hbox{$\delta_1:=\delta_0(\gb)$}, $\gh:=(1-\gb)/\gb$, $\delta:=\dm\delta_1$. Nous appliquons ensuite ce rŽsultat aux couples $\big(r,f\ov g\big)$ et $\big(r,\ov fg\big)$ (respectivement $(r,|f|^2)$, $(r,|g|^2)$) en nŽgligeant la contribution impliquant le paramtre~$\gc$. \par 
Le cas $r=\1$ du \ref{fgbar} relve de la thŽorie des fonctions ÇsimulatricesÈ ({\it pretentious} en anglais) telle que dŽveloppŽe depuis quelques annŽes par Granville, Soundararajan et d'autres auteurs --- voir par exemple \citer{BGS13}, \citer{GS08}, \citer{Ko13}. Si $f$ et $g$ sont des \fmus\ ˆ valeurs dans le disque unitŽ et si \eqref{rh} est satisfaite avec $r=\1$ pour $h=|f|^2$, $|g|^2$ et $\re f\ov g$,\note{Cette hypothse est en particulier impliquŽe, pour une valeur convenable de $\eta_x$, par la convergence des trois sŽries $\sum_{p}\{1-\re h(p)\}/p$.} il rŽsulte en particulier de \eqref{f-g2} avec $\gb=\ft14$, $A=\dm$, que
$$M\big(x;|f-g|^2\big)={\e^{-\gamma} x\over \log x}\bigg\{\sum_{P^+(n)\leqslant x}{|f(n)-g(n)|^2\over n}+O\big(\eta_x^{2/15}\big)\bigg\}.\eqdef{Mf-g}$$
Le terme principal de \eqref{Mf-g} vaut alors $$M(x;|f|^2)+M(x;|g|^2)-x\e^{-{\bb D}^2(x;f,g)+O(1)}\qquad (x\to\infty),$$
o 
${\bb D}^2(x;f,g):=\re Z(x;1-f\ov g)$ est la pseudo-norme de la thŽorie des fonctions simulatrices. 
\par \medskip
Au titre d'une autre illustration du \ref{fr}, nous pouvons prŽciser un thŽorme de \Halasz\ \citer{Ha71}, \citer{Ha72}, relatif aux lois locales de la rŽpartition des facteurs premiers d'un entier dans un ensemble quelconque. ƒtant donnŽ  un ensemble $E$ de \nps, notons $\Omega(n;E)$ le nombre des facteurs premiers appartenant ˆ $E$, comptŽs avec multiplicitŽ, d'un entier $n$,  et posons \vskip-2mm 
$$E(x):=\sum_{p\leqslant x,\,p\in E}{1\over p}\cdot$$ 
\Halasz\ a montrŽ que, pour tout $\kappa\in]0,1[$ fixŽ et uniformŽment pour $$0\leqslant m\leqslant (2-\kappa)E(x),\eqdef{varm0}$$ le nombre $N_m(x;E)$ des entiers $n\leqslant x$ tels que $\Omega(n;E)=m$ vŽrifie
$$N_m(x;E)=x\e^{-E(x)}{E(x)^m\over m!}\Big\{1+O\Big({|m-E(x)|\over E(x)}+{1\over \sqrt{E(x)}}\Big)\Big\}.\eqdef{HaNm}$$
\par \goodbreak\goodbreak
\break
\Sarkozy\ \citer{Sa77} a ensuite Žtabli que le terme principal de \eqref{HaNm} fournit en fait l'ordre de grandeur du membre de gauche dans  l'intervalle $$\kappa E(x)\leqslant m\leqslant (2-\kappa)E(x).\eqdef{varm}$$
 Autrement dit, sous la contrainte 
\eqref{varm}, nous avons, pour $x$ assez grand,
$$N_m(x;E)\asymp x\e^{-E(x)}{E(x)^m\over m!}\eqdef{ogSa}$$
 ds que $E(x)$ tend vers l'infini avec $x$.\note{Le cas o $E(x)$ est bornŽ relve de techniques de crible ŽlŽmentaires.} Cet encadrement a ŽtŽ ultŽrieurement prŽcisŽ et gŽnŽralisŽ par Balazard \citer{Ba89}.
\goodbreak\par \goodbreak
Soit $\kappa\in]0,1[$. Un argument de convolution standard fournit 
\vskip-5mm
$$S(x;r,E):=\sum_{n\leqslant x}r^{\Omega(n;E)}=x\e^{(r-1)E(x)}\bigg\{1+O\bigg(|r-1|+{1\over (\log x)^{1/2}}\bigg)\Bigg\},\eqdef{SxrE}$$
uniformŽment pour $\kappa\leqslant r\leqslant 2-\kappa$.\note{Une estimation essentiellement Žquivalente dŽcoule d'ailleurs du \ref{cascv} appliquŽ 	aux fonctions $f(n):=r^{\Omega(n;E)}$, $r(n):=\{\max(1,r)\}^{\Omega(n)}$. On peut en effet supposer $\sigma:=|r-1|\leqslant 1-\kappa$ arbitrairement petit et choisir alors  $\gb:=\dm$, $A:=1+\sigma$, $\varrho:=\max(1,r)$, $\gh:=1$, $\varepsilon:=\sigma^4/\{1+\sigma^4(\log x)^{5\sigma}\}$, et $\delta:=\log (1/\sigma)/\log (1/\varepsilon)\leqslant \ft14$.} De plus, il rŽsulte par exemple du thŽorme 1.1 de \citer{Te16} que, dans les mmes conditions, 
$$S(x;r,E)\asymp x\e^{(r-1)E(x)}.\eqdef{SxrE-og}$$\par \goodbreak
Nous obtenons le rŽsultat suivant, qui prŽcise \eqref{HaNm} lorsque $m-E(x)\gg \sqrt{E(x)}$ et implique \eqref{ogSa}.  Notant $\N(E)$ l'ensemble des entiers dont tous les facteurs premiers sont dans $E$, nous posons, lorsque $E$ est fini,
$$F(z;E):=\sum_{n\in\N(E)}{z^{\Omega(n)}\over n}=\prod_{p\in E}\Big(1-{z\over p}\Big)^{-1}\qquad (z\in\CC,\,|z|<2),\eqdef{FzE}$$ et
$\gt(x;E):=\sqrt{\{\log E(x)\}/E(x)}.$
\goodbreak
\Propc{OmE}{Soient $E\subset [1,x]$ un ensemble de \nps\ tel que $\lim_{x\to\infty}E(x)=\infty$, $\kappa\in]0,1[$, et $K>0$.\par 
{\rm(i)} Sous l'hypothse $\kappa\leqslant r\leqslant 2-\kappa$,   $z=r\e^{i\vartheta}$, $-\pi\leqslant \vartheta\leqslant \pi$,  nous avons, uniformŽment pour $x$ assez grand,
$$S(x;z,E)\ll x\,\e^{(r-1-\kappa\vartheta^2/180)E(x)}.\eqdef{majSxzE}$$
\par 
{\rm(ii)} Sous l'hypothse $\kappa\leqslant r\leqslant 2-\kappa$,   $z=r\e^{i\vartheta}$, $|\vartheta|\leqslant K\gt(x;E)$,  nous avons, uniformŽment pour $x$ assez grand,
$$S(x;z,E)=S(x;r,E)\bigg\{{F(z;E)\over F(r;E)}+O\bigg(|\vartheta|\e^{-c\vartheta^2E(x)}+{1\over (\log x)^c}\bigg)\bigg\},\eqdef{SzE/SrE}$$
o $c=c(\kappa,K)>0$.\par 
{\rm(iii)} Nous avons, uniformŽment pour $x$ assez grand,
$$N_m(x;E)\ll x\e^{-E(x)}{E(x)^m\over m!}\qquad (0\leqslant m\leqslant (2-\kappa)E(x)).\eqdef{majunivNm}$$
De plus, sous la condition $\kappa E(x)\leqslant m\leqslant (2-\kappa)E(x)$ et avec $r:=m/E(x)$, nous avons
$$N_m(x;E)=S(x;r,E) {E(x)^m\over m! \e^{m}}\bigg\{1+O\bigg({1\over \sqrt{E(x)}}\bigg)\bigg\}.\eqdef{NmxE}$$}
Nous donnons la dŽmonstration au paragraphe \ref{prgdemOmE}. Les assertions  (i) et (ii) impliquent immŽdiatement  \eqref{NmxE} en appliquant la formule de Cauchy sur le cercle $|z|=r:=m/E(x)$ et en observant que $F(z;E)/F(r;E)=\e^{(z-r)E(x)}\{1+O(\vartheta)\}$. Compte tenu de \eqref{majunivNm}, on obtient~\eqref{HaNm}  en  Žvaluant $S(x;r,E)$ par \eqref{SxrE}.  Enfin, l'estimation \eqref{ogSa} rŽsulte immŽdia\-tement de \eqref{SxrE-og} et~\eqref{NmxE}.
\par \goodbreak
 \Halasz\ a annoncŽ la possibilitŽ de prouver \eqref{NmxE} --- voir \citer{El80}, \page312 --- en utilisant une variante du thŽorme 3 de \citer{Ha71} dans laquelle les hypothses incluent une condition de type~\eqref{cdpr2}, mais valable pour tous les \nps\ n'excŽdant pas $x$. Notre approche ne nŽcessite qu'une majoration en moyenne de type \eqref{cdpr1}.
\par 
Notons encore que le \ref{OmE} demeure valable, {\it mutatis mutandis}, lorsque l'on remplace $\Omega(n;E)$ par la fonction fortement additive $\omega(n;E)=\sum_{p|n,\,p\in E}1$. La borne supŽrieure $2-\kappa$ peut alors tre remplacŽe par $1/\kappa$.\par  De mme, l'extension ˆ l'approximation des lois locales conjointes des fonctions $\Omega(n;E_j)$ ou $\omega(n;E_j)$ $(1\leqslant j\leqslant k)$ relatives ˆ des ensembles fixŽs de \nps\ disjoints  tels que $\min_j E_j(x)\to\infty$ est immŽdiate. En appliquant ˆ la fonction $f(n):=\prod_{1\leqslant j\leqslant  k}z_j^{\omega(n;E_j)}$ le traitement dŽtaillŽ au paragraphe \ref{prgdemOmE}, nous obtenons ainsi
$$\sum_{\di{n\leqslant x}{\omega(n;E_j)=m_j\, (1\leqslant j\leqslant k)}}1=\Bigg\{1+O\Bigg(\sum_{1\leqslant j\leqslant k}{1\over \sqrt{E_j(x)}}\Bigg)\Bigg\}\prod_{1\leqslant j\leqslant k}{E_j(x)^{m_j}\over m_j!\,\e^{m_j}}\sum_{n\leqslant x}\prod_{1\leqslant j\leqslant k}r_j^{\omega(n;E_j)}\eqdef{omEj}$$
uniformŽment pour $\kappa\leqslant r_j:=m_j/E_j(x)\leqslant 1/\kappa$ $(1\leqslant j\leqslant k)$. De plus, le membre de gauche de \eqref{omEj} est Žgal ˆ
$$x\e^{O(\gT)}\prod_{1\leqslant j\leqslant k}{E_j(x)^{m_j}\over m_j!\,\e^{E_j(x)}}$$
uniformŽment pour $\max_jr_j\leqslant 1/\kappa$, avec $\gT:=\sum_{1\leqslant j\leqslant k}\big\{|r_j-1|+1/\sqrt{E_j(x)}\,\big\}$. Ces estimations prŽcisent le thŽorme 1.3 de l'article \citer{Ma17b}, Žgalement mis en ligne  postŽrieurement ˆ la complŽtion et ˆ la diffusion du prŽsent travail.
\medskip\goodbreak
Nos deux dernires applications concernent la rŽpartition des fonctions additives rŽelles relativement ˆ des mesures pondŽrŽes sur l'ensemble des entiers n'excŽdant pas $x$. Nous nous restreignons ici au cas standard d'une loi limite gaussienne, mais nos rŽsultats sont susceptibles de fournir des estimations analogues dans des situations significativement plus gŽnŽrales. 
\par \goodbreak
 ƒtant donnŽes une fonction positive ou nulle $r\in\M_0(x;A,B)$ et une \fad\ rŽelle~$h$, nous notons $z\mapsto F_x(z;h,r)$ la fonction de rŽpartition de la variable alŽatoire $h(n)$ sur l'ensemble des entiers n'excŽdant pas $x$, ŽquipŽ de  la mesure associant ˆ chaque entier~$n$ le poids $r(n)/M(x;r)$, autrement dit
$$F_x(z;h,r):={1\over M(x;r)}\sum_{\di{n\leqslant x}{h(n)\leqslant z}}r(n).$$
\vskip-5mm
Nous posons
$$E_h(x;r):=\sum_{p\leqslant x}{r(p)h(p)\over p},\quad D_h(x;r)^2:=\sum_{p\leqslant x}{r(p)h(p)^2\over p},$$
 et notons
$\Phi(z):=\big(1/\sqrt{2\pi}\big)\int_{-\infty}^z\e^{-u^2/2}\d u$
la fonction de rŽpartition de la loi normale.\par 
\vskip-2mm
\Propc{repad}{Soient $A$, $B$,  des constantes positives, $x\geqslant 2$, $r\in\M_0(x;A,B)$, et $h$ une \fad\ rŽelle. Supposons que:\smallskip \qquad {\rm(i)} \quad$\dsp\min_{\exp\sqrt{\log x}<p\leqslant x}r(p)\gg1$
 \ ;\qquad {\rm(ii)} \quad $D_h(x;r)\gg1$;\par \smallskip
\qquad {\rm(iii)}\quad $\dsp\max_{p\leqslant x}{|h(p)|\over D_h(x;r)}\leqslant \mu_x\leqslant 1$;
\qquad {\rm(iv)}\quad $\dsp\sum_{\pnu\leqslant x}\sum_{\nu\geqslant 2}{r(\pnu)|h(\pnu)|\log \pnu\over \pnu}\ll1.$
\par 
Alors
$$F_x\Big(E_h(x;r)+zD_h(x;r);h,r\Big)=\Phi(z)+O\bigg(\mu_x+{1\over D_h(x;r)}\bigg).\eqdef{apprep}$$
De plus, si l'hypothse (i) est remplacŽe par  \eqref{moyrp} avec $\varepsilon_1:=1/(\log x)^{1/4}$ et $\gb>0$ arbitraire, l'estimation \eqref{apprep} persiste ˆ condition de remplacer $\mu_x$ par $\mu_x\sqrt{\log (1+1/\mu_x)}$.}
\par \medskip
\goodbreak
Le thŽorme 20.1 de \citer{El80} fournit, dans le cas $r=\1$, une Žvaluation de mme type que \eqref{apprep} pour une combinaison linŽaire finie de fonctions $h_j(n+a_j)$ o les $h_j$ $(1\leqslant j\leqslant k)$ sont fortement additives et les $a_j$ des entiers fixŽs. Pour $k=1$, le terme d'erreur de \eqref{apprep} est un peu plus prŽcis que celui de \citer{El80} lorsque, par exemple, $z$ est bornŽ et $D_h(x;r)\mu_x\gg1$. Cette dernire condition est certainement remplie ds que $\max_{p\leqslant x}|h(p)|\gg1$.
\medskip\goodbreak
Notons $\varphi$ la fonction indicatrice d'Euler et dŽsignons par $\Omega(n)$ le nombre total des facteurs premiers d'un entier naturel $n$, comptŽs avec multiplicitŽ. Comme consŽquence spŽcifique du rŽsultat prŽcŽdent, nous pouvons dŽduire facilement une extension d'un thŽorme d'Erd\H os \& Pomerance \citer{EP85} concernant la rŽpartition des nombres $\Omega(\varphi(n))$ et obtenir un terme d'erreur identique ˆ celui de Balazard \& Smati \citer{BS89}. Pour la simplicitŽ de l'ŽnoncŽ, nous restreignons plus qu'il n'est nŽcessaire les hypothses concernant la fonction pondŽrale $r$. 
\par \goodbreak
\Propc{Omphi}{Soient $A$, $B$ des constantes positives, $x\geqslant 2$, $\varrho:=\varrho_x>0$, $\varepsilon_1:=1/(\log x)^{1/4}$, $r\in\M_0(x;A,B)$. Supposons que :\par \quad {\rm(i)}
\ $\dsp\min_{\exp\sqrt{\log x}<p\leqslant x}r(p)\gg1$
  ;
\quad {\rm(ii)} $\dsp\sum_{p\leqslant y}{r(p)\log p\over p}=\varrho\log y+O\big(\varepsilon_1\log y\big)$ \quad $(\e^{1/\varepsilon_1}\leqslant y\leqslant x)$.\par 
Alors, notant $h:=\Omega\circ\varphi$, nous avons
$$F_x\bigg( \dm\varrho (\log_2x)^2+{z\varrho(\log_2x)^{3/2}\over \sqrt{3}};h,r\bigg)=\Phi(z)+O\bigg({1\over \sqrt{\log_2x}}\bigg).\eqdef{repOmphi}$$
\par 
De plus, si l'hypothse (i) est remplacŽe par  \eqref{moyrp}  avec $\varepsilon_1:=1/(\log x)^{1/4}$ et $\gb>0$ arbitraire, l'estimation \eqref{repOmphi} persiste ˆ condition de multiplier le terme d'erreur par $\sqrt{\log_3x}$.
}\goodbreak
\medskip
Le terme d'erreur de \eqref{repOmphi} pose un intŽressant problme ouvert: en accord avec l'estimation de concentration obtenue par Marie-Jeanne \& Tenenbaum \citer{MJT98}, on attend $\ll1/(\log_2x)^{3/2}$, une majoration qui demeure pour l'instant hors d'atteinte des techniques disponibles. 
\bigskip\medskip
\paraunn{Preuve du \ref{casdiv}}\par \smallskip
\paradeuxn{Lemme de Gallagher}
Un rŽsultat de Gallagher \citer{Ga70} (th.\thinspace1), tel qu'ŽnoncŽ, par exemple, au lemme III.4.9 de \citer{Te15} fournit une majoration gŽnŽrique pour la norme
quadratique d'un polyn™me de Dirichlet. Une inŽgalitŽ bien connue de Montgomery et Vaughan (\citer{MV74}, cor.$\thinspace$2) en constitue une forme plus prŽcise. Le rŽsultat suivant est une consŽquence immŽdiate de la majoration initiale de Gallagher. Le lemme 2.1 de \citer{BT12} Žtend le rŽsultat de Gallagher dans une autre direction.
\par 
Nous notons $\e(x):=\e^{2\pi i x}$ $(x\in\r)$.
\Propl{Ga+}{Soient $N\in\N^*$, $\{\lambda_n\}_{n=1}^N$ une suite finie de nombres
rŽels distincts.  Pour tous $\{a_n\}_{n=1}^N\in\CC^N$, $T>0$, nous avons
$$ \int_{-T}^T\Vbs{\sum_{1\leqslant n\leqslant N}a_n\e(\lambda_n t)}^2\dt\ll T\sum_{1\leqslant n\leqslant
N}|a_n|^2\sum_{|\lambda_m-\lambda_n|\leq 1/T}1,\eqdef{majL2Pol}$$
o la constante implicite est absolue.}\medskip
\goodbreak
\rem La forme usuelle sous laquelle  est utilisŽe la majoration de Gallagher est (\cf\ \citer{Te15}, lemme III.4.9)
$$\int_{-T}^T\Vbs{\sum_{1\leqslant n\leqslant N}a_n\e(\lambda_n t)}^2\dt\ll\sum_{1\leqslant n\leqslant
N}|a_n|^2\Big\{T+{1\over \delta_n}\Big\}\eqdef{Gal}$$
o l'on a posŽ $\delta_n:=\min_{m\neq n}|\lambda_m-\lambda_n|$. En pratique, la majoration \eqref{majL2Pol} est souvent plus prŽcise. C'est notamment le cas lorsque $\delta_n\ll \delta_m$ pour $|\lambda_m-\lambda_n|\leqslant 1/T$.
\medskip\goodbreak
\nid
Pour Žtablir \eqref{majL2Pol}, nous observons que, posant
$$A(x):=T\sum_{\di{1\leqslant n\leqslant N}{|x-\lambda_n|\leqslant 1/4T}}a_n,\qquad  S(t):=\sum_{1\leqslant n\leqslant
N}a_n\e(\lambda_n t),$$ nous avons
$$\wh A(t):=\int_\r A(x)\e(-tx)\d x=S(-t){\sin(\pi t/2T)\over \pi t/2T}\cdot $$
D'aprs la formule de Plancherel, nous pouvons donc Žcrire
$$\int_{-T}^T|S(t)|^2\d t\ll \int_\r\Vbs{S(-t){\sin(\pi t/2T)\over \pi t/2T}}^2\d t=\int_\r|A(x)|^2\d x. $$
Posons alors $N_k:=\sum_{|2T\lambda_n-k|\leqslant 1}|a_n|$ $(k\in\z)$ et observons que $|A(x)|\leqslant TN_k$ lorsque $|x-k/2T|\leqslant 1/4T$. Il suit
$$\int_\r|A(x)|^2\d x\ll T\sum_{k\in\z}N_k^2.$$
D'aprs l'inŽgalitŽ de Cauchy-Schwarz, nous avons
$$N_k^2\leqslant \sum_{|2T\lambda_n-k|\leqslant 1}|a_n|^2\sum_{|2T\lambda_m-k|\leqslant 1}1\leqslant\sum_{|2T\lambda_n-k|\leqslant 1}|a_n|^2\sum_{|\lambda_m-\lambda_n|\leqslant 1/T}1. $$
En sommant sur $k$, nous obtenons bien l'inŽgalitŽ annoncŽe.
\qed
\bigskip\goodbreak
\paradeuxn{RŽduction au cas exponentiellement multiplicatif}
Soit $g$ la fonction exponentiellement multiplicative co•ncidant avec $f$ sur l'ensemble des nombres premiers n'excŽdant pas $x$ et nulle sur les \nps\ $>x$, autrement dit
$$g(\pnu):=\normalbaselineskip=18pt\cases{f(p)^\nu/\nu!& si $p\leqslant x$, $\nu\geqslant 1$,\cr0 &si $p>x$.\cr}\eqdef{defg}$$
Nous avons donc, avec la notation \eqref{MphiH},
$$G(s):=\sum_{n\geqslant 1}{g(n)\over n^s}=\e^{v_f(s;x)}.\eqdef{defG}$$
\par 
Supposons la majoration \eqref{majeff} acquise pour la fonction $g$. Nous allons montrer qu'elle vaut encore pour $f$.
\par 
Nous avons $f=g*h$
avec
$$h(p^\nu)=\sum_{j+k=\nu}(-1)^j{f(p)^j\over j!}f(p^k)\qquad (p\leqslant x,\,\nu\geqslant 1), \eqdef{defh}$$ 
de sorte que $h(p)=0$ pour tout $p$, $h(p^\nu)=0$ si $p>x$, et
$$\eqalign{\sum_{p,\,\nu\geqslant 2}{|h(p^\nu)|\over p^\nu}&\leqslant \sum_{p}\sum_{k\geqslant 0}{{|f(p^k)|\over
p^k}}\sum_{j\ge\max(0, 2-k)}{|f(p)|^j\over j!p^j}
\cr&\ll\sum_{p}{|f(p)|^2\over
p^2}+\sum_{p,\,k\geqslant 2}{|f(p^k)|\over p^k}\ll 1.\cr}\eqdef{majShpnu}
$$
Cela implique
$$\leqalignno{&\prod_{p}\sum_{\nu\geqslant 0}
{|h(p^\nu)|\over p^\nu}\ll 1, &\eqdef{prodh1}\cr
&\prod_{p}\sum_{\nu\geqslant 0}
{h(p^\nu)\over p^\nu}=\prod_{p\leqslant x}\e^{-f(p)/p}\sum_{\nu\geqslant 0}{f(p^\nu)\over p^\nu}\cdot
&\eqdef{prodh2}\cr}$$
\par
 
De plus, pour $2\leqslant y\leqslant x$, nous pouvons Žcrire d'une part
$$\eqalign{&\sum_{y<n\leqslant x}{|h(n)|\over n}\leqslant {1\over (\log y)^2}\sum_{n\leqslant x}{|h(n)|(\log n)^2\over n}\cr
&\qquad \le{1\over (\log y)^2}\bigg\{\sum_{mp^\nu\leqslant x}{|h(m)h(p^\nu)|(\log p^\nu)^2\over mp^\nu}+\sum_{\di{mp^\nu q^\mu\leqslant x}{p\neq q}}{|h(m)h(p^\nu)h(q^\mu)|(\log p^\nu)(\log q^\mu)\over mp^\nu q^\mu}\bigg\}
\cr&\qquad \ll
{1\over (\log y)^2}\bigg\{\sum_{p^\nu\leqslant x}{|h(p^\nu)|(\log p^\nu)^2\over p^\nu}+\bigg(\sum_{p^\nu\leqslant x}{|h(p^\nu)|\log p^\nu\over p^\nu}\bigg)^2\bigg\},\cr } $$
et, d'autre part,
$$\eqalign{\sum_{p^\nu\leqslant x}{|h(p^\nu)|(\log p^\nu)^2\over p^\nu}&\leqslant \sum_{\di{p^{j+k}\leqslant x}{j+k\ge
2}}{|f(p)|^j|f(p^k)|(k+j)^2(\log p)^2\over j!p^{k+j}} \cr
&\le\sum_{p^k\leqslant x}{|f(p^k)|(\log p^k)^2\over p^k}\sum_{j\geqslant \max(0,2-k)}{(j+1)^2|f(p)|^j\over j!p^j}\ll
1,
\cr}$$
d'aprs \eqref{C1} avec $j=1$. Donc
$$Q(y):= \sum_{n>y}{|h(n)|\over n}\ll {1\over (\log y)^2}\qquad (2\leqslant y\leqslant x).\eqdef{majQ}$$
\par \goodbreak
Cela dit, nous avons
$$M(x;f)=\sum_{n\leqslant x}h(n)M\Big({x\over n};g\Big). $$
Il suit
$$\eqalign{&M(x;f)\cr&\ll \sum_{n\leqslant x}{x|h(n)|\over n\log (2x/n)}\bigg\{\int_{1/\log (3x/n)}^1{H_T(\alpha)\over \alpha}\d\alpha+{\e^{Z(x;r)}\over\min(\sqrt{T}, \log (2x/n))}+{\e^{Z(x/n;r)}\log_2(3x/n)\over T}\bigg\}\cr&\ll x\int_{1/\log(3x)}^1{H_T(\alpha)\over \alpha}\Theta_1(3x\e^{-1/\alpha})\d\alpha+x\e^{Z(x;r)}\Big({\Theta_1(x)\over \sqrt{T}}+\Theta_2(x)\Big)+{x\e^{Z(x;r)}\Theta_1(x)\log_2x\over T},\cr}\eqdef{majconv}$$
o l'on a posŽ
$$\Theta_j(y):=\sum_{n\leqslant y}{|h(n)|\over n\{\log (2x/n)\}^j}\quad(j=1,2,\quad y\geqslant 2).$$
\goodbreak\noi
En scindant la somme ˆ $\sqrt{x}$, nous dŽduisons de \eqref{majQ} que
$$\Theta_j(y)\ll{1\over (\log x)^j}\qquad (2\leqslant y\leqslant x).$$
Cela implique bien l'estimation requise.
\par 
Si nous adjoignons l'hypothse \eqref{hypsupdiv}, l'assertion relative au cas exponentiellement multiplicatif nous permet d'omettre dans \eqref{majconv} le dernier terme de l'accolade.
\goodbreak\par\bigskip
\paradeuxn{Preuve  dans le cas exponentiellement multiplicatif }
Soit $g$ la fonction multiplicative dŽfinie par \eqref{defg}. La premire Žtape consiste ˆ majorer 
$$K(x):=\sum_{n\leqslant x}g(n)\log n$$
en fonction d'une moyenne sur $[1,x]$ de $t\mapsto M(t;g)$. \par 
\Propl{eqfoncK}{Soient $A>0,\, B>0$. UniformŽment pour $x\geqslant 2$, $f,r\in\M_0(x;A,B)$, $|f|\leqslant r$, nous avons
$$|K(x)|\leqslant Ax\int_{1}^x\big|M(t;g)\big|{\dd t\over t^2}+O\Big({x\e^{Z(x;r)}\over \log x}\Big).\eqdef{majfoncK}$$}
\nid Nous avons 
$$\eqalign{K(x)=\sum_{n\leqslant x}g(n)\sum_{\pnu|n}\log p=\sum_{\pnu\leqslant x}(\log p)\sum_{m\leqslant x/\pnu}g(m\pnu).\cr}\eqdef{decK}$$
La somme intŽrieure vaut
$$\sum_{\di{m\leqslant x/\pnu}{p\,\nmid\, m}}g(\pnu)g(m)+\sum_{m\leqslant x/p^{\nu+1}}g(mp^{\nu+1})=g(\pnu)M\Big({x\over \pnu};g\Big)-R_\nu(p)$$
avec
$$\eqalign{R_\nu(p):=&\sum_{m\leqslant x/p^{\nu+1}}\big\{g(\pnu)g(mp)-g(mp^{\nu+1})\big\}\cr=&\ \sum_{k\geqslant 0}\big\{g(\pnu)g(p^{k+1})-g(p^{\nu+k+1})\big\}M\Big({x\over p^{\nu+k+1}};g_p\Big)\cr=&\ \sum_{k\geqslant 0}{1\over (k+\nu+1)!}\Big\{{\nu+k+1\choose k}-1\Big\}g(p)^{\nu+k+1}M\Big({x\over p^{\nu+k+1}};g_p\Big),\cr}$$
o $g_p$ dŽsigne la fonction multiplicative co•ncidant avec $g$ sur l'ensemble des entiers premiers ˆ $p$ et nulle sur l'ensemble des multiples de $p$.\par \goodbreak Compte tenu de la majoration de Halberstam--Richert \citer{HR79}\note{Voir le th. III.3.5 de \citer{Te15} pour une version simplifiŽe suffisante ici.} 
$$M(y;|g|)\ll{y\,\e^{Z(y;r)}\over \log y}\qquad (2\leqslant y\leqslant x), \eqdef{triv1}$$
nous pouvons Žcrire
$$\eqalign{\sum_{\pnu\leqslant x}R_\nu(p)\log p&\ll\sum_{j\geqslant 2}\sum_{p^j\leqslant x}{(1+|g(p)|)^j\log p\over j!}M\Big({x\over p^j};g_p\Big)\cr&\ll x\e^{Z(x;r)}\sum_{j\geqslant 2}\sum_{p^j\leqslant x}{(1+|g(p)|)^j\log p\over j!\,p^j\log (2x/p^j)}\cr&\ll{x \e^{Z(x;r) }\over \log x}+x\e^{Z(x;r)}\sum_{p> \sqrt{x}}{(\log p)^2\over p^{3/2}x^{1/4}} \ll{x \e^{Z(x;r) }\over \log x}\cdot\cr}$$
\par \goodbreak
En reportant dans \eqref{decK}, nous obtenons donc
$$K(x)=\sum_{d\leqslant x}\Lambda(d)g(d)M\Big({x\over d};g\Big)+O\Big({x\e^{Z(x;r)}\over \log x}\Big).\eqdef{eqK1}$$
\par 
La contribution au terme principal des entiers $d=\pnu$ avec $\nu\geqslant 2$ n'excde pas
$$\ll x\e^{Z(x;r)}\sum_{\di{\pnu\leqslant x}{\nu\geqslant 2}}{|f(p)|^\nu\log p\over \nu!\,p^\nu\log (2x/\pnu)}\ll{x\e^{Z(x;r)}\over \log x},$$
o nous avons de nouveau fait appel ˆ \eqref{triv1} et estimŽ la somme en $\pnu$ en la scindant ˆ $\sqrt{x}$. Introduisant la fonction de TchŽbychev $\vartheta(t):=\sum_{p\leqslant t}\log p$ et notant $R(t):=\vartheta(t)-t$, nous pouvons donc Žcrire 
$$\eqalign{|K(x)|&\leqslant A\int_{1}^x\Big|M\Big({x\over t};g\Big)\Big|\d\vartheta(t)+O\Big({x\e^{Z(x;r)}\over \log x}\Big)\cr&=A\int_{1}^x\Big|M\Big({x\over t};g\Big)\Big|\d t+A\int_{1-}^x\Big|M\Big({x\over t};g\Big)\Big|\d R(t)+O\Big({x\e^{Z(x;r)}\over \log x}\Big)\cr
&=Ax\int_{1}^x|M(t;g)\Big|{\d t\over t^2}+A\int_{1-}^x\Big|M\Big({x\over t};g\Big)\Big|\d R(t)+O\Big({x\e^{Z(x;r)}\over \log x}\Big).\cr}$$
\par  La dernire intŽgrale peut tre ŽvaluŽe par sommation d'Abel en notant que $$|\d|M(t;g)||\leqslant \d M(t;|g|).$$
Nous avons
$$\eqalign{\int_{1-}^x\Big|M\Big({x\over t};g\Big)\Big|\d R(t)&=|M(x;g)|-\int_1^xR\Big({x\over t}\Big)\d |M(t;g)|\cr&\ll{x\e^{Z(x;r)}\over \log x}+\sum_{n\leqslant x}{x|g(n)|\over n(\log 2x/n)^2}\cr&\ll{x\e^{Z(x;r)}\over \log x}+\sum_{2^k\leqslant x}{2^k\over k^2+1}\sum_{x/2^{k+1}<n\leqslant x/2^k}|g(n)|\cr&\ll
{x\e^{Z(x;r)}\over \log x}+\sum_{2^k\leqslant x}{x\e^{Z(x;r)}\over (k^2+1)\log (2x/2^k)}\ll {x\e^{Z(x;r)}\over \log x}\cdot\cr}$$
Cela complte la preuve de \eqref{majfoncK}.\qed\goodbreak
\Propl{LineqMg}{Soient $A>0$, $B>0$. Pour $x\geqslant 2$, \hbox{$f,r\in\M_0(x;A,B)$}, $|f|\leqslant r$, nous avons uniformŽment 
$$\leqalignno{|M(x;g)|&\leqslant {Ax\over \log x}\int_1^x|M(t;g)|{\dd t\over t^2}+O\Big({x\e^{Z(x;r)}\over (\log x)^2}\Big),&\eqdef{ineqMg}\cr
\int_1^x{|M(t;g)|\log t\over t^2}\dt&\leqslant 2\int_1^x{|K(t)|\over t^2}\d t+O(1).&\eqdef{ineqMg1}\cr}$$}
\nid
Nous avons
$$M(x;g)\log x-K(x)=\sum_{n\leqslant x}g(n)\int_n^x{\dd t\over t}=\int_1^xM(t;g){\dd t\over t}\ll\int_1^x{\e^{Z(t;r)}\over \log 2t}\d t\ll{x\e^{Z(x;r)}\over \log x},$$
o l'avant-dernire majoration rŽsulte de \eqref{triv1}. Cela implique \eqref{ineqMg} en repor\-tant dans~\eqref{majfoncK}. De plus, pour $x>\e^2$, nous pouvons Žcrire
$$\eqalign{\int_{\e^2}^x{|M(t;g)|\log t\over t^2}\d t&\leqslant \int_{\e^2}^x{|K(t)|\over t^2}\d t+\int_{\e^2}^x{\d t\over t^2}\int_1^t{|M(u;g)|\over u}\d u\cr&=\int_{\e^2}^x{|K(t)|\over t^2}\d t+\int_1^x{|M(u;g)|\over u}\int_{\max(u,\e^2)}^x{\dd t\over t^2}\d u\cr&\leqslant \int_{\e^2}^x{|K(t)|\over t^2}\d t+\int_{\e^2}^x{|M(u;g)|\over u^2}\d u +O(1). \cr}$$
Cela implique bien \eqref{ineqMg1}.
\qed
\Propl{MinH}{Soient $A>0,\, B>0$. Sous les conditions $x\geqslant 2$, \hbox{$f\in\M_0(x;A,B)$}, $T\geqslant 1$, $1/\log x\leqslant \alpha\leqslant \dm$, nous avons uniformŽment 
$$H_T(\alpha)\gg1.\eqdef{minH}$$}
\nid
Nous avons
$$\int_{-1/2}^{1/2}|v_f(1+\alpha+i\tau)|^2\d\tau\ll\sum_{p}{|g(p)|^2\over p^{2+2\alpha}}\sum_{|\log (p'/p)|\leqslant 1/2}1\ll\sum_{p}{1\over p\log p}\ll1,$$
o la premire majoration rŽsulte du \ref{Ga+}, et la seconde du thŽorme de Brun--Titchmarsh. Il s'ensuit que $\max_{\di{\sigma=1+\alpha}{|\tau|\leqslant 1/2}}|\e^{v_f(s)}|\gg 1.$
\qed
\bigskip\goodbreak
Nous sommes ˆ prŽsent en mesure d'aborder la phase finale de la dŽmonstration. L'argument repose essentiellement sur la technique ŽlaborŽe par Hal‡sz \citer{Ha68}, avec, ainsi qu'il a ŽtŽ mentionnŽ dans l'introduction, certains raffinements issus de \citer{Mo78} et \citer{Te15}.
\par 
Nous observons d'abord que l'inŽgalitŽ de Cauchy--Schwarz fournit
$$\int_1^x{|K(t)|\over t^2}\d t\leqslant \bigg(\int_1^x{|K(t)|^2\over t^3}\d t\log x\bigg)^{1/2}\eqdef{K/t2}$$
et que la formule de Plancherel permet d'Žcrire, pour tout $\alpha>0$, avec la notation \eqref{defG},
$$\int_1^x{|K(t)|^2\over t^{3+2\alpha}}\d t={1\over 2\pi}\int_\r\abs{G'(1+\alpha+i\tau)\over 1+\alpha+i\tau}^2\d\tau.\eqdef{Plan}$$

\par 
La dernire intŽgrale peut tre estimŽe gr‰ce au \ref{Ga+}. En notant que $G'(s)=v_f'(s)\e^{v_f(s)}$ et que $$G(s)\ll\e^{Z(\exp(1/\{\sigma-1\});r)}\qquad (\sigma>1),\eqdef{majtrivG}$$ il suit, dans un premier temps,
$$\eqalign{\int_{-T}^T&\abs{G'(1+\alpha+i\tau)\over 1+\alpha+i\tau}^2\d\tau\ll \sum_{|k|\leqslant T}\int_{-1/2}^{1/2}{\abs{G'(1+\alpha+ik+i\tau)}^2\over k^2+1 }\d\tau\cr&\ll\sum_{|k|\leqslant T}{1\over k^2+1}\max_{\di{\sigma=1+\alpha}{|\tau-k|\leqslant \fs12}}\abs{\e^{2v_f(s)}}\sum_{p}{|g(p)|^2(\log p)^2\over p^{2+2\alpha}}\sum_{|\log (p'/p)|\leqslant 1/2}1.\cr}$$
En vertu de l'inŽgalitŽ de Brun--Titchmarsh,  la somme intŽrieure en $p'$  est $\ll p/\log p$. La somme en $p$ est donc
$$\ll \sum_{p}{\log p\over p^{1+2\alpha}}\ll{1\over \alpha }\cdot$$
D'o
$$\int_{-T}^T\abs{G'(1+\alpha+i\tau)\over 1+\alpha+i\tau}^2\d\tau\ll{H_T(\alpha)^2\over \alpha}\cdot$$
Nous traitons l'intŽgrale complŽmentaire de manire similaire en faisant appel ˆ la majoration triviale \eqref{majtrivG}.
Nous avons
$$\eqalign{\int_{|\tau|>T}&\abs{G'(1+\alpha+i\tau)\over 1+\alpha+i\tau}^2\d\tau\ll \sum_{k\neq0,-1}\int_{kT}^{(k+1)T}{\abs{G'(1+\alpha+i\tau)}^2\over k^2T^2 }\d\tau\cr&\ll\sum_{k\neq0,-1}{\e^{2Z(\exp(1/\alpha);r)}\over k^2T}\sum_{p}{|g(p)|^2(\log p)^2\over p^{2+2\alpha}}\sum_{|\log (p'/p)|\leqslant 1/T}1.\cr}$$
En vertu de l'inŽgalitŽ de Brun--Titchmarsh,  la somme intŽrieure en $p'$  est $$\ll (1+p/T)/\log (2+p/T).$$ La somme en $p$ est donc
$$\ll \sum_{p\leqslant T}{(\log p)^2\over p^2}+\sum_{T<p\leqslant T^2}{(\log p)^2\over Tp}+\sum_{p}{\log p\over Tp^{1+2\alpha}}\ll 1+{1\over \alpha T}\cdot$$
En rassemblant nos estimations nous pouvons donc Žnoncer que
$$\int_{\r}\abs{G'(1+\alpha+i\tau)\over 1+\alpha+i\tau}^2\d\tau\ll {H_T(\alpha)^2\over \alpha}+ {\e^{2Z(\exp(1/\alpha);r)}\over T}\Big(1+{1\over \alpha T}\Big).$$\par
Choisissons $\alpha:=1/\log x$ et reportons successivement dans \eqref{Plan}, \eqref{K/t2} et \eqref{ineqMg1}. Nous obtenons
$$\int_1^x{|M(t;g)|\log t\over t^2}\d t\ll H_T\Big({1\over \log x}\Big)\log x+{\e^{Z(x;r)}\sqrt{\log x}\over \sqrt{T}}+{\e^{Z(x;r)}\log x\over T},$$  
puisque, en vertu de \eqref{minH}, la contribution du terme d'erreur de \eqref{ineqMg1} est dominŽe par le premier terme du majorant.
Nous exploitons cette estimation sous la forme
$$\int_{\sqrt{y}}^y|M(t;g)|{\dd t\over t^2}\ll H_T\Big({1\over \log y}\Big)+{\e^{Z(y;r)}\over \sqrt{T\log y}}+{\e^{Z(y;r)}\over T}\qquad (\e\leqslant y\leqslant x).\eqdef{intMryy}$$
Ainsi
$$\eqalign{\int_\e^x|M(t;g)|{\dd t\over t^2}&\ll\int_\e^x|M(t;g)|{\dd t\over t^2}\int_{\sqrt{t}}^t{\dd y\over y\log y}=\int_{\sqrt{\e}}^{x}{\dd y\over y\log y}\int_y^{\min(x,y^2)}|M(t;g)|{\dd t\over t^2}\cr&\ll\int_{\sqrt{\e}}^{x}H_T\Big({1\over \log y}\Big){\dd y\over y\log y}+{\e^{Z(x;r)}\over \sqrt{T}}+{\e^{Z(x;r)}\log_2x\over T}\cr&\ll\int_{1/\log x}^1{H_T(\alpha)\over \alpha}\d\alpha+{\e^{Z(x;r)}\over \sqrt{T}}+{\e^{Z(x;r)}\log_2x\over T}\cdot\cr}$$
En reportant dans \eqref{ineqMg}, nous obtenons bien \eqref{majeff}.\par 
Sous l'hypothse supplŽmentaire \eqref{hypsupdiv}, le dernier terme de \eqref{intMryy} est $$\ll{\e^{Z(x;r)}\over T}\Big({\log y\over \log x}\Big)^c.$$
En reportant dans les calculs, nous constatons alors que la contribution de cette quantitŽ est bien dominŽe par celle du second terme de la majoration \eqref{intMryy}.
\bigskip
\medskip
\par 
\paraunn{Preuve du \ref{cascv}}

\paradeuxn{Lemmes}

Nous aurons plusieurs fois l'usage du rŽsultat auxiliaire suivant.
\Propl{lem1}{Dans les hypothses du \ref{cascv}, nous avons
$$\e^{Z(y;r)}\ll\e^{Z(x;r)}\Big\{\epsilon^{\varrho}+{\Big({\log y\over\log 2x}\Big)^{\varrho}}\Big\}\qquad
(1\leqslant y\leqslant x).\eqdef{L2}$$}

\nid 
Posons
$$Z_1(t):=\sum_{p\leqslant t}{r(p)\log p\over p}\qquad (t\geqslant 2),\eqdef{Tu}$$
de sorte que, par \eqref{C2}, nous avons, pour une constante positive $c_1$ convenable, $$Z_1(t)\geqslant \big(\varrho-c_1\epsilon\big)\log t \qquad (t> x^{\epsilon}).$$
De plus, 
$$Z(x;r)-Z(y;r)=\int_y^x{Z_1(t)\over t(\log 2t)^2}\d t+O(1)\qquad (2\leqslant y\leqslant x).\eqdef{Ex-Ey}
$$
 Pour
$x^{\epsilon}<y\leqslant x$, nous pouvons donc Žcrire
$$Z(x;r)-Z(y;r)\ge\big(\varrho-c_1\epsilon\big)\log \Big({\log x\over\log y}\Big)+O(1)=\varrho\log\Big({\log x\over \log
y}\Big)+O(1).\eqdef{Ex-Ey/2}$$
 Or, si $y\leqslant x^{\epsilon}$, nous avons
$$\eqalign{Z(x;r)-Z(y;r)&\geqslant Z(x;r)-Z\big(x^{\epsilon};r\big)\ge\varrho\log(1/\epsilon)+O(1).\cr }\eqdef{Ex-Ey/3}
$$ Cela implique bien le rŽsultat annoncŽ.\qed
\medskip\goodbreak
Le rŽsultat qui suit est une consŽquence immŽdiate du thŽorme II.5.2 ou du thŽo\-rme~II.5.4 de
\citer{Te15}.
\Propl{lem2}{Soient $A$, $\gb$ des constantes positives.  Pour $\gb\leqslant \varrho\leqslant A$, $y\geqslant 2$, nous avons uniformŽment
$$T_\varrho(y):=\sum_{n\leqslant y}\tau_\varrho(n)= {y\over \Gamma(\varrho)}(\log y)^{\varrho-1}\Big\{1+O\Big({1\over
\log y}\Big)\Big\}.
\eqdef{Trho}$$}
\par\medskip
Soient $f$ et $r$ des \fmus\ satisfaisant aux hypothses du \ref{cascv} et  $g$ la fonction exponentiellement multiplicative dŽfinie par
$$g(p^\nu)=\normalbaselineskip=15pt\cases{f(p)^\nu/\nu!& si $p\leqslant x$,\cr \varrho^\nu/\nu!& si $p>x$,\cr} \qquad (p\ge
2,\,\nu\geqslant 1).\eqdef{g}$$
Nous posons $$Z(y;g):=\sum_{p\leqslant y}{g(p)\over
p}\qquad (y\geqslant 2).\eqdef{ZZ}$$
Ainsi, $Z(y;g)=Z(y;f)$ tel que dŽfini en \eqref{defZZ} pour $y\leqslant x$. Nous dŽfinissons en
outre
$$\eqalign{\L^*_x&:={\e^{-\gamma\varrho+Z(x;f)}\over (\log x)^\varrho}\qquad (x\geqslant 2), \cr
G(s)&:=\sum_{n\geqslant 1} {g(n)\over n^s}=\exp\bigg\{\sum_{p}{g(p)\over p^s}\bigg\}\qquad (\sigma>1),\cr
g_x(n)&:=g(n)-\L^*_x\tau_\varrho(n)\qquad (n\geqslant 1)}\eqdef{L*G}
$$
de sorte que, par \eqref{Trho},
$$M(x;g_x)=M(x;g)-{\e^{-\gamma\varrho}x\,\e^{Z(x;f)}\over \Gamma(\varrho)\log x}+O\Big({\epsilon^2 x\,\e^{Z(x;f)}\over
\log x}\Big) .\eqdef{Msggx}$$ 
\par \goodbreak
Rappelons encore la notation $\delta:=w_f\delta_1\leqslant \ft13\min(2,\varrho)$.
\Propl{lem3}{Pour $x\geqslant 2$, nous avons
$$M(x;g_x)\ll{x\over \log x}\bigg\{\int_1^x{|M(y;g_x)|\over
y^2}\dy+\epsilon^{\delta}\e^{\re Z(x;f)}
\bigg\}.\eqdef{maj1Mx}
$$ }
\par
\nid Posons $\varrho_0:=\min(1,\varrho)/3$. Observons d'abord qu'il rŽsulte de \eqref{triv1},  \eqref{L2}, et \eqref{Trho},  que
$$\eqalign{M(y;|g_x|)&\ll {y\,\e^{Z(x;r)}\over \log y}\rstt{\varrho}\cr&\ll y\varepsilon^{\min(2,\varrho)}\e^{Z(x;r)}\ll y\varepsilon^{\delta+\varrho_0}\e^{Z(x;f)}\cr}\qquad (2\leqslant y\leqslant x),
\eqdef{triv2}$$
o la seconde estimation rŽsulte d'un calcul d'extremum et la dernire de \eqref{pqr}.
\par\goodbreak
Posons ensuite
$$N_x(y):=\sum_{n\leqslant y}g_x(n)\log n\qquad (2\leqslant y\leqslant x). \eqdef{Nx}$$
Nous avons
$$\eqalign{M(y;g_x)\log y-N_x(y)&=\sum_{n\leqslant y}g_x(n)\log \Big({y\over n}\Big)=\int_1^y{M(z;g_x)\over z}\d z\ll y\varepsilon^{\delta+\varrho_0}\e^{Z(x;f)},\cr}
$$ 
et donc
$$M(y;g_x)\ll{|N_x(y)|\over \log y}+{y\varepsilon^{\delta+\varrho_0}\e^{\re Z(x;f)}\over \log y}\qquad (2\leqslant y\le
x).\eqdef{maj1My}
$$
\par
\goodbreak
Nous allons estimer $N_x(y)$ pour $2\leqslant y\leqslant x$ en utilisant l'identitŽ
$$g_x(ab)=g(a)g_x(b)+\L^*_x \wt g(a)\tau_\varrho(b)\qquad \big((a,b)=1)\big),$$
o nous avons posŽ $\wt g(n):=g(n)-\tau_\varrho(n)$ $(n\geqslant 1)$. Ainsi $$\eqalign{N_x(y)&=\sum_{\di{mp^\nu\leqslant y}{p\,\nmid\,m}}g_x(mp^\nu)\log p^\nu\cr
&=\sum_{p^\nu\leqslant y}g(p^\nu)\log
p^\nu\sum_{\di{m\leqslant y/p^\nu}{p\,\nmid\,m}}g_x(m)+\L^*_x\sum_{p^\nu\leqslant y}\wt g(p^\nu)\log
p^\nu\sum_{\di{m\leqslant y/p^\nu}{p\,\nmid\,m}}\tau_\varrho(m).\cr}\eqdef{R} $$
Ici et dans la suite, une somme portant sur des puissances de nombres premiers $p^\nu$ est entendue comme
une somme double o l'exposant satisfait systŽmatiquement $\nu\geqslant 1$. La contribution au membre de
droite de \eqref{R} des puissances $p^\nu$ avec $\nu\geqslant 2$ est estimŽe en faisant appel
ˆ \eqref{Trho} et \eqref{triv2}. Elle est
$$\eqalign{&\ll y\e^{Z(x;r)}\rstt{\varrho}\sum_{\di{p^\nu\leqslant y}{\nu\ge
2}}{\{|g(p^\nu)|+\tau_\varrho(p^\nu)\}\log p^\nu\over p^\nu\log (2y/p^\nu)}\cr
&\ll {y\e^{Z(x;r)}\over \log y}\,\rstt{\varrho}\ll y\varepsilon^{\delta+\varrho_0}\e^{Z(x;f)},\cr}$$ 
o la seconde estimation dŽcoule, par sommation d'Abel, de l'hypothse de multiplicativitŽ
exponentielle effectuŽe sur $g$.
\par\goodbreak
Nous pouvons donc Žcrire
$$N_x(y)=W_1+W_2+O\Big(y\varepsilon^{\delta+\varrho_0}\e^{Z(x;f)}\Big)\qquad (2\leqslant y\leqslant x),\eqdef{est1Ry}$$
avec $$\eqalign{W_1&:=\sum_{p\leqslant y}g(p)\log
p\sum_{\di{m\leqslant y/p}{p\,\nmid\,m}}g_x(m), \qquad 
W_2:=\L^*_x\sum_{p\leqslant y}\wt g(p)\log
p\sum_{\di{m\leqslant y/p}{p\,\nmid\,m}}\tau_\varrho(m).}$$
 La  somme intŽrieure de $W_1$ vaut
$M(y/p;g_x)-U_p(y/p)$
avec, pour $p\leqslant z\leqslant y$,
$$\eqalign{U_p(z):=\sum_{\di{m\leqslant z}{p|m}}g_x(m)&\ll\sum_{\nu\le(\log z)/\log
p}\bigg\{|g(p^\nu)|M\Big({z\over p^\nu};|g_x|\Big)+\L_x^*\tau_\varrho(p^\nu)T_\varrho\Big({z\over p^\nu}\Big)\bigg\}\cr&\ll {z\,\e^{Z(x;r)}\over
p\log (2z/p)}\Big\{\epsilon^{\varrho}+\Big({\log 2z/p\over \log x}\Big)^{\varrho}\Big\}\ll {z\,\varepsilon^{\delta+\varrho_0}\e^{Z(x;f)}\over p},\cr}$$
o l'on a fait appel ˆ \eqref{triv2} et \eqref{Trho}.\note{Nous avons Žgalement utilisŽ ici
l'inŽgalitŽ ŽlŽmentaire $$\sum_{1\leqslant \nu\leqslant (\log z)/\log p}{z\over p^\nu}\log^{\varrho-1}
\({2z\over p^\nu}\)\ll {z\over p}\log^{\varrho-1} \({2z\over p}\)\qquad (2\leqslant p\leqslant z),$$
valable uniformŽment pour $0\leqslant \varrho\ll1$, et dont nous omettons la
dŽmonstration.} La contribution ˆ
$W_1$ des termes impliquant $U_p(y/p)$ est donc
$$\ll\sum_{p\leqslant \sqrt y}{y\,\varepsilon^{\delta+\varrho_0}\e^{Z(x;f)}\log p\over p^2}\ll y\,\varepsilon^{\delta+\varrho_0}\e^{Z(x;f)}.
$$
\par 
Une variante de \eqref{Trho} prouvŽe par la mme mŽthode d'intŽgration complexe (\cf\ \citer{Te15}, thŽorme
III.5.2), fournit par ailleurs, uniformŽment pour $2\leqslant p\leqslant z$, $2\gb\leqslant \varrho\leqslant A$,
$$\sum_{\di{m\leqslant z}{p\,\nmid\,m}} \tau_\varrho(m)={z\over \Gamma(\varrho)}\Big(1-{1\over
p}\Big)^\varrho(\log 2z)^{\varrho-1}\Big\{1+O\Big({1\over \log 2z}\Big)\Big\}.$$
\goodbreak\noi
On en dŽduit que l'on peut remplacer la somme intŽrieure de $W_2$ par
$${y\over \Gamma(\varrho)p}\Big(\log {2y\over p}\Big)^{\varrho-1}, $$
avec une erreur englobŽe par le terme rŽsiduel de \eqref{est1Ry}. Nous avons donc
obtenu jusqu'ici
$$N_x(y)=S_x(y)+W_2^*+O\Big(y\,\varepsilon^{\delta+\varrho_0}\e^{Z(x;f)}\Big)\qquad (2\leqslant y\leqslant x), \eqdef{est2Ry}$$
o l'on a posŽ
$$S_x(y):=\sum_{p\leqslant y}g(p)(\log
p)M\Big({y\over p};g_x\Big),\qquad W_2^*:={y\L^*_x\over \Gamma(\varrho)}\sum_{p\leqslant y}{\wt g(p)\log
p\over p}\Big(\log {2y\over p}\Big)^{\varrho-1}. \eqdef{Sx}$$
\par
 Nous allons Žtablir
l'estimation
$$W^*_2\ll\epsilon^{\delta}y\,\e^{Z(x;f)}\Big\{\varepsilon^{\varrho_0}+\Big({\log y\over \log x}\Big)^{\varrho_0}\Big\}\qquad
(2\leqslant y\leqslant x).\eqdef{W*2}
$$
 Pour cela, nous observons d'abord que la majoration triviale $g(p)\ll1$
fournit
$$W^*_2\ll y\e^{Z(x;f)}\Big({\log y\over \log x}\Big)^\varrho, $$
et donc \eqref{W*2} si $y\le
x^{2\epsilon}.$  La mme majoration $g(p)\ll1$ permet d'Žcrire, lorsque $y>
x^{2\epsilon},$
$$\eqalign{{1\over (\log x)^\varrho}&\sum_{p\leqslant x^\varepsilon}{\wt g(p)\log p\over p}\Big(\log {2y\over p}\Big)^{\varrho-1}\ll\epsilon\Big({\log y\over \log
x}\Big)^{\varrho-1}\ll\epsilon\Big({\log y\over \log x}\Big)^{3\varrho_0-1}
\cr&\ll\epsilon^{2\varrho_0}\Big({\log y\over \epsilon\log x}\Big)^{2\varrho_0-1}\Big({\log y\over \log
x}\Big)^{\varrho_0}\ll\epsilon^{\delta}\rst,\cr}
$$ 
puisque $2\varrho_0\geqslant \delta$. Cette estimation est bien compatible avec~\eqref{W*2}. Il reste ˆ Žvaluer la contribution comp\-lŽ\-men\-taire, relative aux \nps\ de l'intervalle $]x^\varepsilon,y]$. Ë cette fin, nous employons l'inŽgalitŽ $|\im g(p)|^2\leqslant 2A\{r(p)-\re g(p)\}$ sous la forme
$$\sum_{x^\varepsilon<p\leqslant z}{\wt g(p)\log p\over p}=\sum_{x^\varepsilon<p\leqslant z}{\{r(p)-\varrho\}\log p\over p}+O\bigg(\sum_{x^\varepsilon<p\leqslant z}{\{r(p)-\re g(p)\}^{w_f}\log p\over p}\bigg).\eqdef{gp}$$
Il rŽsulte alors des relations \eqref{cdpr1} et
\eqref{C2}  que le membre de droite  est $\ll \varepsilon^\delta\log z$. Nous utiliserons Žgalement le fait que la dernire somme en $p$ est ˆ termes positifs ou nuls.\par 
Reportons \eqref{gp} dans l'expression
$$\gG(y):=\sum_{x^\epsilon<p\leqslant y}{\wt g(p)\log
p\over p}\Big(\log {2y\over p}\Big)^{\varrho-1}\qquad (x^{2\varepsilon}<y\leqslant x).$$
Compte tenu de \eqref{C2}, nous pouvons traiter la contribution du terme principal de \eqref{gp}, disons $\gG^+(y)$, par sommation d'Abel. Notant $a_p:=r(p)-\varrho$, nous avons en effet, si $\varrho>1$,
$$\eqalign{\gG^+(y)&=(\varrho-1)\sum_{x^\varepsilon<p\leqslant y}{a_p\log p\over p}\int_1^{2y/p}(\log t)^{\varrho-2}{\dd t\over t}\cr&=(\varrho-1)\int_1^{2y/x^\varepsilon}(\log t)^{\varrho-2}\sum_{x^\varepsilon<p\leqslant 2y/t}{a_p\log p\over p}{\dd t\over t}\ll\varepsilon(\log y)^\varrho\ll\varepsilon^\delta(\log y)^\varrho,\cr}$$
et, si $\varrho<1$, puisque $\delta\leqslant \varrho$,
$$\eqalign{\gG^+(y)&=(1-\varrho)\sum_{x^\varepsilon<p\leqslant y}{a_p\log p\over p}\int_{2y/p}^\infty{\dd t\over t(\log t)^{2-\varrho}}\cr&=(1-\varrho)\int_2^{2y/x^\varepsilon}\sum_{2y/t<p\leqslant y}{a_p\log p\over p}{\dd t\over t(\log t)^{2-\varrho}}+\sum_{x^\varepsilon<p\leqslant y}{a_p\log p\over p}\bigg(\log {2y\over x^\varepsilon}\bigg)^{\varrho-1}\cr
&\ll\int_2^{y^{\varepsilon}}{\dd t\over t(\log t)^{1-\varrho}}+(1-\varrho)\int_{y^{\varepsilon}}^{2y/x^\varepsilon}\varepsilon(\log y){\dd t\over t(\log t)^{2-\varrho}}+\varepsilon(\log y)^\varrho\cr
&\ll\varepsilon^\varrho(\log y)^\varrho\ll\varepsilon^\delta(\log y)^\varrho.\cr}$$
Le cas $\varrho=1$ Žtant trivial, nous avons donc  $\gG^+(y)\ll\varepsilon^\delta(\log y)^\varrho$ en toute circonstance. \goodbreak
\par \goodbreak
ConsidŽrons ˆ prŽsent la contribution $\gG^-(y)$ du terme d'erreur de \eqref{gp} ˆ la quantitŽ $\gG(y)$. Nous avons trivialement $\gG^-(y)\ll\varepsilon^\delta(\log y)^\varrho$ si $\varrho\geqslant 1$. Lorsque $2\gb\leqslant \varrho<1$, nous employons l'inŽgalitŽ de Hšlder avec exposants $\gh/w_f$ et $\gh/(\gh-w_f)$. Nous obtenons la majoration
$$\eqalign{&\ll\bigg\{\sum_{x^\varepsilon<p\leqslant y}{\{r(p)-\re g(p)\}^\gh\log p\over p}\bigg\}^{w_f/\gh}\bigg\{\sum_{x^\epsilon<p\leqslant y}{\log
p\over p}\Big(\log {2y\over p}\Big)^{(\varrho-1)\gh/(\gh-w_f)}\bigg\}^{1-w_f/\gh}\cr
&\ll \varepsilon^{\delta}(\log y)^\varrho,\cr}$$
puisque $${(\varrho-1)\gh\over \gh-w_f}=-{(1-\gb)(1-\varrho)\over 1-\gb-w_f(\varrho-\gb)}\geqslant\gb-1. $$
Cela achve la preuve de \eqref{W*2}.
\par
En reportant
\eqref{W*2} dans
\eqref{est2Ry}, il vient donc
$$N_x(y)\ll |S_x(y)|+\epsilon^{\delta}y\,\e^{\re Z(x;f)}\rst\qquad
(2\leqslant y\leqslant x).\eqdef{est3Ry}$$
Ë ce stade, nous observons que lorsque $|y-x|\leqslant h$ avec $\sqrt x\leqslant h\leqslant \dm x$, on a
$$|N_x(x)-N_x(y)|\ll(\log x)\sum_{x-h\leqslant n\leqslant x+h}\{|g(n)|+\L_x^*\tau_\varrho(n)\} \ll h\e^{Z(x;r)},
$$
d'aprs le thŽorme de Shiu \citer{Sh80}.  Nous pouvons donc Žcrire, compte tenu de
\eqref{est3Ry},
$$N_x(x)\ll{1 \over h}\int_{x-h}^{x}|S_x(y)|\d y+h\,\e^{Z(x;r)}+\epsilon^{\delta}x\, \e^{\re Z(x;f)}. \eqdef{Nxx}$$
Choisissons
$h:=x/\log x\leqslant \epsilon^2 x$. Au vu de \eqref{pqr}, le second terme du membre de droite de~\eqref{Nxx} est
donc dominŽ par le troisime. De plus, d'aprs \eqref{Sx}, nous avons
$$\eqalign{
\int_{x-h}^{x} |S_x(y)|\d y &\ll \int_{x-h}^{x} \sum_{p\leqslant y}
\log p \Big|M\Big({y\over  p};g_x\Big)\Big| \d y\ll \sum_{p \leqslant x} \log p \int_{x-h}^{x}
\Big|M\Big({y\over p};g_x\Big)\Big| \d y \cr
&= \sum_{p \leqslant x} p \log p
\int_{(x-h)/p}^{x/p} |M(y;g_x)| \d y \cr
&= \int_1^x |M(y;g_x)| \sum_{(x-h)/y < p \leqslant x/y}
p \log p \d y  \cr
&\leqslant x \int_1^x {|M(y;g_x)|\over y}
\sum_{(x-h)/y < p \leqslant x/y} \log p \d y .\cr}$$
D'aprs le thŽorme de Brun--Titchmarsh, la somme intŽrieure est
$$\ll {h\log(2x/y)\over y\log(h/y)}$$
ds que $h \geqslant 2y$. Nous utilisons cette majoration lorsque $y\leqslant h/\log x=x/(\log x)^2$. La
contribution correspondante ˆ la dernire intŽgrale est donc
$$\ll h \int_1^x {|M(y;g_x)|\over y^2}\d y .$$
Lorsque $h/\log x < y \leqslant x$, nous majorons trivialement $|M(y;g_x)|$ par \eqref{triv2}. La
contribution correspondante est
$$\eqalign{
&\ll {\e^{Z(x;r)}\over \log x}\sum_{p \leqslant x/(h/\log x)}
\log p \int_{(x-h)/p}^{x/p} \d y\cr
&\ll {h\e^{Z(x;r)}\over \log x}\sum_{p \leqslant  \log ^2x} {\log p
\over p} \ll  {h\e^{Z(x;r)}\log_2x\over \log x}\cr&\ll {x\e^{Z(x;r)}\over \log x}\ll x\epsilon^{2-\beta\gb/2}\e^{Z(x;f)}\ll x\varepsilon^{\delta}\e^{Z(x;f)}.\cr}$$
En reportant nos estimations dans \eqref{Nxx} puis \eqref{maj1My} avec $y=x$, nous obtenons bien~\eqref{maj1Mx}.
\qed
\bigskip\goodbreak
\paradeuxn{Preuve dans le cas exponentiellement multiplicatif}
Sous l'hypothse supplŽmentaire de multiplicativitŽ exponentielle, nous Žtablirons le \ref{cascv} en trois
Žtapes. La premire consiste
ˆ
Žtudier la valeur moyenne  de la fonction exponentiellement multiplicative $g$ dŽfinie en \eqref{g}. Avec les notations $G(s)$ et $\L_x^*$ dŽfinies en \eqref{L*G},  nous posons
$$
\HH_x(\alpha)^2:=\sum_{k\in\z}{1\over k^2+1}\sup_{s\in I_k(\alpha)}\big|G(s)-\L^*_x\zeta(s)^\varrho\big|^2,\eqdef{H}
$$
o, ici et dans la suite, $I_k(\alpha)$ dŽsigne le segment $1+\alpha+ik+i[-\dm,\dm]$. Par ailleurs, nous conservons la notation $\varrho_0:=\min(1,\varrho)/3$ introduite dans la preuve du \ref{lem3}.
\Propp{prop1}{Sous les hypothses du \ref{cascv}, nous avons
$$\Big|M(x;g)-{\e^{-\gamma\varrho}x\,\e^{Z(x;f)}\over \Gamma(\varrho)\log x}\Big|
\ll
{x\over \log x}\Bigg\{
\int_{1/\log x}^{1/(\epsilon\log x)}{\HH_x(\alpha)\over\alpha}\d\alpha
+\epsilon^{\delta}\e^{\re Z(x;f)}
\Bigg\}.\eqdef{HM}$$}
\nid Observons d'emblŽe que, gr‰ce ˆ \eqref{Msggx}, nous pouvons remplacer le membre de gauche de \eqref{HM} par $|M(x;g_x)|$. 
Nous dŽduirons \eqref{HM} de \eqref{maj1Mx}. Ë cette fin,
nous prouvons en premier lieu que, pour $x^\epsilon\le
y\leqslant x$, nous avons
$$\eqalign{{1\over (\log y)^{3/2}}\int_1^y&{|M(z;g_x)|(\log z)^{3/2}\over z^2} \d z\cr&\ll \HH_x\Big({1\over \log
y}\Big)+\epsilon^{\delta}\e^{\re Z(x;f)}\rst.\cr}\eqdef{maj2My}$$  Il est suffisant d'Žtablir cette
majoration lorsque
$|M(z;g_x)|(\log z)^{3/2}$ est remplacŽ dans l'intŽgrande par $|N_x(z)|\sqrt{\log z}$: en vertu de \eqref{maj1My}, l'erreur
impliquŽe est
$$\eqalign{&\ll{\varepsilon^{\delta+\varrho_0}\e^{Z(x;f)}\over (\log y)^{3/2}}\int_1^y{\sqrt{\log
2z}\over z}\d z\ll\varepsilon^{\delta+\varrho_0}\e^{Z(x;f)}.\cr}
$$
Au vu de \eqref{est3Ry},  nous pouvons encore remplacer
$N_x(z)$ par
$S_x(z)$ (\cf\ \eqref{Sx}) avec une erreur
acceptable. Maintenant, l'inŽgalitŽ de Cauchy--Schwarz
$$\int_1^y {|S_x(z)|\sqrt{\log z}\over z^2} \d z \leqslant \biggl\{\int_1^y {|S_x(z)|^2\over z^3} \d z \int_1^y {\log
z
\over z}\d z\biggr\}^{1/2}$$
nous permet de ramener la preuve de \eqref{maj2My} ˆ celle de
$$\int_1^\infty {|S_x(z)|^2 \over z^{3+2\alpha}}\d z \ll {\HH_x(\alpha)^2\over \alpha} \qquad (\alpha >
0).\eqdef{maj1R2}$$
En effet, la majoration souhaitŽe dŽcoulera alors du choix $\alpha = 1/\log y$.
\par\goodbreak
La fonction $G$ Žtant dŽfinie en \eqref{L*G}, posons 
$$\Phi(s):=G(s)-\L^*_x\zeta(s)^\varrho. $$

La relation
$$\int_0^{+\infty} S_x(\e^v) \e^{-v\sigma} \e^{-iv\tau} \d v =-{G'(s)\over G(s)}{\Phi(s)\over s} 
\qquad (\sigma > 1)$$ nous permet d'Žcrire la formule de Plancherel
$$\leqalignno{\int_1^\infty {|S_x(z)|^2 \over z^{3+2\alpha}}\d z &= \int_0^\infty {|S_x(\e^v)|^2 \over \e^{2v(1+\alpha)}}
\d v = {1\over 2\pi i} \int_{1+\alpha-i\infty}^{1+\alpha+i\infty} \Big|{G'(s)\over
G(s)}{\Phi(s)\over s}\Big|^2
\d s&\eqdef{majI1}\cr
&\ll \sum_{k\in \z}{\max_{s\in I_k(\alpha)}
|\Phi(s)|^2 \over k^2+1}
\int_{I_k(\alpha)}\Big|{G'(s)\over
G(s)}\Big|^2 |\rd s|.\cr}$$
En vertu de la multiplicativitŽ exponentielle de $g$, nous avons
$${G'(s)\over G(s)} = - \sum_{p} {g(p)\log p \over p^s},$$
et donc, par le \ref{Ga+},
$$\eqalign{\int_{k-{1\over 2}}^{k+{1\over 2}}
\Big|{G'\over G}
(1+\alpha+i\tau)\Big|^2 \d \tau
= \int_{-{1\over 2}}^{1\over 2} \Big|{G'\over G}
(1+\alpha+ik+i\tau)\Big|^2 \d \tau\ll \sum_{p}{\log p\over p^{1+2\alpha}}
\ll {1\over \alpha}\cdot\cr}$$
En reportant dans \eqref{majI1}, il suit
$$\int_1^\infty {|S_x(z)|^2 \over z^{3+2\alpha}}\d z\ll {\HH_x(\alpha)^2\over \alpha}\cdot\eqdef{majI} $$
 \par
Nous avons ainsi Žtabli \eqref{maj1R2} et donc \eqref{maj2My}.
\par 
Nous sommes maintenant en mesure d'achever la preuve de \eqref{HM}.
Nous utilisons \eqref{maj2My} sous la forme
$$\int_{\sqrt y}^y {|M(z;g_x)|\over z^2} \d z \ll \HH_x \Big({1\over \log y}\Big)+
\epsilon^{\delta}\e^{\re Z(x;f)}\rst\qquad (x^\epsilon\le
y\leqslant x).$$ Il suit
$$\eqalign{
\int_{x^{2\epsilon}}^x{|M(y;g_x)|\over y^2}\d y &\ll
\int_{x^{2\epsilon}}^x {|M(y;g_x)|\over y^2}\int_{\sqrt{y}}^{y} {\d t \over t \log t}\d y
\ll \int_{x^{\epsilon}}^{x} {\d t \over t \log t}
\int_{t}^{\min(t^2,x)} {|M(y;g_x)|\over y^2}\d y\cr
&\ll \int_{x^{\epsilon}}^{x} \Big\{\HH_x\Big({1\over \log
t}\Big)+\epsilon^{\delta}\e^{\re Z(x;f)}\Big[\epsilon^{\varrho_0}+\Big({\log t\over
\log x}\Big)^{\varrho_0}\Big]\Big\} {\d t
\over t
\log t}\cr&\ll\int_{1/\log x}^{1/(\epsilon\log x)} {\HH_x(\alpha)\over \alpha} \d
\alpha+\epsilon^{\delta}\e^{\re Z(x;f)}.\cr}$$ Par ailleurs, la majoration triviale
$$\eqalign{\int_1^{x^{2\epsilon}}{|M(y;g_x)|\over y^2}\d y&\le\sum_{n\leqslant x^{2\epsilon}}{|g_x(n)|\over n}
\cr&\ll\e^{Z(x^{2\epsilon};r)}+\epsilon^\varrho\,\e^{Z(x;r)}\ll\epsilon^\varrho\e^{Z(x;r)}\ll\epsilon^{\delta}\e^{Z(x;f)}\cr}
$$ dŽcoule de \eqref{Trho}, \eqref{L2} et \eqref{pqr}.
En reportant dans \eqref{maj1Mx}, nous obtenons bien \eqref{HM}.\qed
\bigskip
Nous devons maintenant estimer l'intŽgrale figurant au membre de droite de \eqref{HM}. Une majoration de $G(s)$ est
ˆ cet Žgard essentielle. Nous rappelons les notations $\gpp:=\pi\varrho/A$ et $\beta:=1-(\sin\gpp)/\gpp$ introduites dans l'ŽnoncŽ du \ref{cascv}. 
\Propp{prop2}{Sous les conditions
$s=1+\alpha+i\tau$, $1\leqslant \alpha\log x\leqslant 1/\epsilon$, $x\geqslant 2$, nous avons uniformŽment 
$$G(s)\ll\normalbaselineskip30pt\cases{\Big(\dsp{\alpha\over\alpha+|\tau|}\Big)^{\beta\varrho}{\e^{Z(x;f)}\over
(\alpha\log x)^\varrho},& si\quad 
$|\tau|\leqslant 1/(\epsilon\log x)$,\cr
\dsp{\epsilon^{\beta\varrho-\beta\gb/2}\,\e^{Z(x;f)}\over (\alpha\log x)^{(1-\beta)\varrho}},& si\quad 
$ 1/(\epsilon\log x)<|\tau|\leqslant 1/\epsilon^{2\varrho}$.\cr}\eqdef{majG}
$$}
\par
\nid Nous pouvons, sans perte de gŽnŽralitŽ, supposer $\tau>0$. Nous avons
$$G(s)=\exp\Big\{\sum_{p}{g(p)\over p^{1+\alpha+i\tau}}\Big\}\ll\e^{S(\alpha,\tau)}\eqdef{majGbase}$$
o l'on a posŽ
$$S(\eta,\tau):=\re\sum_{p\le\exp(1/\eta)}{g(p)\over p^{1+i\tau}}\qquad (0\leqslant \eta\leqslant 1).$$
\par\goodbreak
L'approximation $p^{i\tau}=1+O(\tau\log p)$ implique clairement $S(\eta,\tau)=S(\eta,0)+O(\tau/\eta)$, et en
particulier
$$S(\alpha,\tau)=S(\alpha,0)+O(1)=\re Z\big(\e^{1/\alpha};f\big)+O(1)\qquad (|\tau|\leqslant \alpha).\eqdef{petittau}$$
De plus, par sommation d'Abel, l'hypothse \eqref{cdpr1} implique  pour \hbox{$x^\epsilon\leqslant v\leqslant w\leqslant x$},
$$\eqalign{\re Z(w;f)-\re Z(v;f)&= Z(w;r)-Z(v;r)+O(1)\cr&\geqslant \varrho\log \Big({\log w\over \log v}\Big)+O(1),\cr}\eqdef{Ev-Eu}$$
o la dernire inŽgalitŽ dŽcoule de \eqref{C2} par une seconde sommation d'Abel.\note{Nous omettons les
dŽtails, qui sont essentiellement identiques ˆ ceux de la preuve de \eqref{Ex-Ey/2}.}
\par
Appliquons \eqref{Ev-Eu} avec $v=\exp(1/\alpha)$, $w=x$,\note{Ce choix est licite puisque 
$1\leqslant \alpha\log x\leqslant 1/\epsilon$.} et reportons dans
\eqref{petittau} puis
\eqref{majGbase}. Il suit
$$G(s)\ll {\e^{Z(x;f)}\over(\alpha\log x)^\varrho}\qquad (|\tau|\leqslant \alpha), $$
ce qui Žtablit bien \eqref{majG} lorsque $|\tau|\leqslant \alpha$.
\par
Supposons maintenant
$$\alpha<\tau\leqslant 1/(\epsilon\log x).\eqdef{cas1} $$
Gr‰ce aux remarques ŽnoncŽes plus haut concernant la quantitŽ $S(\eta,\tau)$, nous pouvons Žcrire 
$$S(\alpha,\tau)=S(\tau,\tau)+D=S(\tau,0)+D+O(1)\eqdef{decSat}$$
avec $$D:=\sum_{\exp(1/\tau)<p\leqslant \exp(1/\alpha)}{\re(g(p)/p^{i\tau})\over
p}\cdot
$$\goodbreak
Pour majorer $D$, nous observons que, pour tous $z\in\CC,\,\theta\in\r$, nous avons 
$$\re\big(z\e^{i\theta}\big)\leqslant |z|\cos\theta+\sqrt{2|z|(|z|-\re z)}.\eqdef{inegrez} $$
En effet, 
$\re\big(z\e^{i\theta}\big)=|z|\cos\theta-(|z|-\re z)\cos\theta-\sin
\theta\,\im z$ et, notant $z=r\e^{i\varphi}$, nous pouvons Žcrire
$$\eqalign{-(|z|-\re z)\cos\theta-&\sin
\theta\,\im z
=-r\big\{(1-\cos\varphi)\cos\theta+\sin\varphi\sin\theta\big\}\cr
&=-2r\sin\dm\varphi\big\{\sin\dm\varphi\cos\theta+\cos\dm\varphi\sin\theta\big\}\cr
&=-2r\sin\dm\varphi\sin(\dm\varphi+\vartheta)\leqslant 2r|\sin\dm\varphi|=\sqrt{2r(|z|-\re z)}.\cr}
$$
En appliquant \eqref{inegrez} avec $z=g(p)$, $\theta=-\tau\log p$, et en notant que $|z|\leqslant r(p)\leqslant A,$ nous dŽduisons de
\eqref{inegrez} et \eqref{cdpr} que, dans le domaine \eqref{cas1}, nous avons
$$\eqalign{D\leqslant D_0+\sqrt{2A}U+O(1)\cr}$$
o, notant $ \D(\alpha,\tau):=]\exp(1/\tau),\exp(1/\alpha)]$, nous avons posŽ
$$\eqalign{D_0&:=\sum_{p\in\D(\alpha,\tau)}{|g(p)|\cos(\tau\log p)\over p}\cr
& = S(\alpha,0)-S(\tau,0)+ \sum_{p\in\D(\alpha,\tau)}{|g(p)|-\re g(p)\over p}-\sum_{p\in\D(\alpha,\tau)}{|g(p)|\{1-\cos(\tau\log p)\}\over p}\cr
&\leqslant S(\alpha,0)-S(\tau,0)-D_1+O(1) ,\cr}$$
avec 
$$D_1:=\sum_{p\in\D(\alpha,\tau)}{r(p)\{1-\cos(\tau\log p)\}\over p},$$
et
$$\eqalign{
U&:=\sum_{x^\epsilon<p\le\exp(1/\alpha)}{\sqrt{|g(p)|-\re g(p)}\over p}\cr
&\ll\bigg\{\sum_{x^\epsilon<p\leqslant x}{1\over p}\bigg\}^{1-1/2\gh}\bigg\{\sum_{x^\epsilon<p\leqslant x}{\{r(p)-\re g(p)\}^\gh\over
p}\bigg\}^{1/2\gh}
\ll\epsilon^{\delta_1/2}\log(1/\epsilon).\cr}\eqdef{majU}
 $$ 
En reportant dans \eqref{decSat} et en faisant de nouveau appel ˆ \eqref{Ev-Eu} sous la forme $$S(\alpha,0)\leqslant \re Z(x;f)-\varrho\log(\alpha\log x)+O(1),$$ nous obtenons que, toujours sous l'hypothse \eqref{cas1},
$$S(\alpha,\tau)\leqslant \re Z(x;f)-\varrho\log (\alpha\log x)-D_1+O(1) .\eqdef{gptau1}$$
\par\goodbreak
Cela Žtant, posons
$\lambda:=\cos\gpp$. Il est facile de vŽrifier que l'on a 
$$\eqalign{r(1-\cos\theta)&\geqslant (1-\lambda)r-(\cos \theta-\lambda)^+r\cr&\geqslant (1-\lambda)r-A(\cos
\theta-\lambda)^+\cr}\qquad (0\leqslant r\leqslant A,\,\theta\in\r).
$$ D'o
$$D_1\geqslant (1-\lambda)\sum_{p\in\D(\alpha,\tau)}{r(p)\over p}
-A\sum_{p\in\D(\alpha,\tau)}{(\cos(\tau\log p)-\lambda)^+\over p}\cdot \eqdef{gptau2}$$
 Par \eqref{Ev-Eu}, nous avons
$$
\eqalign{\sum_{p\in\D(\alpha,\tau)}{r(p)\over p}
\geqslant \varrho\log\Big({\tau\over \alpha}\Big)+O(1).\cr}
$$
 La seconde somme de \eqref{gptau2} peut tre ŽvaluŽe par sommation
d'Abel ˆ l'aide du thŽorme des nombres premiers (\cf, par exemple, \citer{Te15}, Lemme III.4.13). 
Nous obtenons
$$\sum_{p\in\D(\alpha,\tau)}{(\cos(\tau\log p)-\lambda)^+\over p}=m\log\Big({\tau\over \alpha}\Big)+O(1),$$
avec 
$$m:={1\over 2\pi}
\int_{-\gpp}^{\gpp}(\cos\theta-\lambda)\d\theta={\sin\gpp-\gpp\lambda\over\pi}=(1-\beta-\lambda){\varrho\over A}.$$ En
reportant dans \eqref{gptau2}, il vient
$$D_1\geqslant \beta\varrho\log (\tau/\alpha)+O(1), $$
ce qui, compte tenu de \eqref{gptau1} et \eqref{majGbase}, Žtablit bien \eqref{majG} pour l'intervalle
\eqref{cas1}.\note{Notons que cette minoration de $D_1$ est optimale: lorsque $\varrho=2\gb$, elle devient une ŽgalitŽ pour 
$$r(p):=\normalbaselineskip11pt\cases{A& si $\cos(\tau\log p)\geqslant \cos\gpp$,\cr 0& si $\cos(\tau\log
p)<\cos\gpp.$}$$ }
\par\goodbreak
Examinons ˆ prŽsent le cas
$$ 1/(\epsilon\log x)<\tau\le1/\epsilon^{2\varrho}.\eqdef{cas2} $$
Nous devons estimer $S(\alpha,\tau)$. Ë cette fin, nous majorons $\re (g(p)/p^{i\tau})$ par $r(p)$ si $p\le
x^\epsilon$ et par \eqref{inegrez} lorsque $x^\epsilon<p\leqslant \exp(1/\alpha)$. En employant la majoration 
\eqref{majU} pour $U$, il vient
$$\eqalign{S(\alpha,\tau)&\leqslant \sum_{p\leqslant x^\varepsilon}{r(p)\over p}+\sum_{x^\varepsilon<p\leqslant \e^{1/\alpha}}{|g(p)|\cos(\tau\log p)\over p}+O(1)\cr&\leqslant Z\big(\e^{1/\alpha};r\big)-D_2+O(1)\cr}\eqdef{majS2} $$
avec
$$D_2:=\sum_{x^\epsilon<p\leqslant \exp(1/\alpha)}{r(p)\{1-\cos(\tau\log p)\}\over p}\cdot $$
Nous appliquons ˆ $D_2$ la mŽthode de minoration utilisŽe pour $D_1$, en observant que l'hypothse 
$\epsilon\log x>\sqrt{\log x}$ implique, pour $x$ assez grand, $|\tau|\leqslant 1/\epsilon^{2\varrho}\leqslant \e^{\sqrt{\epsilon\log
x}}$ d'o, gr‰ce au lemme III.4.13 de \citer{Te15},  
$$
\sum_{x^\epsilon<p\leqslant \exp(1/\alpha)}{(\cos(\tau\log p)-\lambda)^+\over
p}=-m\log(\alpha\epsilon\log x)+O(1).
$$
Nous obtenons ainsi
$$D_2\geqslant -\beta\varrho\log (\alpha\epsilon\log x)+O(1). $$
\par
Maintenant, nous dŽduisons de \eqref{Ev-Eu} et de \eqref{pqr} que 
$$Z(\e^{1/\alpha};r)\leqslant \re Z(x;f)-\varrho\log (\alpha\log x)+\ft12\beta\gb\log (1/\epsilon)+O(1).$$
\par
En reportant les deux dernires inŽgalitŽs dans \eqref{majS2} puis \eqref{majGbase}, nous obtenons bien la
seconde majoration de \eqref{majG}.
\qed\goodbreak
\medskip
Nous sommes maintenant en mesure de fournir un premier jeu d'estimations pour la fonction $\HH_x(\alpha)$ dŽfinie
par \eqref{H}. Nous posons
$$W(s):=\sum_{p}{g(p)\over p^{s}}-\varrho\log \zeta(s)-\sum_{p\leqslant x}{g(p)\over
p}+\varrho\log_2 x+\gamma\varrho\qquad (\re s>1)\eqdef{W}$$
et, pour $0\leqslant u\leqslant 1 $,
$$\I_x(\alpha,u):=\sup_{\di{\sigma=1+\alpha}{|\tau|\leqslant \alpha/\epsilon^u}}{\big|\e^{W(s)}-1|\over (1+|\tau|/\alpha)^\varrho},\qquad
\J_x(\alpha,u):=\min\big(1,\I_x(\alpha,u)\big).\eqdef{ix}
$$
\Propp{prop3}{Pour $1\leqslant \alpha\log x\leqslant 1/\epsilon$, $0\leqslant u\leqslant 1$, nous avons
$$\HH_x(\alpha)\ll{\e^{Z(x;f)}\over (\alpha\log
x)^\varrho}\Big\{\J_x(\alpha,u)+\epsilon^{\beta u\varrho}+{(\alpha\epsilon\log x)^{\beta\varrho}\over
\epsilon^{\beta\gb/2}}\Big\}.\eqdef{redH}
$$}
\par
\nid L'estimation classique
$$\zeta(s)\ll{1\over |s-1|}+\log (1+|s|)\qquad (\sigma>1)\eqdef{Z}$$
implique
$$\eqalign{\sup_{\di{\sigma=1+\alpha}{\alpha/\epsilon^u\leqslant |\tau|\leqslant 1}}|\L^*_x\zeta(s)^\varrho\big|^2+\sum_{|k|\ge
1}{\sup_{s\in I_k(\alpha)} \big|\L^*_x\zeta(s)^\varrho\big|^2\over
k^2+1}&\ll{\e^{2Z(x;f)}\over (\log
x)^{2\varrho}}\Big\{{\epsilon^{2u\varrho}\over \alpha^{2\varrho}}+1\Big\}.\cr}\eqdef{ZH}$$
Par ailleurs, la majoration triviale
$$G(s)\ll\e^{Z(\exp(1/\alpha);r)}\ll{\e^{Z(x;r)}\over (\alpha\log
x)^\varrho}\ll{\epsilon^{-\beta\gb/2}\e^{Z(x;f)}\over (\alpha\log x)^\varrho}\qquad (\sigma=1+\alpha)$$ 
fournit  
$$\sum_{|k|>1/\epsilon^{2\varrho}}{\sup_{s\in I_k(\alpha)} |G(s)|^2\over k^2+1}\ll
{\epsilon^{2\varrho-\beta\gb}\e^{2Z(x;f)}\over (\alpha\log x)^{2\varrho}}\eqdef{GH1}$$ alors que \eqref{majG} implique
$$\eqalign{\sup_{\di{\sigma=1+\alpha}{\alpha/\epsilon^u\leqslant |\tau|\leqslant 1}}|G(s)\big|^2+&\sum_{1\leqslant |k|\le
1/\epsilon^{2\varrho}}{\sup_{s\in I_k(\alpha)} \big|G(s)\big|^2\over
k^2+1}\cr&\ll{\e^{2Z(x;f)}\over (\alpha\log x)^{2\varrho}}\Big\{\epsilon^{2\beta\varrho
u}+{(\alpha\epsilon\log x)^{2\beta\varrho}\over\epsilon^{\beta\gb}}\Big\}.\cr}\eqdef{GH2}
$$
\par
On constate facilement que les membres de droite de \eqref{ZH} et \eqref{GH1} sont dominŽs par celui de
\eqref{GH2}. Nous pouvons donc Žcrire
$$\HH_x(\alpha)\ll\sup_{\di{\sigma=1+\alpha}{|\tau|\le
\min(1,\alpha/\epsilon^u)}}\big|G(s)-\L^*_x\zeta(s)^\varrho\big|+{\e^{\re Z(x;f)}\over (\alpha\log
x)^{\varrho}}\Big\{\epsilon^{\beta\varrho u}+{(\alpha\epsilon\log x)^{\beta\varrho}\over \epsilon^{\beta\gb/2}}\Big\}.$$
Lorsque $\sigma=1+\alpha,\,|\tau|\le
\min(1,\alpha/\epsilon^u)$, nous avons
$$G(s)-\L^*_x\zeta(s)^\varrho={\e^{Z(x;f)-\gamma\varrho}\zeta(s)^\varrho\over (\log x)^\varrho}\Big\{\e^{W(s)}-1\Big\}\ll{\e^{Z(x;f)}\I_x(\alpha,u)\over
(\alpha\log x)^\varrho}\cdot $$
De plus, \eqref{majG} et \eqref{Z} impliquent Žgalement que le membre de gauche est
$$\ll{\e^{Z(x;f)}\over (\alpha\log x)^{\varrho}}\Big\{1+{(\alpha\epsilon\log x)^{\beta\varrho}\over
\epsilon^{\beta\gb/2}}\Big\} .$$
Cela termine la preuve de \eqref{redH}.\qed
\medskip\goodbreak
Il nous faut ˆ ce stade disposer d'estimations pour la quantitŽ $\I_x(\alpha,u)$ dŽfinie en \eqref{ix}.
\Propp{prop4}{Pour $0\leqslant u\leqslant  1-\delta$, $1\leqslant \alpha\log
x\leqslant \epsilon^{u+\delta-1}$, nous avons
$$\I_x(\alpha,u)\ll \epsilon^\delta\big\{1+\log (\alpha\log x)\big\}. \eqdef{majJx}$$}
\par
\nid La fonction $\zeta(s)\exp\big\{-\sum_{p}1/p^{s}\big\}$ est analytique en $s=1$. On a donc, au
voisinage de $s=1$ et pour $\re s>1$,
$$\log \zeta(s)=\sum_{p}{1\over p^s}-b+O(s-1) \quad \hbox{avec}\quad b=\sum_{p}\Big\{{1\over p}+\log \Big(1-{1\over
p}\Big)\Big\}.$$  D'aprs la formule de Mertens, nous avons Žgalement
$$\sum_{p\leqslant x}{1\over p}=\log_2x+b+\gamma+O\Big({1\over \log x}\Big)\qquad (x\geqslant 2) .$$
Soit alors $s=1+\alpha+i\tau$, o $\alpha$ et $\tau$ satisfont aux conditions de l'ŽnoncŽ. Comme nous avons  $1/\log x\leqslant |s-1|$ dans le domaine de variation en $\alpha$ considŽrŽ, 
nous obtenons, pour $|\tau|\le
\alpha/\epsilon^u\leqslant \epsilon^{1+\delta}\ll1$, 
$$W(s)=\sum_{p}{g(p)-\varrho\over p^s}-\sum_{p\leqslant x}{g(p)-\varrho\over
p}+O(s-1).\eqdef{Wat}$$
Posons $a_p:=\varrho-g(p)$, de sorte que $a_p=0$ pour $p>x$, et dŽsignons par $W_1(s)$ le terme principal de \eqref{Wat}, soit 
$$W_1(s):=\sum_{p\leqslant x}{a_p\over p}\Big(1-{1\over p^{s-1}}\Big). $$
\par
Nous scindons la somme en $p$ sous la forme $W_1(s)=W_{11}(s)+W_{12}(s)$,  correspondant aux
domaines  de sommation respectifs
$p\leqslant x^\epsilon$ et $x^\epsilon<p\leqslant x$.
\par
Nous estimons $W_{11}(s)$ en employant la majoration triviale $a_p\ll1$ et l'inŽgalitŽ $$|1-p^{1-s}|\leqslant |s-1|\log
p.$$ Il suit
$$\eqalign{W_{11}(s)&\ll|s-1|\sum_{p\leqslant x^\epsilon}{\log p\over p}\ll|s-1|\epsilon\log x\cr&\ll(\epsilon\alpha\log
x)(1+|\tau|/\alpha)\ll\epsilon^{\delta}.\cr} \eqdef{W11}$$
\par
Par \eqref{gp}, \eqref{cdpr1} et \eqref{C2}, nous disposons d'une dŽcomposition canonique $a_p=b_p+O(c_p)$ $(p\leqslant x)$, avec $c_p\geqslant 0$, et
$$\eqalign{B(y)&:=\sum_{x^\varepsilon<p\leqslant \exp y}{b_p\log p\over p}\ll\varepsilon y,\cr C(y)&:=\sum_{x^\varepsilon<p\leqslant \exp y}{c_p\log p\over p}\ll\varepsilon^{\delta} y\cr}\qquad (\varepsilon\log x<y\leqslant \log x).$$
On en dŽduit que
$$\eqalign{\int_{\epsilon\log x}^{\log x}&{1-\e^{(1-s)y}\over y}\d B(y)\cr
&\ll\epsilon+\epsilon\int_{\epsilon\log x}^{\log x}\abs{(s-1)y\e^{(1-s)y}+\e^{(1-s)y}-1\over y}\d y\cr
&\ll\epsilon\bigg\{1+|s-1|^2\int_{\varepsilon \log x}^{1/|s-1|}y\d y+\int_{1/|s-1|}^{\log x}\bigg({|s-1|\over \e^{\alpha y}}+{1\over y}\bigg)\d
y\bigg\}\cr&\ll\epsilon\Big\{1+{|\tau|\over \alpha}+\log (\alpha\log x)\Big\}\ll\varepsilon^{1-u}\ll\varepsilon^{\delta}\{1+\log (\alpha\log x)\},\cr}\eqdef{W12}
$$
et
$$\eqalign{\sum_{x^\varepsilon<p\leqslant x}{c_p\over p}\abs{1-{1\over p^{s-1}}}&\ll|s-1|\sum_{x^\varepsilon<p\leqslant \exp(1/|s-1|)}{c_p\log p\over p}+\sum_{\exp(1/|s-1|)<p\leqslant x}{c_p\over p}\cr&\ll\varepsilon^\delta+\varepsilon^\delta\log (|s-1|\log x).\cr}\eqdef{somcp}$$
\par
En regroupant les estimations \eqref{Wat}, \eqref{W11},\eqref{W12} et \eqref{somcp}, nous pouvons donc Žnoncer que, sous les hypothses effectuŽes, nous avons
$$W(s)\ll\epsilon^{\delta}\{1+\log (1+|\tau|/\alpha)+\log (\alpha\log x)\},\eqdef{majW}$$ puisque le terme d'erreur de
\eqref{Wat} est $\ll \alpha+|\tau|\ll \alpha/\varepsilon^u\ll\epsilon^{1+\delta}.$\par 
En reportant dans \eqref{ix}, nous obtenons bien l'estimation \eqref{majJx} annoncŽe.
\qed
\goodbreak
\Propp{prop5}{
Nous avons
$$\int_{1/\log x}^{1/(\epsilon\log x)}{\HH_x(\alpha)\over
\alpha}\d\alpha\ll\epsilon^{\delta}\e^{Z(x;f)}.\eqdef{P5} $$ }
\par
\nid 
Fixons une fois pour toutes $u:=\ft13\leqslant 1-\delta$. 
\par
DŽsignons par $\K(\epsilon,x)$ l'intŽgrale figurant au membre de gauche de \eqref{P5}. La relation \eqref{redH}
fournit d'abord
$$\K(\epsilon,x)\ll\e^{
Z(x;f)}\bigg\{\int_{1/\log x}^{1/(\epsilon\log
x)}{\J_x(\alpha,u)\over (\alpha\log x)^\varrho}{\d\alpha\over \alpha}+\epsilon^{\beta
\varrho/3}+\epsilon^{\beta\varrho-\beta\gb/2}\log (1/\varepsilon)\bigg\}.$$
Les deux derniers termes dans l'accolade sont $\ll\varepsilon^\delta$ puisque nos hypothses impliquent $\beta\varrho/3\geqslant 2\gb\beta/3\geqslant \delta$ et $\beta\varrho-\beta\gb/2\geqslant 2\delta$.\par 
En
majorant
$\J_x(\alpha,u)$ par $O(1)$ dans le domaine $\epsilon^{-1/3}<\alpha\log x\le
1/\epsilon$, nous obtenons de plus
$$\int_{\epsilon^{-1/3}/\log x}^{1/\epsilon\log x}{\J_x(\alpha,u)\over
(\alpha\log x)^\varrho}{\d\alpha\over \alpha}\ll\epsilon^{\varrho/3}\ll\varepsilon^\delta. 
$$
\par 
Pour estimer l'intŽgrale complŽmentaire, relative au domaine $1\leqslant \alpha\log x\leqslant \varepsilon^{-1/3}$, nous appliquons \eqref{majJx} avec $u=\ft13$. Cela complte la preuve de \eqref{P5}.
\qed
\medskip\goodbreak
En concatŽnant les Propositions \ref{sprop1} et \ref{sprop5}, nous obtenons le rŽsultat suivant.
\Propc{cor4}{Dans les hypothses du \ref{cascv}, la fonction exponentiellement multiplicative $g$ dŽfinie en \eqref{g} vŽrifie
$$M(x;g)={\e^{-\gamma\varrho+Z(x;f)}x\over \Gamma(\varrho)\log
x}\big\{1+O\big(\epsilon^\delta\big)\big\}. $$}
\par
\medskip\goodbreak
\paradeuxn{ComplŽtion de l'argument}\drefdeux{compcascv}
 Soit $f$ une \fmu\ satisfaisant aux hypothses du \ref{cascv}. Consi\-dŽrons la fonction exponentiellement multiplicative $g$ dŽfinie en \eqref{g}. Nous pouvons,
sans affecter la gŽnŽralitŽ, supposer que $f(p^\nu)=g(p^\nu)=\varrho^\nu/\nu!$ pour $p>x$ et $f(\pnu)=0$ si $p\leqslant x$, $\pnu>x$. Nous avons alors $f=g*h$ o $h$ est dŽfinie par \eqref{defh}. Notons que, bien que les dŽfinitions  \eqref{defg} et \eqref{g} ne co•ncident pas, la relation de convolution dŽfinit bien la mme fonction~$h$. Ainsi, $h(p)=0$ pour tout \np\ $p$, $h(p^\nu)=0$ si $p>x$, et la majoration \eqref{majQ} demeure valide.
\par 
Nous avons
$$M(x;f)=\sum_{n\leqslant x}h(n)M\Big({x\over n};g\Big). \eqdef{conv}$$
\par
 Appliquons le \ref{cor4} ˆ $M(x/n;g)$ pour $n\le
\sqrt{x}$. Nous obtenons
$$M\Big({x\over n};g\Big)={\e^{-\gamma\varrho+Z(x;f)}x\over \Gamma(\varrho)n\log
x}\Big\{1+O\Big(\epsilon^\delta+{\log n\over \log x}\Big)\Big\}\qquad (n\le
\sqrt{x}),$$
o l'on a utilisŽ  l'estimation
$$Z(x;f)-Z(x/n;f)\ll\sum_{x/n<p\leqslant x}{1\over p}\ll{\log n\over \log x}\qquad (n\le\sqrt{x}). $$
Gr‰ce ˆ \eqref{conv} et \eqref{triv2}, il suit 
$$M(x;f)={\e^{-\gamma\varrho+Z(x;f)}x\over \Gamma(\varrho)\log
x}\bigg\{\sum_{n\leqslant \sqrt{x}}{h(n)\over
n}+Y_1+Y_2+Y_3\bigg\}\eqdef{Mf}$$
avec, compte tenu de \eqref{prodh1} et \eqref{Ev-Eu},
$$\eqalign{Y_1&\ll \varepsilon^\delta\sum_{n\leqslant \sqrt{x}}{|h(n)|\over
n}\ll \varepsilon^\delta,\cr
Y_2&\ll {1 \over \log x}\sum_{n\leqslant \sqrt{x}}{|h(n)|\log n\over n},\cr
Y_3&\ll\epsilon^{-\gb/2}\sum_{\sqrt{x}<n\leqslant x}{|h(n)|\over
n}\Big\{{\varepsilon^{\varrho}\log x\over
\log (2x/n)}+{\Big({\log 2x/n)\over \log x}\Big)^{\varrho-1}}\Big\}.
\cr}$$
\par 
En vertu de \eqref{majQ}, nous pouvons Žcrire
$$\eqalign{Y_2&\ll{1\over \log x}\int_1^{\sqrt{x}}\log y\,\d Q(y)\ll Q\big(\sqrt{x}\big) +{1\over \log x}\int_1^{\sqrt{x}}{Q(y)\over
y}\d y\ll{\log_2x\over \log x}\ll\epsilon.\cr}$$
De mme
$$\eqalign{Y_3&\ll\epsilon^{-\gb/2}\int_{\sqrt{x}}^{x}\Big\{{\epsilon^\varrho\log x\over
\log (2x/y)}+{\Big({\log 2x/y)\over \log x}\Big)^{\varrho-1}}\Big\}\d Q(y)\cr&\ll\epsilon^{3\varrho/4}
+\epsilon^{3\varrho/4}\int_{\sqrt{x}}^{x}{Q(y)\log x\over y(\log 2x/y)^2}
\d y+|\varrho-1|\varepsilon^{-\gb/2}\int_{\sqrt{x}}^{x}{(\log 2x/y)^{\varrho-2}Q(y)\over y(\log x)^{\varrho-1}}\d
y\cr&\ll\epsilon^{3\varrho/4}+{\log_2x\over (\log x)^{\min(1,\varrho)}}\ll\epsilon^{\delta}.\cr}$$ En reportant dans
\eqref{Mf} o la somme en $n$ est Žtendue jusqu'ˆ l'infini au prix d'un terme d'erreur supplŽmentaire
$\ll1/\log x\ll\epsilon$ et en tenant compte de \eqref{prodh2}, nous obtenons bien~\eqref{moy}. 
\bigskip
\paraunn{Preuve du \ref{fr}}
\paradeuxn{InŽgalitŽ de Tur‡n-Kubilius pondŽrŽe}
Le rŽsultat suivant, qui possde un intŽrt propre, n'est pas {\it stricto sensu} indispensable ˆ la dŽmonstration. Il permettra cependant de simplifier la prŽsentation. Un rŽsultat de mme nature, mais relatif ˆ des hypothses sensiblement diffŽrentes des n™tres, a ŽtŽ Žtabli par Bir— \& Szamuely \citer{BS96}.
\Propl{TKpond}{Soient $A$, $B$ des constantes positives, $x\geqslant 1$, $\lambda$ une \fmu\ positive ou nulle de la classe $\M_0(x;A,B)$, et $\vartheta$ une fonction additive complexe. Posons
$$\Theta(x):=\sum_{\pnu\leqslant x}{\lambda(\pnu)\vartheta(\pnu)\over \pnu},\ \Theta_0(x):=\sum_{\pnu\leqslant x}{\lambda(\pnu)|\vartheta(\pnu)|\over \pnu},\  \gS(x):=\sum_{\pnu\leqslant x}{\lambda(\pnu)|\vartheta(\pnu)|^2\over \pnu}\cdot\eqdef{Th}$$
Nous avons alors
$$\sum_{n\leqslant x}{\lambda(n)\over n}\big|\vartheta(n)-\Theta(x)\big|^2\ll\sum_{n\leqslant x}{\lambda(n)\over n}\Big\{\gS(x)+\Theta_0(x)\sqrt{\gS(x)}\Big\}.\eqdef{TKp}$$}
\goodbreak
\nid Nous pouvons supposer sans perte de gŽnŽralitŽ que $\vartheta$ est ˆ valeurs rŽelles positives. Le rŽsultat gŽnŽral s'en dŽduit en introduisant classiquement les parties positives et nŽgatives des parties rŽelle et imaginaire de $\vartheta$. \par 
Observons d'abord qu'il rŽsulte du thŽorme 1.1 de \citer{Te16} que 
$$L:=\sum_{n\leqslant x}{\lambda(n)\over n}\asymp\prod_{p\leqslant x}\Big(1+{\lambda(p)\over p}\Big).\eqdef{SP}$$
Ensuite, nous dŽveloppons le membre de gauche de \eqref{TKp} sous la forme
$$M_2-2\Theta(x)M_1+\Theta(x)^2L,\eqdef{devcar}$$
o l'on a posŽ
$$M_j:=\sum_{n\leqslant x}{\lambda(n)\vartheta(n)^j\over n}\quad(j=1,2).$$
\par 
Nous avons
$$M_1=\sum_{\pnu\leqslant x}{\lambda(\pnu)\vartheta(\pnu)\over \pnu}\sum_{\di{m\leqslant x/\pnu}{p\,\nmid\,m}}{\lambda(m)\over m}=\Theta(x)L-\sum_{\pnu\leqslant x}{\lambda(\pnu)\vartheta(\pnu)\over \pnu}\big\{R_1(\pnu)+R_2(\pnu)\big\},$$
avec
$$R_1(\pnu):=\sum_{k\geqslant 1}{\lambda(\pk)\over \pk}\sum_{\di{n\leqslant x/p^{\nu+k}}{p\,\nmid\, n}}{\lambda(n)\over n},\qquad R_2(\pnu):=\sum_{x/\pnu<n\leqslant x}{\lambda(n)\over n}\cdot$$
Nous estimons $R_1(\pnu)$ en majorant trivialement la somme intŽrieure par $L$. Il suit
$$R_1(\pnu)\leqslant L\Big({\lambda(p)\over p}+\varepsilon_p\Big),$$
o $\varepsilon_p\geqslant 0$ est le terme gŽnŽral d'une sŽrie convergente. La contribution ˆ \eqref{devcar} est donc
$$\eqalign{&\ll L\sum_{\pnu\leqslant x}{\lambda(\pnu)\vartheta(\pnu)\over \pnu}\sum_{\pmu\leqslant x}{\lambda(\pmu)\vartheta(\pmu)\over \pmu}\Big({\lambda(p)\over p}+\varepsilon_p\Big)\cr
&\ll L\Theta_0t(x)\sum_{p\leqslant x}\Big(\sum_{\nu\leqslant (\log x)/\log p}{\lambda(\pnu)\vartheta(\pnu)^2\over \pnu}\Big)^{1/2}\Big({\lambda(p)\over p}+\varepsilon_p\Big)^{3/2}\ll L\Theta_0(x)\gS(x)^{1/2}\cr}$$
Pour Žvaluer $R_2(\pnu)$, nous utilisons \eqref{triv1} et \eqref{SP} sous la forme
$$\sum_{x/\e^{k+1}<n\leqslant x/\e^k}\lambda(n)\ll {xL\over \e^k\{\log(2x)-k\}}\qquad \big(0\leqslant k\leqslant\log \pnu\big).$$
Par sommation sur $k$, nous obtenons
$$R_2(\pnu)\ll L\log \Big({\log 2x\over \log (2x/\pnu)}\Big).$$
Pour les grandes valeurs de $\pnu$, cette majoration est moins prŽcise que l'inŽgalitŽ triviale $R_2(\pnu)\leqslant L$ ; elle sera nŽanmoins suffisante pour la suite. 
Posant $\nu_p:=\fl{(\log x)/\log p}$, nous obtenons en effet que la contribution globale correspondante, disons $S_2$, ˆ \eqref{devcar} vŽrifie 
$$\eqalign{S_2\ll L\Theta(x)\sum_{p\leqslant x}\sum_{\nu\leqslant \nu_p}{\lambda(\pnu)\vartheta(\pnu)\over \pnu}\sum_{\mu\leqslant \nu_p}{\lambda(\pmu)\vartheta(\pmu)\over \pmu}\log \Big({\log 2x\over \log (2x/\pmu)}\Big).\cr}$$
En observant que
$$\log \Big({\log 2x\over \log (2x/\pmu)}\Big)\ll{\log_2x\over \log x}\log \pmu\qquad(\pmu\leqslant x),$$ et en traitant les deux sommes intŽrieures par l'inŽgalitŽ de Cauchy-Schwarz, nous pouvons Žcrire
$$\eqalign{S_2&\ll L{\log_2x\over \log x}\sum_{p\leqslant x}\sum_{\nu\leqslant \nu_p}{\lambda(\pnu)\vartheta(\pnu)^2\over \pnu}\bigg\{\sum_{\nu\leqslant \nu_p}{\lambda(\pnu)\over \pnu}\sum_{\mu\leqslant \nu_p}{\lambda(\pmu)(\log \pmu)^2\over \pmu}\bigg\}^{1/2}\cr
&\ll L\sum_{\pnu\leqslant x}{\lambda(\pnu)\vartheta(\pnu)^2\over \pnu}=L\gS(x),\cr}$$
puisque la quantitŽ entre accolades est $\ll\log x$.\par 
Nous dŽduisons de ce qui prŽcde que le membre de gauche de \eqref{TKp} vaut
$$M_2-\Theta(x)^2L+O\Big(L\gS(x)+L\Theta_0(x)\sqrt{\gS(x)}\Big).$$
\vskip-1.5mm
Or
$$\eqalign{M_2&=\sum_{n\leqslant x}{\lambda(n)\over n}\bigg\{\sum_{\pnu\| n}\vartheta(\pnu)^2+\sum_{\di{\pnu\|n,\,\pmu\| n}{p\neq q}}\vartheta(\pnu)\vartheta(\qmu)\bigg\}\cr&\ll\sum_{\pnu\leqslant x}{\lambda(\pnu)\vartheta(\pnu)^2\over \pnu}\sum_{\di{n\leqslant x/\pnu}{p\,\nmid\, n}}{\lambda(n)\over n}+\sum_{\di{\pnu\leqslant x,\,\pmu\leqslant x }{p\neq q}}{\lambda(\pnu)\lambda(\qmu)\vartheta(\pnu)\vartheta(\qmu)\over \pnu\qmu}\sum_{\di{n\leqslant x/\pnu\qmu}{(n,pq)=1}}{\lambda(n)\over n}\cr
&\leqslant \gS(x)L+\Theta(x)^2L.\cr}$$
Cela achve la dŽmonstration.
\qed
\goodbreak
\bigskip
\paradeuxn{ComplŽtion de l'argument}
Nous pouvons supposer $\varepsilon$ arbitrairement petit. En effet, dans le cas contraire, le rŽsultat dŽcoule de l'inŽgalitŽ $|f|\leqslant r$ et de la majoration de Halberstam--Richert \citer{HR79}.\par 
Par ailleurs, nous pouvons supposer sans perte de gŽnŽralitŽ que $f(\pnu)=0$ si $\pnu>x$.\par 
Soit $$\K:=\bigg]{\log (1/\varepsilon_1)\over \log (1+\varepsilon_1)}, {\log_2 x\over \log (1+\varepsilon_1)}-1\bigg]\cap\N.$$ Appliquons \eqref{moyrp} avec $y=y_k:=\exp\{(1+\varepsilon_1)^k\}$ pour $k\in\K$. Nous obtenons
$$\sum_{y_k<p\leqslant y_{k+1}}{r(p)\log p\over p}=\gb_k\varepsilon_1\log y_k\qquad \big(k\in\K\big), \eqdef{moylocrp}$$ 
avec $$4\gb\leqslant \gb_k\leqslant 2A+O\Big(\e^{-\sqrt{\log y_k}}\Big).$$
\goodbreak
Nous
posons alors $s(p):=2\gb r(p)/\gb_k$ pour $y_k<p\leqslant y_{k+1}$, $k\in\K$, et $s(p):=0$ pour $p\in[2,x]\sset\cup_{k\in\K}\,]y_k,y_{k+1}]$. Nous avons clairement $0\leqslant s(p)\leqslant \dm r(p)$ pour $p\leqslant x$. Par sommation, il suit
$$\sum_{p\leqslant y}{\{s(p)-2\gb\}\log p\over p}\ll\varepsilon_1\log y\qquad (\e^{1/\varepsilon}<y\leqslant x).\eqdef{regs}$$  
Introduisons alors $\vartheta(\pnu):=\arg f(\pnu)\in\,]-\pi,\pi]$, avec la convention $\vartheta(\pnu)=0$ si $f(\pnu)=0$, et dŽfinissons trois \fas\ multiplicatives $t$, $\psi$, et $\varphi$, par les formules
$$t(\pnu):=\normalbaselineskip=15pt\cases{r(p)-s(p)& si $\nu=1$\cr
r(\pnu)& si $\nu\geqslant 2$\cr},\quad\psi(\pnu):=\normalbaselineskip=15pt\cases{f(p)-s(p)\e^{i\vartheta(p)} & si $\nu=1$\cr t(\pnu)\e^{i\vartheta(\pnu)}& si $\nu\geqslant 2$},\quad \varphi*\psi=f.$$ 
Nous prolongeons alors $s$ en une \fmu\ par la formule $s*t=r$.
\par \goodbreak
Notons d'emblŽe que l'inŽgalitŽ $s(p)\leqslant \dm r(p)$ implique que $s(p)\leqslant t(p)$ pour tout \np\ $p$. En particulier, nous avons donc $|\psi|\leqslant t$.  De plus, il rŽsulte de la relation $\varphi(p)=s(p)\e^{i\vartheta(p)}$, valable pour tout \np\ $p$, que, notant $\gh:=(1-\gb)/\gb$, nous avons, pour $x^\varepsilon<y\leqslant x$,
$$\eqalign{\sum_{x^\varepsilon<p\leqslant y}{\{s(p)-\re\varphi(p)\}^\gh\log p\over p}&=\sum_{x^\varepsilon<p\leqslant y}{s(p)^\gh\{1-\cos\vartheta(p)\}^\gh\log p\over p}\cr&\ll \varepsilon^{\delta_1\gh}\log y=\varepsilon_1^{2\delta_1\gh}\log y,\cr}\eqdef{sp-phip}$$
puisque $s(p)\{1-\cos\vartheta(p)\}\leqslant r(p)-\re f(p)$ pour tout $p$.\par 
Notant $s^*$ la \fmu\ dŽfinie par $s^*(p)=s(p)$, $s^*(\pnu)=|\varphi(\pnu)|$ $(\nu\geqslant 2)$, il rŽsulte de
 \eqref{regs} que, quitte ˆ altŽrer la constante $B$, le couple $(\varphi,s^*)$ satisfait pour tout $y\in[x^{\varepsilon_1},x]$, aux hypothses du \ref{cascv} avec $\varepsilon_1$ au lieu de $\varepsilon$, $2\delta_1$ au lieu de $\delta_1$, $\varrho:=2\gb$. En effet, nous avons bien $2\delta_1\leqslant 2\delta_0(\gb)=\ft23\beta\gb$ o $\beta$ est dŽfini par \eqref{defpar} avec $\varrho=2\gb$, et $\varepsilon_1\geqslant 1/\sqrt{\varepsilon\log x}\geqslant 1/\sqrt{\log y}$. Nous pouvons donc Žcrire
$$M(y;\varphi)={\e^{-2\gamma\gb}y\over \Gamma(2\gb)\log y}\bigg\{\prod_{p}\sum_{\pnu\leqslant y}{\varphi(\pnu)\over \pnu}+O\big(\varepsilon^\delta\e^{Z(y;\varphi)}\big)\bigg\}\qquad (x^{\varepsilon_1}\leqslant y\leqslant x)$$
avec $\delta:=w_f\delta_1$. Rappelons la notation $\gc:=\gb/A$. En notant que $$s(p)-\re\varphi(p)\geqslant \gc\{|f(p)|-\re f(p)\}\qquad (p\leqslant x),$$  il suit
$$\eqalign{M(x;f)&={\e^{-2\gamma\gb}x\over \Gamma(2\gb)}\sum_{x^{\varepsilon_1}<d\leqslant x^{1-\varepsilon_1}}\bigg\{\prod_{p}\sum_{\pnu\leqslant x/d}{\varphi(\pnu)\over \pnu}\bigg\}{\psi(d)\over d \log (x/d)}\cr&\quad+\varepsilon^\delta \e^{-\gc Z(x;|f|-f)}x\,\gR_0+x\,\gR_1+x\,\gR_2\cr}\eqdef{decMf}$$
avec $$\eqalign{&\gR_0\ll\sum_{x^{\varepsilon_1}<d\leqslant x^{1-\varepsilon_1}}{\e^{Z(x/d;s)}t(d)\over \log (x/d)},\cr
&\gR_1\ll\sum_{d\leqslant x^{\varepsilon_1}}{t(d)\,\e^{Z(x;s)}\over d\log x},\quad \gR_2\ll \sum_{x^{1-\varepsilon_1}<d\leqslant x}{t(d)\,\e^{Z(x/d;s)}\over d\log (2x/d)}\cdot\cr}$$
\par \goodbreak
Ë ce stade, observons que la formule \eqref{regs} et l'inŽgalitŽ $t(p)\geqslant s(p)$ impliquent, par sommation d'Abel, 
$$\eqalign{Z(x;s)-Z(y;s)&= \{2\gb+O(\varepsilon_1)\}\log \Big({\log x\over \log y}\Big)\cr
Z(x;t)-Z(y;t)&\geqslant  \{2\gb+O(\varepsilon_1)\}\log \Big({\log x\over \log y}\Big)+O(1)}\qquad (\e^{1/\varepsilon}<y\leqslant x).\eqdef{minZphi}$$
\par 
De plus
$$M(y;t)\ll{y\,\e^{Z(y;t)}\over \log 2y}\qquad (1\leqslant y\leqslant x).\eqdef{Myt}$$\par 
Une sommation d'Abel tenant compte de \eqref{minZphi} et \eqref{Myt} permet donc d'Žcrire
$$\eqalign{\gR_0&\ll{\e^{Z(x;s)}\over (\log x)^{2\gb}}\sum_{x^{\varepsilon_1}<d\leqslant x^{1-\varepsilon_1}}{t(d)\over d(\log x/d)^{1-2\gb}}\cr&\ll{\e^{Z(x;r)}\over (\log x)^{4\gb}}\int_{x^{\varepsilon_1}}^{x^{1-\varepsilon_1}}{\dd y\over y(\log y)^{1-2\gb}(\log x/y)^{1-2\gb}}\ll {\e^{Z(x;r)}\over \log x}\cdot\cr}$$
De mme, nous obtenons d'une part, gr‰ce ˆ \eqref{minZphi},
$$\eqalign{\gR_1&\ll{\e^{Z(x;s)}\over \log x}\int_1^{x^{\varepsilon_1}}{\e^{Z(y;t)}\over y\log (2y)}\d y\cr&\ll{\e^{Z(x;s)+Z(\exp(1/\varepsilon);t)}\log (1/\varepsilon)\over \log x}+{\e^{Z(x;s)}\over \log x}\int_{\exp(1/\varepsilon)}^{x^{\varepsilon_1}}{\e^{Z(y;t)}\over y\log y}\d y\cr
&\ll{\e^{Z(x;r)}\over \log x}\Big({1\over \varepsilon\log x}\Big)^{3\gb/2}\log (1/\varepsilon)+{\e^{Z(x;r)}\over (\log x)^{1+3\gb/2}}\int_{\exp(1/\varepsilon)}^{x^{\varepsilon_1}}{\dd y\over y(\log y)^{1-3\gb/2}}\cr
&\ll{\e^{Z(x;r)}\varepsilon^{3\gb/2}\log (1/\varepsilon)\over \log x}\ll {\e^{Z(x;f)}\varepsilon^{\gb/2}\over \log x},\cr}$$
o nous avons fait appel ˆ \eqref{pqr},
et, d'autre part,
$$\eqalign{\gR_2&\ll{\e^{ Z(x;s)}\over (\log x)^{2\gb}}\sum_{x^{1-\varepsilon_1}<d\leqslant x}{t(d)\over d(\log 2x/d)^{1-2\gb}}\cr&\ll{\e^{Z(x;r)}\over (\log x)^{1+2\gb}}\int_{x^{1-\varepsilon_1}}^x{\dd y\over y(\log 2x/y)^{1-2\gb}}\ll{x\e^{Z(x;f)}\varepsilon^{\gb/2}\over \log x}.\cr}$$
\par 
 Ces majorations sont bien compatibles avec l'estimation annoncŽe \eqref{Mf/Mr}.
\par 
Il reste ˆ Žvaluer le terme principal de \eqref{decMf}. \par 
Commenons par observer que la condition \eqref{cdpr} avec $\gh=1$ permet, sans altŽrer les termes d'erreur, de remplacer dans \eqref{decMf} le produit en $p$ par 
$$\prod_{p\leqslant x^\varepsilon}\sum_{\pnu\leqslant x/d}{\varphi(\pnu)\over \pnu}\prod_{p>x^\varepsilon}\sum_{\pnu\leqslant x/d}{s(\pnu)\over \pnu}\cdot\eqdef{prodphi/s}$$
De plus, lorsque $d\leqslant x^{1-\varepsilon_1}$, la seconde condition \eqref{C0} permet de supprimer la contrainte $\pnu\leqslant x/d$ dans le produit de gauche de \eqref{prodphi/s}. Par ailleurs, nous pouvons Žgalement, toujours sans dŽtŽrioration des termes rŽsiduels, substituer l'expression $\e^{Z(x/d;s)-Z(x^\varepsilon;s)}$ au produit de droite.
 Introduisant la notation
$$\gE(p;g):=\sum_{\nu\geqslant 0}{g(\pnu)\over \pnu}$$
pour le facteur eulŽrien d'indice $p$ d'une \fmu\ gŽnŽrique $g$, nous avons donc Žtabli~que
 $$\eqalign{M(x;f)&={\e^{-2\gamma\gb}x\over \Gamma(2\gb)}\prod_{p\leqslant x^\varepsilon}\gE(p;\varphi)\sum_{x^{\varepsilon_1}<d\leqslant x^{1-\varepsilon_1}}{\psi(d)\e^{ Z(x/d;s)-Z(x^\varepsilon;s)}\over d \log (x/d)}\cr&\quad+O\Big({\varepsilon^{\delta}x\e^{Z(x;r)-\gc Z(x;|f|-f)}\over \log x}\Big).\cr}\eqdef{dec2Mf}$$
\par 
Pour estimer la somme intŽrieure en $d$, nous observons que l'on a, pour une constante convenable~$c_0$, 
$$Z(x/d;s)=2\gb\log_2(x/d)+c_0\gb+O(\varepsilon).\eqdef{Zxd}$$
Par sommation d'Abel et en utilisant \eqref{Myt}, on vŽrifie aisŽment que la contribution du terme d'erreur est acceptable. \par  \goodbreak
DŽsignons ensuite par $a$ (resp. $b$) un entier gŽnŽrique dont  les facteurs premiers ne dŽpassent pas $x^\varepsilon$ (resp. excdent tous  $x^\varepsilon$), de sorte que chaque entier $d$ possde une dŽcomposition unique sous la forme $d=ab$. Soit alors $\varepsilon_2:=\varepsilon^{3/4}$. La majoration de Rankin
$$\sum_{a>x^{\varepsilon_2}}{t(a)\over a}\leqslant \e^{-\varepsilon_2/\varepsilon}\sum_{a\leqslant y}{t(a)\over a^{1-1/(\varepsilon\log x)}}\ll\e^{-\varepsilon^{-1/4}+Z(y;t)-Z(x^\varepsilon;t)}\eqdef{rankin}$$
permet de montrer que la contribution au terme principal de \eqref{dec2Mf} des entiers $d=ab$ tels que $a>x^{\varepsilon_2}$ est Žgalement nŽgligeable. Cela permet, toujours au prix d'une erreur acceptable, de remplacer $x/d$ par $x/b$ dans la sommation. On utilise ensuite la majoration
$$\sum_{b\leqslant v}\psi(b)\ll {v\e^{Z(v;t)-Z(x^\varepsilon;t)}\over \log v}\qquad (x^\varepsilon<v\leqslant y)$$ pour Žtablir que l'on peut remplacer le domaine de sommation par $x^{\varepsilon_1}<b\leqslant x^{1-\varepsilon_1}$, et enfin \eqref{rankin} pour Žtendre la sommation en $a$ ˆ toutes les valeurs possibles.  
\par 
Puisque $\varphi*\psi=f$, nous obtenons
$$\eqalign{M(x;f)&={\e^{-2\gamma\gb}x\over \Gamma(2\gb)}\prod_{p\leqslant x^\varepsilon}\gE(p;f)\sum_{x^{\varepsilon_1}<b\leqslant x^{1-\varepsilon_1}}{\psi(b)\e^{ Z(x/b;s)-Z(x^\varepsilon;s)}\over b \log (x/b)}\cr&\quad+O\Big({\varepsilon^{\delta}x\e^{Z(x;r)-\gc Z(x;|f|-f)}\over \log x}\Big).\cr}\eqdef{dec3Mf}$$
\par 
Prolongeons la fonction $\vartheta$ prŽcŽdemment dŽfinie sur les puissances de \nps\ en une fonction additive.
Nous avons ainsi $$\psi(n)=|\psi(n)|\e^{i\vartheta(n)+i\pi\kappa(n)}\qquad (n\geqslant 1)$$ 
o l'on a posŽ
$$\kappa(n):=\sum_{\di{p\|n}{|f(p)|<s(p)}}1.$$
 Notons encore 
$$\Theta(y):=\sum_{p>x^\varepsilon,\,\pnu\leqslant y}{t(\pnu)\{\vartheta(\pnu)+\pi\kappa(\pnu)\}\over \pnu}\cdot$$
Comme $\kappa(p)t(p)\leqslant r(p)-|f(p)|$, nous dŽduisons
 de la validitŽ de \eqref{cdpr} avec $\gh=1$ que
$$\Theta(y)\ll\varepsilon^{\delta}\qquad (x^{\varepsilon}<y\leqslant x).\eqdef{majThy}$$ 
Posons 
$$\Psi(y):=\sum_{b\leqslant y}{\psi(b)\over b}\cdot$$
Il suit, pour tout $y\in[x^{\varepsilon_1},x^{1-\varepsilon_1}]$,
$$\eqalign{
\e^{-i\Theta(y)}\Psi(y)&=\sum_{\di{b\leqslant y}{t(b)\neq0}}{t(b)\over b}\bigg\{1-\bigg(1-{|\psi(b)|\over t(b)}\bigg)\bigg\}\e^{i\vartheta(b)+i\pi\kappa(b)-i\Theta(y)}\cr
&=\sum_{b\leqslant y}{t(b)\over b}+O\big(\gR_3+\gR_4\big),\cr}$$
o l'on a posŽ
$$\gR_3:=\sum_{\di{b\leqslant y}{t(b)\neq0}}{t(b)\over b}\sum_{\pnu\|b}\Big(1-{|\psi(\pnu)|\over t(\pnu)}\Big),\quad\gR_4:=\sum_{b\leqslant y}{t(b)|\vartheta(b)+\pi\kappa(b)-\Theta(y)|\over b}\cdot$$
Nous estimons $\gR_3$ par interversion de sommation en notant que, pour tout $p$, $$t(p)-|\psi(p)|=r(p)-s(p)+||f(p)|-s(p)|\leqslant 2\{r(p)-|f(p)|\}.
$$  Nous obtenons $\gR_3\ll\varepsilon^{\delta_1}\e^{Z(y;t)-Z(x^\varepsilon;t)}$. Nous majorons ensuite $\gR_4$ en faisant appel au \ref{TKpond} aprs application de l'inŽgalitŽ de Cauchy--Schwarz. Au vu de \eqref{majThy}, nous obtenons ainsi
$$\Psi(y)=\sum_{b\leqslant y}{t(b)\over b}+O\Big(\varepsilon^{\delta}\e^{Z(y;t)-Z(x^\varepsilon;t)}\Big)\qquad (x^{\varepsilon_1}\leqslant y\leqslant x^{1-\varepsilon_1}).$$ 
Reportons ˆ prŽsent dans \eqref{dec3Mf} en effectuant une nouvelle sommation d'Abel pour substituer $t$ ˆ $\psi$ dans le terme principal. Compte tenu de \eqref{Zxd}, cette manipulation n'altre pas les termes d'erreur. Nous obtenons
$$\eqalign{M(x;f)&={\e^{-2\gamma\gb}x\over\Gamma(2\gb)}\prod_{p\leqslant x^\varepsilon}\gE(p;f)\sum_{x^{\varepsilon_1}<b \leqslant x^{1-\varepsilon_1}}{t(b)\e^{Z(x/b;s)-Z(x^\varepsilon;s)}\over b\log (x/b)}\cr&\quad+O\bigg({\varepsilon^{\delta}x\e^{ Z(x;r)-\gc Z(x;|f|-f)}\over \log x}\bigg).\cr}$$
\par 
En appliquant cette Žvaluation pour $f=r$, nous obtenons bien l'estimation annon\-cŽe~\eqref{Mf/Mr}.
\par\goodbreak
\medskip\bigskip
\paraunn{Preuve du \ref{comprf}}\drefun{comp}
Nous indiquons brivement ici comment modifier la dŽmonstration du \ref{fr} sous les hypothses indiquŽes dans l'ŽnoncŽ du \ref{comprf}.\par 
En application du \ref{cascv} avec l'hypothse \eqref{cdpr1}, la formule \eqref{decMf} demeure valable. Il en va de mme de la majoration des termes d'erreur $\gR_j$ $(0\leqslant j\leqslant 2)$.
\par 
La substitution indiquŽe en \eqref{prodphi/s} produit ˆ prŽsent des termes d'erreur 
$$\eqalign{&\ll \varepsilon^{\delta}x\prod_{p\leqslant x^\varepsilon}\gE(p;\varphi)\sum_{x^{\varepsilon_1}<d\leqslant x^{1-\varepsilon_1}}{t(d)\e^{Z(x/d;s)-Z(x^\varepsilon;s)}\over d\log (x/d)}\log \Big({\log x\over \log (x/d)}\Big)\cr&\ll \varepsilon^{\delta}x\,\e^{Z(x;\varphi)}\sum_{x^{\varepsilon_1}<d\leqslant x^{1-\varepsilon_1}}{t(d)\e^{Z(x;s)}\over d\{\log (x/d)\}^{1-2\gb}(\log x)^{2\gb}}\log \Big({\log x\over \log (x/d)}\Big)\cr
&\ll {\varepsilon^\delta x\,\e^{Z(x;r)-\gc Z(x;|f|-f)}\over \log x},
\cr}$$ 
ce qui est acceptable.
\par 
On vŽrifie Žgalement que la formule \eqref{dec3Mf} demeure valable.
\par 
Le point le plus dŽlicat consiste ˆ introduire la fonction $t$ dans la somme en $b$ de \eqref{dec3Mf}.
Ë cette fin, nous Žvaluons
$$\R(y):=\sum_{b\leqslant y}{\psi(b)\over b}-\sum_{b\leqslant y}{t(b)\over b}\qquad (x^{\varepsilon_1}<y\leqslant x^{1-\varepsilon_1})\eqdef{defgR}$$
avec pour objectif la majoration 
$$\R(y)\ll \varepsilon^\delta\sum_{b\leqslant y} {t(b)\over b}\cdot\eqdef{majgR}$$
Cela sera rendu possible gr‰ce ˆ l'hypothse additionnelle concernant les nombres $r(p)$, qui implique $t(p)\geqslant 2\gb$ pour tous les \nps\ considŽrŽs.
Nous pouvons clairement restreindre les deux sommations de \eqref{defgR} aux entiers sans facteur carrŽ avec une erreur acceptable. Comme l'hypothse \eqref{cdpr1} fournit la borne $\ll\varepsilon^\delta\log (1/\varepsilon)$ pour la contribution des \nps, nous pouvons Žgalement convenir que la sommation ne porte que sur des entiers $b$ composŽs.\par  Posons donc $w(n):=\{\psi(n)-t(n)\}\mu(n)^2$, de sorte que $w(p)\ll\{r(p)-\re f(p)\}^{w_f}$, et introduisons les fonctions sommatoires
$$\gW(y):=\sumast_{b\leqslant y}w(b),\qquad \gT(y):=\sumast_{b\leqslant y}t(b),$$
o l'astŽrisque indique que les entiers sommŽs ont au moins deux facteurs premiers.
Pour tout $z\in[x^\varepsilon,y]$, nous avons 
$$\eqalign{\sumast_{b\leqslant y}w(b)\log (y/b)&\ll \gT(z)\log y+\gT(y)\log (y/z)\cr&\ll{z\e^{Z(z;t)-Z(x^\varepsilon;t)}\log y\over \log z}+{y\e^{Z(y;t)-Z(x^\varepsilon;t)}\log (y/z)\over \log y}\cdot\cr}$$
En choisissant $z:=y/(\log x)^3$, nous obtenons
$$\sumast_{b\leqslant y}w(b)\log (y/b)\ll {y\e^{Z(y;t)-Z(x^\varepsilon;t)}\log_2y\over \log y}$$
et donc
$$\gW(y)={\gV(y)\over \log y}+O\bigg({y\e^{Z(y;t)-Z(x^\varepsilon;t)}\log_2y\over (\log y)^2}\bigg)\eqdef{estgM}$$
o l'on a posŽ
$$\eqalign{\gV(y):=\sumast_{b\leqslant y}w(b)\log b&\ll\sum_{x^\varepsilon<p\leqslant y/x^\varepsilon}{y|w(p)|\log p\over p\log (y/p)}\e^{Z(y/p;t)-Z(x^\varepsilon;t)}\cr
&\ll{y\e^{Z(y;t)-Z(x^\varepsilon;t)}\over \log y}\sum_{x^\varepsilon<p\leqslant y/x^\varepsilon}{|w(p)|\log p\over p}\bigg({\log y\over \log (y/p)}\bigg)^{1-2\gb}\cr
&\ll\varepsilon^\delta y\e^{Z(y;t)-Z(x^\varepsilon;t)},}$$
la dernire somme en $p$ ayant ŽtŽ estimŽe par l'inŽgalitŽ de Hšlder avec exposants $\gh/w_f$ et $\gh/(\gh-w_f)$, puis recours ˆ \eqref{cdpr1} et sommation d'Abel.
\par 
En reportant dans \eqref{estgM}, nous obtenons bien l'estimation annoncŽe \eqref{majgR}.
\par 
Cela permet d'opŽrer une sommation d'Abel dans \eqref{dec3Mf} pour obtenir
$$\eqalign{M(x;f)&={\e^{-2\gamma\gb}x\over \Gamma(2\gb)}\prod_{p\leqslant x^\varepsilon}\gE(p;f)\prod_{x^\varepsilon<p\leqslant x}\sum_{x^{\varepsilon_1}<b\leqslant x^{1-\varepsilon_1}}{t(b)\e^{ Z(x/b;s)-Z(x^\varepsilon;s)}\over b \log (x/b)}\cr&\quad+O\Big({\varepsilon^{\delta}x\e^{Z(x;r)-\gc Z(x;|f|-f)}\over \log x}\Big).\cr}\eqdef{dec4Mf}$$
En appliquant cette relation pour $f=r$, nous obtenons bien \eqref{Mf/Mr}, comme souhaitŽ.
\bigskip 
\paraunn{Preuves des corollaires}\drefun{demcor}
\paradeuxn{Preuve du \ref{majm}} Observons d'emblŽe que l'insertion de l'hypothse \eqref{rp/p} dans la dŽmonstration  du thŽo\-rme~1.1 de \citer{Te16}\note{Plus prŽcisŽment, la conclusion du thŽorme 1.1 de \citer{Te16} est inchangŽe si l'on remplace la condition (v)  par $\sum_{y<p\leqslant (1+\lambda)y}f(p)\log p\gg y$ $(x^{\tau} <y\leqslant x^{1-\sigma})$, 	avec les notations de \citer{Te16}.} fournit sans autre modification l'Žvaluation
$$M(x;r)\asymp {x\over \log x}\e^{Z(x;r)}.\eqdef{Mxrstand}$$
Le mme rŽsultat dŽcoule du thŽorme 2 de \citer{El16}.
\par 
Comme au corollaire III.4.12 de \citer{Te15}, l'hypothse $r(p)\ll1$ implique
$$Z(y;r)-m_f(y;T)\leqslant Z(x;r)-m_f(x;T)+O(1)\qquad (2\leqslant y\leqslant x).$$
Il s'ensuit que, pour $s=1+\alpha+i\tau$, $1/\log x\leqslant \alpha\leqslant 1$, $|\tau|\leqslant T$, nous pouvons Žcrire, avec la notation~\eqref{MphiH},
$$\re v_f(s)\leqslant Z(\e^{1/\alpha};r)-m_f\big(\e^{1/\alpha};T\big)+O(1)\leqslant Z(x;r)-m_f(x;T)+O(1).\eqdef{revf}$$
Nous utilisons cette majoration pour $\alpha\leqslant \alpha_0:=\e^{m_f(x;T)/\gb}/\log x$ et employons la majoration triviale $$\re v_f(s)\leqslant Z(\e^{1/\alpha};r)+O(1)\leqslant Z(x;r)-\gb\log (\alpha\log x)+O(1)$$ lorsque $\alpha_0<\alpha\leqslant 1$. Il suit
$$H_T(\alpha)\ll\normalbaselineskip=15pt\cases{\e^{Z(x;r)-m_f(x;T)} &$(\alpha\leqslant \alpha_0)$\cr\dsp{\e^{Z(x;r)}\over (\alpha\log x)^\gb}&$(\alpha_0<\alpha\leqslant 1$).\cr}$$
L'intŽgrale de \eqref{majeff} est donc
$$\ll{\{1+m_f(x;T)\}\e^{Z(x;r)}\over \e^{m_f(x;T)}}\cdot$$
Cela implique bien la majoration annoncŽe en reportant dans \eqref{majeff} compte tenu de \eqref{Mxrstand}.
\smallskip\medskip\goodbreak
\paradeuxn{Preuve du \ref{OmE}}\drefdeux{demOmE}
ConsidŽrons la \fmu\ $f(n):=z^{\Omega(n;E)}$. Commenons par une minoration de la quantitŽ $m_f(x;T)$ dŽfinie en \eqref{defmfyT} avec $T:=\log x$. 
  Nous avons, pour une certaine valeur de $\tau\in[0,T]$, 
$$m_f(x;T)\geqslant r_1\sum_{p\leqslant x}{1-\cos(\vartheta_p-\tau\log p)\over p},$$
o l'on a posŽ $r_1:=\min(r,1)$, $\vartheta_p:=\1_E(p)\vartheta$. Notant $h_\vartheta$ la fonction indicatrice de $$[-\pi,\pi]\sset([-\ft13\pi,\ft13\pi]\cup[\vartheta-\ft13\pi,\vartheta+\ft13\pi])\ (\mod 2\pi),$$ nous pouvons donc Žcrire
$$m_f(x;T)\geqslant \dm r_1\sum_{p\leqslant x}{h_\vartheta(\tau\log p)\over p}\cdot$$
Le membre de droite relve du \tnp\ via une sommation d'Abel dont les dŽtails sont explicitŽs par exemple au lemme III.4.13 de \citer{Te15}. Pour tout $w\in[2,x]$, nous avons
$$\sum_{w<p\leqslant x}{h_\vartheta(\tau\log p)\over p}\geqslant \ft16\log \Big({\log x\over \log w}\Big)+O\bigg({1\over |\tau|\log w}+{1+|\tau|\over \e^{\sqrt{\log w}}}\bigg).\eqdef{somhth}$$
Si $1\leqslant |\tau|\leqslant T:=\log x$, nous choisissons $w:=\exp\{(\log_2x)^2\}$ et obtenons pour $x$ assez grand
$$m_f(x;T)\geqslant \ft1{12}r_1\log_2x+O(r_1\log_3x).$$
DŽfinissons ensuite $v$ par $\log v=\e^{-\fs13(1-\cos\vartheta)E(x)}\log x$.  Si $1/\log v< |\tau|\leqslant 1$, nous choisissons $w=v$ dans \eqref{somhth} et obtenons
$$\sum_{v<p\leqslant x}{h_\vartheta(\tau\log p)\over p}\geqslant\ft1{18}(1-\cos\vartheta)E(x)+O(1)\geqslant {\vartheta^2\over 9\pi^2}E(x)+O(1).$$
Enfin, si $|\tau|\leqslant 1/\log v$, nous avons trivialement
$$\eqalign{\sum_{p\leqslant v}{1-\cos(\vartheta_p-\tau\log p)\over p}&\geqslant\sum_{\di{p\leqslant v}{p\in E}}{1-\cos\vartheta+O(\tau\log p)\over p}\cr&\geqslant (1-\cos\vartheta)E(x)+O(1)-\sum_{\di{v<p\leqslant x}{p\in E}}{1-\cos\vartheta\over p}\cr&\geqslant \ft13(1-\cos\vartheta)E(x)+O(1)\geqslant {2\vartheta^2\over 3\pi^2}E(x)+O(1).\cr}$$
En rassemblant nos rŽsultats, nous obtenons
$$m_f(x;T)\geqslant \min\Big(\ft1{12}r_1\log_2x+O(r_1\log_3x),{\vartheta^2r_1\over 18\pi^2} E(x)\Big)+O(1)$$
et finalement 
$$m_f(x;T)\geqslant {\vartheta^2r_1\over 18\pi^2} E(x)+O(1).\eqdef{minmf}$$
\par 
En reportant dans \eqref{majfm}, nous obtenons bien 
\eqref{majSxzE}. Cela Žtablit l'assertion (i).
\par \goodbreak
Prouvons l'assertion (ii). Soit $K$ une constante positive.  Nous appliquons le \ref{comprf} sous la forme de la validitŽ de \eqref{Mf/Mr} pour la fonction $f(n):=z^{\Omega(n;E)}$ lorsque $|\vartheta|\leqslant K\gt(x;E)$, avec $r(n):=r^{\Omega(n;E)}$, $A:=1$, $\delta_1:=\eta\delta_0(\gb)$, $\delta:=\dm\delta_1$, $\gb:=\ft18\kappa$, $\varepsilon:=|\vartheta|^{1/\delta}+1/\sqrt{\log x}$, o $\eta=\eta(\kappa,K)$ est une constante assez petite. Comme
$$r(p)-\re f(p)\leqslant r(1-\cos\vartheta)\ll \vartheta^2\ll\varepsilon^{2\delta},$$
 la condition \eqref{cdpr2} est bien satisfaite. 
\par
\goodbreak Nous avons de mme, toujours pour \hbox{$|\vartheta|\leqslant K\gt(x;E)$},
$$\eqalign{\sum_{p\leqslant x}{r(p)-\re f(p)\over p}&=(1-\cos \vartheta)E(x)\leqslant \dm \vartheta^2E(x)\leqslant \dm K^2\log E(x)\ll \eta K^2\log (1/\varepsilon).\cr}$$
 La condition \eqref{pqr} est donc Žgalement satisfaite pour un choix convenable de $\eta$. La formule~\eqref{Mf/Mr} fournit alors l'Žvaluation
$$S(x;z,E)=S(x;r,E)\prod_{p\in E}{1-r/p\over 1-z/p}+O\Big(x\,\e^{(r-1)E(x)-\gc(r-z)E(x)}\varepsilon^{\delta}\Big).\eqdef{Sz/Sr}$$  Cela implique bien \eqref{SzE/SrE} puisque  $S(x;r,E)\asymp x\,\e^{(r-1)E(x)}$. Comme indiquŽ dans l'intro\-duction, nous obtenons \eqref{NmxE} par intŽgration sur le cercle $|z|=r=m/E(x)$.
\par \goodbreak
Pour achever la dŽmonstration du point (iii), il reste ˆ Žtablir \eqref{majunivNm}. Nous pouvons supposer $r\leqslant \kappa\leqslant \dm$. Nous utilisons la majoration de Montgomery \& Vaughan \citer{MV01}
$$\sum_{n\leqslant x}{f(n)\over n}\ll{1\over \log x}\int_{1/\log x}^1 {H_1(\alpha)\over \alpha}\d\alpha\eqdef{MV1}$$
o l'on a posŽ
$$H_1(\alpha)^2:=\sum_{k\in\z}\max_{\di{1+\alpha\leqslant \sigma\leqslant 2}{|\tau-k|\leqslant 1/2}}\abs{\e^{v_f(s)}\over s-1}^2.$$
D'aprs \eqref{revf} et \eqref{minmf}, nous avons $$\max_{|\tau|\leqslant \log x}\re v_f(s)\leqslant Z(x;r)-\vartheta^2rE(x)/180+O(1),$$
donc
$$H_1(\alpha)\ll {\e^{Z(x;r)-r\vartheta^2E(x)/180}\over \alpha}+{\e^{Z(x;r)}\over \sqrt{\log x}}\qquad (1/\log x\leqslant \alpha\leqslant  1).$$\goodbreak
Il rŽsulte alors de \eqref{MV1} que
$$\sum_{n\leqslant x}{z^{\Omega(n;E)}\over n}\ll\e^{Z(x;r)-r\vartheta^2E(x)/180}+{\e^{Z(x;r)}\over \log x}.$$
Comme $Z(x;r)=(r-1)E(x)+\log_2x+O(1)$, nous obtenons, par intŽgration sur le cercle $|z|=r=(m+1)/E(x)\leqslant \kappa+1/E(x)$,
$$\sum_{\di{n\leqslant x}{\Omega(n;E)=m}}{1\over n}\ll\e^{-E(x)}{E(x)^m\over m!}\log x.\eqdef{majNmpond}$$
\par 
La suite de la preuve est essentiellement identique ˆ celle de \Halasz\ \citer{Ha72}. Nous reproduisons les dŽtails pour la commoditŽ du lecteur.
\par \goodbreak
Nous avons  $$\eqalign{N_m(x;E)&={1\over \log x}\sum_{\di{n\leqslant x}{\Omega(n;E)=m}}\log n+O\bigg({x\over \log x}\sum_{\di{n\leqslant x}{\Omega(n;E)=m}}{1\over n}\bigg)\cr&={1\over \log x}\sum_{\di{dp\leqslant x}{\Omega(dp;E)=m}}\log p+O\Big(x\e^{-E(x)}{E(x)^m\over m!}\Big)\cr&\ll{1\over \log x}\sum_{\di{d\leqslant x}{\Omega(d;E)=m\,{\rm{ ou }}\,m-1}}{x\over d}+x\e^{-E(x)}{E(x)^m\over m!}\ll x\e^{-E(x)}{E(x)^m\over m!}\cr}$$
o nous avons appliquŽ \eqref{majNmpond} pour $m$ et $m-1$, ce dernier terme n'Žtant prŽsent que si $m\geqslant 1$. 
\bigskip
\paradeuxn{Preuve du \ref{repad}} Nous pouvons, sans restreindre la gŽnŽralitŽ, supposer que $\mu_x$ est arbitrairement petit.
Avec la perspective de recourir ˆ l'inŽgalitŽ de Berry-Esseen, nous Žvaluons la fonction caractŽristique
$$\varphi(t):={1\over M(x;r)}\sum_{n\leqslant x}r(n)\e^{it\{h(n)-E\}/D}\qquad (t\in\r),$$
o, pour simplifier l'Žcriture, nous avons posŽ $E:=E_h(x;r)$, $D:=D_h(x;r)$.
\par 
Appliquons d'abord le \ref{comprf} ˆ la fonction $f(n):=r(n)\e^{ith(n)/D}$. Nous avons $$\leqalignno{\max_{p\leqslant x}\{r(p)-\re f(p)\}&\ll t^2\mu_x^2\qquad (|t|\leqslant 1/\mu_x),&\eqdef{r-reith}\cr
\sum_{p\leqslant x}{r(p)-\re f(p)\over p}&\leqslant \dm t^2\qquad (t\in\r).&\eqdef{reith/p}\cr}$$
Soit alors $\eta$ une constante absolue suffisamment petite. Il rŽsulte de ce qui prŽcde que les conditions \eqref{pqr} et \eqref{cdpr1} sont satisfaites pour le choix $\delta:=\eta\delta_0(\gb)$ et $\varepsilon^{\delta}=|t|\mu_x$ si $\gb$ est assez petit et
$${1\over (\log x)^{\eta\delta_0(\gb)/2}}< |t|\mu_x\leqslant \mu_x\sqrt{{2\over\eta}\log {2\over \mu_x}}.\eqdef{domt}$$
Nous pouvons donc utiliser l'estimation \eqref{Mf/Mr}  dans ce domaine.

La contribution des puissances $\pnu$ avec $\nu\geqslant 2$ au produit apparaissant au membre de droite de \eqref{Mf/Mr} Žtant trivialement $1+O\big(t/D\big)$, nous obtenons, dans le domaine \eqref{domt}, en dŽveloppant $\e^{ith(p)}$ ˆ l'ordre 3,
$$\varphi(t)=\e^{-t^2/2}\Big\{1+O\Big(t^3\mu_x+{t\over D}\Big)\Big\}+O\Big(\e^{-2\gc t^2/\pi^2}t\mu_x\Big).$$ \par 
Lorsque $|t|\mu_x\leqslant (\log x)^{-\eta\delta_0(\gb)/2}$, observant que $E\ll D\mu_x\log_2x$, nous employons la majoration triviale
$$\eqalign{\varphi(t)-\e^{-t^2/2}&\ll t^2+{|t|\over DM(x;r)}\sum_{n\leqslant x}r(n) \{|h(n)|+|E|\}\cr&\ll t^2+|t|\mu_x\log_2x+{|t|\over D}\sum_{\pnu\leqslant x}{r(\pnu)|h(\pnu)|M(x/\pnu;r)\over M(x;r)}\cr&\ll t^2+|t|\mu_x\log_2x+{|t|\over D}\sum_{\pnu\leqslant x}{r(\pnu)|h(\pnu)|\over \pnu}\Big({\log x\over \log (2x/\pnu)}\Big)^{1-\gc},\cr}$$
o $\gc$ est une constante positive convenable: cela rŽsulte  du \theoreme 1.1 de~\citer{Te16} sous la forme $M(x;r)\asymp x\e^{Z(x;r)}/\log x$, de la majoration de Halberstam--Richert~\citer{HR79} pour $M(x/\pnu;r)$, et de \eqref{L2} avec $\varepsilon:=(\log x)^{-1/4}$. Compte tenu des hypothses (ii) et (iv), nous obtenons
$\varphi(t)-\e^{-t^2/2}\ll t^2+|t|\mu_x\log_2x.$ 
\par\smallskip\goodbreak
Lorsque $\sqrt{(2/\eta)\log (2/\mu_x)}\leqslant |t|\leqslant 1/\mu_x$, nous majorons $\varphi(t)$ gr‰ce au  \ref{majm}. Ë cette fin observons d'abord que $$\mu_x\gg1/\sqrt{\log_2x}\eqdef{muxrtlogx}$$ en vertu de l'inŽgalitŽ
$$D^2=\sum_{p\leqslant x}{r(p)h(p)^2\over p}\leqslant \mu_x^2D^2\sum_{p\leqslant x}{r(p)\over p},$$
qui rŽsulte de l'hypothse (iii). Donnons-nous alors une constante absolue $\lambda_0>0$ assez petite pour que $\tau_0:=\e^{\lambda_0/\mu_x^2}/\log x<1/\sqrt{\log x}$. Nous avons 
$$\sum_{p\leqslant x}{r(p)-\re\{f(p)/p^{i\tau}\}\over p}=\sum_{p\leqslant x}{r(p)\{1-\cos(th(p)/D-\tau\log p)\}\over p}\gg{1\over \mu_x^2}$$
uniformŽment pour $\tau_0<|\tau|\leqslant  T_x:=\exp\{(\log x)^{1/4}\}$: cela rŽsulte du \tnp\ (\cf, par exemple, le lemme III.4.13 de \citer{Te15}) en restreignant par exemple la sommation aux \nps\ $p>\exp(1/\tau_0)$ tels que $\tau\log p\in[\dm\pi-\ft1{10},\dm\pi+\ft1{10}]\,(\mod2\pi)$. De plus, lorsque $|\tau|\leqslant \tau_0$,
$$\eqalign{\sum_{p\leqslant x}{r(p)-\re\{f(p)/p^{i\tau}\}\over p}&\geqslant \sum_{p\leqslant \exp(1/\tau_0)}{r(p)(1-\cos\{th(p)/D\})\over p}+O(1)\cr&\geqslant{2t^2\over \pi^2D^2} \sum_{p\leqslant \exp(1/\tau_0)}{r(p)h(p)^2\over p}+O(1)\cr&\geqslant{2t^2\over \pi^2}+ O\Big(1+{t^2\over D^2}\sum_{\exp(1/\tau_0)<p\leqslant x}{\mu_x^2D^2\over p}\Big)\geqslant \ft1{10}t^2\cr}$$
pour un choix convenable de $\lambda_0$. Dans le domaine considŽrŽ, nous avons donc, avec la notation~\eqref{defmfyT},
$$m_f(x;T_x)\geqslant {1\over 5\eta}\log \Big({2\over \mu_x}\Big).$$ Compte tenu de \eqref{muxrtlogx}, nous pouvons donc fixer $\eta$ de sorte que 
$$\varphi(t)\ll\mu_x^2\qquad \big(\sqrt{(2/\eta)\log (2/\mu_x)}\leqslant |t|\leqslant 1/\mu_x\big).$$
\par 
L'inŽgalitŽ de Berry--Esseen (\cf, \eg, \citer{Te15}, \theoreme II.7.16) sur l'intervalle $-T\leqslant t\leqslant T$ avec $T:=1/\mu_x$ fournit donc
$$\eqalign{F_x(E+zD;f,r)-&{1\over \sqrt{2\pi}}\int_{-\infty}^z\e^{-u^2/2}\d u\ll{1\over T}+\int_{-T}^T\abs{{\varphi(t)-\e^{-t^2/2}\over t}}\d t\cr&\ll\mu_x+{\log_2x\over (\log x)^{\eta\delta_0(\gb)/2}}+{1\over D}\cdot\cr}\eqdef{Fx-Phi}$$
Compte tenu de \eqref{muxrtlogx}, on voit que le second terme du membre de droite de \eqref{Fx-Phi} est dominŽ par le premier. \par 
Lorsque l'hypothse (i) est remplacŽe par la condition en moyenne \eqref{moyrp}, nous devons choisir $\varepsilon^{2\delta}:=t^2\mu_x^2\log (1+1/|t|\mu_x)$ pour que la condition \eqref{cdpr} avec $\gh=1$ soit satisfaite et pouvoir ainsi appliquer le \ref{fr}. Les calculs prŽcŽdents demeurent alors valides, {\it mutatis mutandis} et fournissent le rŽsultat indiquŽ.
\par 
Cela achve la dŽmonstration.
\medskip\bigskip\goodbreak
\paradeuxn{Preuve du \ref{Omphi}}
La fonction $h(n):=\Omega(\varphi(n))$ est additive. Le thŽorme de Bombieri--Vinogradov et les estimations de moments centrŽs pour la fonction $p\mapsto\Omega(p-1)$, telles qu'Žtablies par Alladi~\citer{Al85} dans un ŽnoncŽ trs gŽnŽral, permettent aisŽment de montrer que, pour $x\geqslant 3$, nous avons
$$\eqalign{E_h(x;r)&:=\sum_{p\leqslant x}{r(p)\Omega(p-1)\over p}=\dm\varrho(\log_2x)^2+O(\log_2 x),\cr 
D_h(x;r)^2&:=\sum_{p\leqslant x}{r(p)\Omega(p-1)^2\over p}=\ft13\varrho(\log_2x)^3+O\big((\log_2x)^2\big).\cr}$$ Cependant, le \ref{repad} n'est pas applicable directement car la condition (iii) n'est pas satisfaite. Cette difficultŽ peut tre contournŽe en introduisant la fonction additive $$h_1(n):=\sum_{\di{\pnu\|n}{\Omega(p-1)\leqslant 5\log_2x}}\Omega\big(\varphi(\pnu)\big).$$
Nous verrons que $h(n)=h_1(n)$ sauf peut-tre pour au plus $\ll x/\sqrt{\log x}$ entiers $n\leqslant x$. Admettons cela un instant. Au vu du terme d'erreur de \eqref{repOmphi}, il suffit donc d'Žvaluer la fonction de rŽpartition relative ˆ~$h_1$, qui relve du \ref{repad} avec $\mu_x\asymp 1/\sqrt{\log_2x}$, puisque la condition (iv) de cet ŽnoncŽ rŽsulte trivialement de l'hypothse \eqref{C1}. Cela implique bien le rŽsultat annoncŽ.
\par \goodbreak
Pour estimer le cardinal $ U(x)$ de l'ensemble des entiers $n\leqslant x$ tels que $h_1(n)<h(n)$, nous majorons sa fonction sommatoire par
$$\sum_{p|n}(\ft32)^{\Omega(p-1)-5\log_2x}.$$
Il suit
$$\eqalign{U(x)&\leqslant {x\over (\log x)^{5\log (3/2)}}\sum_{p\leqslant x}{(3/2)^{\Omega(p-1)}\over p}\leqslant{x\over (\log x)^{5\log (3/2)}}\sum_{n\leqslant x}{(3/2)^{\Omega(n)}\over n}\cr&\ll x(\log x)^{3/2-5\log (3/2)}\ll {x\over \sqrt{\log x}}\cdot\cr}$$ 
\bigskip
\noi{\bf Remerciements.} L'auteur tient ˆ exprimer sa gratitude ˆ RŽgis de la Bretche pour son aide lors de la prŽparation de cet article.
\bigskip\bigskip\goodbreak
\centerline{\twelvebf Bibliographie}
\bigskip\bigskip\bigskip
\eightpoint{
\parskip=3pt
\bibtem{Al85} K. Alladi, 
Moments of additive functions and the sequence of shifted primes,
{\it Pacific J. Math. \bf118}, \numero2 (1985), 261--275.\par  
\bibtem{Ba89} M. Balazard,
Remarques sur un thŽorme de G. Hal‡sz et A. S‡rkšzy, {\it 
Bull. Soc. Math. France \bf117}, \numero4 (1989), 389--413.\par 
\bibtem{BS89} M. Balazard et A. Smati, Travaux de Pomerance sur la fonction $\varphi$
d'Euler, {\it Publ. Math. d'Orsay} {\bf 87}, (1989/1990), 5--36. \par 
\bibtem{BGS13} A. Balog, A. Granville, \& K. Soundararajan, 
Multiplicative functions in arithmetic progressions, {\it
Ann. Math. QuŽ. \bf37}, \numero1 (2013), 3--30. \par 
\bibtem{BS96} A. Bir\'o \& T. Szamuely, A Tur\'an--Kubilius
inequality with multiplicative weights,
        {\it Acta Math. Hung. \bf 70} (1996), \numeros 1-2, 39--56.\par 
\bibtem{BT12} R. de la Bretche \& G. Tenenbaum, Sur la concentration de certaines fonctions additives, {\it Math. Proc. Camb. Phil. Soc. \bf 152}, \numero1 (2012), 179--189; erratum {\it ibid.}, 191.\par
\bibtem{El80} P.D.T.A. Elliott, {\it Probabilistic number theory : central limit theorems}, Grundlehren der
Math. Wiss.  240, Springer-Verlag, New York, Berlin, Heidelberg, 1980.\par 
\bibtem{El16} P.D.T.A. Elliott, Multiplicative functions mean values : asymptotic estimates,  {\it Functiones et Approximatio, \bf56}, \numero2 (2017), 217-238.\par 
\bibtem{EP85} P. Erd\H os \& C. Pomerance, On the normal number of prime factors of $\varphi(n)$,
{\it Rocky Mountain Journal of Math. (2)} {\bf 15} (1985), 343--352. \par 
\bibtem{Ga70} P.X. Gallagher, A large sieve density estimate near $\sigma =1$, 
{\it Invent. Math. \bf11} (1970), 329--339. \par 
\bibtem{GS03} A. Granville \& K. Soundararajan,
Decay of mean values of multiplicative functions,
{\it Canad. J. Math. \bf55}, \numero6 (2003), 1191--1230.\par  
\bibtem{GS08} A. Granville \& K. Soundararajan,
Pretentious multiplicative functions and an inequality for the zeta-function, in: {\it Anatomy of integers}, 191--197,
CRM Proc. Lecture Notes 46, Amer. Math. Soc., Providence, RI, 2008. \par 
\bibtem{Ha68} G. Hal‡sz, †ber die Mittelwerte multiplikativer zahlentheoretischer Funktionen,
{\it Acad. Math. Acad. Sci. Hungar. \bf19} (1968), 365--403.\par 
\bibtem{Ha71} G. Hal‡sz, On the distribution of additive and the mean values of multiplicative
arithmetic functions, {\it Stud. Sci. Math. Hungar. \bf6} (1971), 211--233.\par 
\bibtem{Ha72} G. Hal‡sz, Remarks to my paper ``On the distribution of additive and the mean values of multiplicative
arithmetic functions", {\it Acta Math. Acad. Scient. Hungar. \bf23} (1972), 425--432.\par
\bibtem{Ha75} G. Hal‡sz, On the distribution of additive functions, \AA\ {\bf27} (1975), 143--152.\par 
\bibtem{HR79} H. Halberstam \& H.-E. Richert, On a result of R.R. Hall, {\it J. Number Theory \rm (1) \bf 11}
(1979), 76--89.\par 
\bibtem{Ha95} R.R. Hall, A sharp inequality of Hal‡sz type for the mean value of a multiplicative arithmetic function, {\it  Mathematika \bf42}, \numero1 (1995), 144--157.\par 
\bibtem{HT91} R.R. Hall \& G. Tenenbaum, Effective mean value estimates for complex multiplicative
functions, {\it Math. Proc. Camb. Phil. Soc.} {\bf 110} (1991),
337--351. 
\bibtem{In09} K.-H. Indlekofer, On a quantitative form of Wirsing's mean-value theorem for multiplicative functions, {\it Publ. Math. Debrecen \bf 75}, \numeros 1-2 (2009), 105-121.\par 
\bibtem{IKW01} K.-H. Indlekofer, I. K‡tai, \& R. Wagner, A comparative result for multiplicative functions, {\it Lith. Math. J. \bf 41}, \numero2 (2001), 143--157.\par 
\bibtem{Ka30} J. Karamata, Sur un mode de croissance des fonctions rŽgulires, {\it Mathematica
(Cluj) \bf 4} (1930), 38--53.\par 
\bibtem{Ko13}  D. Koukoulopoulos, Pretentious multiplicative functions and the prime number theorem for arithmetic progressions, {\it Compos. Math. \bf149} \numero7 (2013), 1129--1149.\par 
\bibtem{LT82} B.V. Levin \& N.M. Timofeev, 
Un thŽorme de comparaison pour les fonctions multiplicatives (en russe),
{\it Acta Arith. \bf42} \numero 1 (1982/83), 21--47; Corrigendum {\it ibid.  \bf42} \numero3 (1983), 325.\par
\bibtem{Ma17a} A.P. Mangerel, A strengthening of theorems of Hal‡sz and Wirsing, prŽpubl.,  	arXiv:1604.00295v3.\par 
\bibtem{Ma17b} A.P. Mangerel, On the number of restricted prime factors of an integer, II, prŽpubl.,	arXiv:1604.01658v2.\par 
\bibtem{MJT98} F.
Marie-Jeanne \& G. Tenenbaum, Une propriŽtŽ de concentration liŽe ˆ la fonction d'Euler, {\it Indag. Math.,} N.S., \bf9 \rm \numero 3 (1998), 382--403. \par 
\bibtem{Mo78} H.L. Montgomery, A note on mean values of multiplicative functions, Report
\numero 17, Institut Mittag-Leffler, Djursholm, 1978, 9 pp.\par 
\bibtem{MV74}  H.L. Montgomery \& R.C. Vaughan,
Hilbert's inequality,
{\it J. London Math. Soc. \rm (2) \bf8} (1974), 73-82. \par
\bibtem{MV01} H.L. Montgomery \& R.C. Vaughan, Mean values of multiplicative functions, {\it Per. Math. Hung. \bf43}, 1-2 (2001), 199--204.\par 
\bibtem{Sa77} A. \Sarkozy, Remarks on a paper of G. Hal‡sz, {\it Period. Math. Hungar. \bf8}, \numero2 (1977), 135--150.\par 
\bibtem{Sh80}  P. Shiu, A Brun--Titchmarsh theorem for multiplicative functions, {\it J. reine
angew. Math. \bf313} (1980), 161--170.\par\goodbreak 
\bibtem{Te15} G. Tenenbaum, {\it Introduction ˆ la thŽorie analytique et probabiliste des nombres},  4me Žd., coll. ƒchelles, Belin, Paris 2015, 592 pp.; English translation: Graduate Studies in Mathematics 163, Amer. Math. Soc. 2015.\par 
\bibtem{Te16} G. Tenenbaum, Fonctions multiplicatives, sommes d'exponentielles, et loi des grands nombres, {\it Indag. Math. \bf27} (2016), 590--600.
\bibtem{Wi67} E. Wirsing, Das asymptotische Verhalten von Summen \"uber multiplikative Funktionen
II, {\it Acta Math. Acad. Sci. Hung. \bf 18} (1967), 411--467.

\smallskip}
\vskip 2mm
{\sevenrm\baselineskip9pt
\obeylines
\bigskip\medskip
GŽrald Tenenbaum
 Institut ƒlie Cartan
 UniversitŽ de Lorraine
 BP 70239
 54506 Vand{\oe}uvre-ls-Nancy Cedex
 France}

\end